\documentclass[10pt]{article}

\usepackage{indentfirst}
\usepackage{amsfonts}
\usepackage{amsmath}
\usepackage{amssymb}
\usepackage{amsbsy}
\usepackage{epic}

\newtheorem{dfn}{Definition} [section]
\newtheorem{theorem}[dfn]{Theorem}
\newtheorem{lemma}[dfn]{Lemma}

\newtheorem{conj}[dfn]{Conjucture}
\newenvironment{pf}{\noindent{\bf Proof.}}
{\enspace\vrule height5pt depth0pt width5pt}
\def\deg{{\rm deg}}
\def\L{{\mathcal L}}
\def\J{{\mathcal J}}
\def\F{{\mathcal F}}

\begin{document}

\title{An upper bound on the fractional chromatic number of triangle-free subcubic graphs}

\author{Chun-Hung Liu\thanks{E-mail:cliu87@math.gatech.edu.} \\
{\small School of Mathematics, Georgia Institute of Technology, Atlanta, Georgia 30332, USA}}

\date{November 20, 2012 (Revised at May 15, 2014)}

\maketitle

\abstract{
An $(a:b)$-coloring of a graph $G$ is a function $f$ which maps the vertices of $G$ into $b$-element subsets of some set of size $a$ in such a way that $f(u)$ is disjoint from $f(v)$ for every two adjacent vertices $u$ and $v$ in $G$.
The fractional chromatic number $\chi_f(G)$ is the infimum of $a/b$ over all pairs of positive integers $a,b$ such that $G$ has an $(a:b)$-coloring.
Heckman and Thomas conjectured that the fractional chromatic number of every triangle-free graph $G$ of maximum degree at most three is at most $2.8$.
Hatami and Zhu proved that $\chi_f(G) \leq 3-3/64 \approx 2.953$.
Lu and Peng improved the bound to $\chi_f(G) \leq 3-3/43 \approx 2.930$.
Recently, Ferguson, Kaiser and Kr\'{a}l' proved that $\chi_f(G) \leq
32/11 \approx 2.909$.
In this paper, we prove that $\chi_f(G) \leq 43/15 \approx 2.867$.
}

\section{Introduction}
A {\it proper $k$-coloring} of a graph $G$ is an assignment of one of $k$ colors to each vertex $v$ of $G$ such that adjacent vertices receive different colors.
The {\it chromatic number} $\chi(G)$ of a graph $G$ is the minimum $k$ such that $G$ has a proper $k$-coloring.
Graph coloring is one of the most celebrated topics in graph theory.
It has been widely explored and has many generalizations.
For any positive integer $a$, let $[a]$ be the set $\{1,2,...,a\}$.
Given a set $S$, we define $2^{S}$ to be the collection of subsets of $S$.
An {\it $(a:b)$-coloring} of a graph $G$ is a function $f: V(G) \rightarrow 2^{[a]}$ such that $\lvert f(v) \rvert = b$ for every vertex $v$, and $f(x) \cap f(y) = \emptyset$ for every pair of adjacent vertices $x$ and $y$.
The {\it fractional chromatic number} $\chi_f(G)$ of $G$ is the infimum of $a/b$ over all pairs of positive integers $a,b$ such that $G$ has an $(a:b)$-coloring.
In particular, every $(k:1)$-coloring is a proper $k$-coloring, so $\chi_f(G) \leq \chi(G)$.

Fractional coloring can be investigated from the point of view of optimization.
An {\it independent set} $I$ of a graph is a subset of vertices such that every pair of vertices in $I$ are non-adjacent.
Observe that the chromatic number is the minimum number of nonempty independent sets that partition the set of vertices.
In other words, the chromatic number is the optimal value of a certain integer programming problem.
In fact, the optimal value of the LP-relaxation of this integer programming problem is the fractional chromatic number, and the infimum in the definition of the fractional chromatic number is attained \cite{su}.
Thus, the infimum can be replaced by minimum.
Furthermore, by taking advantage of the duality of linear programming problems, the fractional chromatic number of a graph $G$ is at most $k$ if and only if for every weighted function defined on vertices of $G$, there exists an independent set $I$ such that the sum of weights of the vertices in $I$ is at least $w/k$, where $w$ is the sum of weights of all vertices in $G$.

Brooks' Theorem implies that $\chi_f(G) \leq \Delta(G)$, where $\Delta(G)$ is the maximum degree of the graph $G$, unless $G$ is the complete graph or an odd cycle.
On the other hand, it is easy to show that $\chi_f(G) \geq \omega(G)$, where $\omega(G)$ is the {\it clique number} which is the maximum size of a subgraph whose vertices are pairwise adjacent in $G$.
So one might expect a better upper bound when $\omega(G)<\Delta(G)$.
Indeed, it is not hard to see that $\chi_f(G) \geq \lvert V(G) \rvert/\alpha(G)$, where $\alpha(G)$ is the maximum size of an independent set in $G$.
So the best upper bound we can expect is $\lvert V(G) \rvert/ \alpha(G)$.

A graph is {\it subcubic} if it has the maximum degree at most three.
Given a family ${\mathcal H}$ of graphs, a graph is {\it ${\mathcal H}$-free} if it does not contain any graph in ${\mathcal H}$ as a subgraph.
A graph is {\it triangle-free} if it is $\{K_3\}$-free.
Staton \cite{s} proved that $\alpha(G) \leq 5\lvert V(G) \rvert/14$ for every triangle-free subcubic graph $G$.
This result is best possible as Fajtlowicz \cite{f} pointed out that the generalized Petersen graph $P(7,2)$ has $14$ vertices but no independent set of size $6$.
Heckman and Thomas \cite{ht} gave a short proof of Staton's Theorem and gave the following conjecture.

\begin{conj} \cite{ht} \label{the conjecture}
The fractional chromatic number of every triangle-free subcubic graph is at most $14/5$.
\end{conj}

At the moment when the first version of this article was submitted, the conjecture remained open, but recently it was confirmed by Dvo\v{r}\'{a}k, Sereni and Volec \cite{dsv}.
Nevertheless, it is still worthwhile to review the progress of this conjecture.
Let $G$ be a triangle-free subcubic graph.
Hatami and Zhu \cite{hz} proved that $\chi_f(G) \leq 3-3/64 \approx 2.953$, and $\chi_f(G) < 14/5$ if the shortest cycle of $G$ has length at least $7$.
Lu and Peng \cite{lp} improved the bound to $\chi_f(G) \leq 3-3/43 \approx 2.930$.
Ferguson, Kaiser and Kr\'{a}l' \cite{fkk} further improved that $\chi_f(G) \leq 32/11 \approx 2.909$.
In this paper, we prove the following theorem.
Note that $43/15 = 14/5 + 1/15 \approx 2.867$.

\begin{theorem} \label{main}
The fractional chromatic number of every triangle-free subcubic graph is at most $43/15$.
\end{theorem}

We remark that the proof of Dvo\v{r}\'{a}k et al. in \cite{dsv} for $\chi_f(G) \leq 14/5$ is not constructive, and they asked whether there exists an integer $t$ such that every subcubic triangle-free graph has a $(14t:5t)$-coloring.
On the other hand, our proof of Theorem \ref{main} is constructive, and we will show that every triangle-free subcubic graph has a $(516:180)$-coloring.

\section{Outline of the proofs and notations}
The idea of our proof of Theorem \ref{main} is not complicated, but it requires a large amount of machinery and number of technical lemmas to implement it.
Now, we give a superficial outline of the proof of Theorem \ref{main}.

For every graph $H$ and every vertex $v$ of $H$, we denote the degree of $v$ by $\deg_H(v)$.
First, given a minimum counterexample of Theorem \ref{main} $G$, we find a proper $3$-coloring $f$ of $G$.
Second, for each color $i \in [3]$, let $H_i$ be the subgraph of $G$ obtained by deleting the vertices $v$ with $f(v)=i$.
We then choose a function $\epsilon_i$ on $V(H_i)$ taking values in $\{0,1\}$.
Next, we find a coloring $g_i$ with $56$ colors such that every vertex $v$ with $f(v)=i$ receives $8$ colors, each remaining vertex $v$ receives $32-4\deg_{H_i}(v)-4\epsilon_i(v)$ colors, and every pair of adjacent vertices receive disjoint sets of colors.
For every vertex $v$, we denote $\epsilon_1(v)+\epsilon_2(v)+\epsilon_3(v)$ by $\epsilon(v)$.
In fact, we will prove that we can choose $\epsilon_1,\epsilon_2,\epsilon_3$ such that $\epsilon(v) \in \{0,1\}$ for every vertex $v$, and the set $\{v: \epsilon(v)=1\}$ is an independent set.
Third, by combining the three colorings $g_1,g_2$ and $g_3$, there exists a coloring $g$ with $168$ colors such that every vertex receives at least
$$(8+32+32)-4(\deg_{H_1}(v)+\deg_{H_2}(v)+\deg_{H_3}(v))-4\epsilon(v)= 72-4\deg_G(v)-4\epsilon(v)$$ 
colors (since $H_1,H_2,H_3$ are pairwise edge-disjoint).
Moreover, every pair of adjacent vertices receive disjoint sets of colors.
Finally, we assign $4$ extra colors to those vertices with $\epsilon$-value $1$ to obtain an $(172:60)$-coloring of $G$.
Notice that $172/60=43/15$, contradicting the assumption that $G$ is a minimum counterexample, so Theorem \ref{main} is proved.
With some extra work, we can construct a $(516:180)$-coloring of the given triangle-free subcubic graph.

We remark that if $\epsilon$ can be removed, then there exists an $(168,60)$-coloring of a minimum counterexample.
This will confirm Conjecture \ref{the conjecture} since $168/60=14/5$, and it will answer Dvo\v{r}\'{a}k et al.'s question by giving a $(504:180)$-coloring of the given triangle-free subcubic graph.
However, it is not clear how to remove $\epsilon$.

The paper is organized as follows.
In Section 3, we give a number of technical lemmas which formally state the notions in the idea just mentioned; then we prove Theorem \ref{main}.
We postpone the proofs of the majority of these lemmas to the other sections.
In Section 4, we investigate the structure of fractionally critical graphs, which play the role of minimum counterexamples of Theorem \ref{main}.
In particular, we show that such a graph has a proper $3$-coloring $f$ such that $H_i$ has a ``good" structure.
This fulfills the first step of the idea we just mentioned as well as prepares us for the remaining steps.
The formal definition of ``good" structure will be stated in Section 3.
In fact, we prove a more general result than we need in this paper in hopes that it might be useful in future work.
In Section 5, we study a list-version of fractional colorings for graphs having those ``good" structures to implement the second and the third steps of the idea just mentioned.
More precisely, the way we construct coloring $g_i$ of $H_i$ is to first assign $8$ colors to each vertex $v$ with $f(v)=i$, and then extend the coloring to the remaining vertices.
When some neighbors of a vertex $v$ are pre-colored, the available colors for $v$ are limited, and it is the reason we study the list-version of fractional colorings in Section 5.
In Section 6, we complete the proof of the lemmas stated in Section 2.
In Section 7, we construct $(516:180)$-colorings for triangle-free subcubic graphs.
Finally, we make some concluding remarks in Section 8.

In the rest of this section, we introduce some terminology.
In this paper, graphs do not contain multiple edges or loops, unless specifically mentioned.
For any subset $S$ of $V(G)$, we define the {\it neighborhood} $N_G(S)$ of $S$ to be the set of vertices which are not in $S$ but are adjacent to a vertex in $S$.
Also, we define $N_G[S]$ to be the set $N_G(S) \cup S$.
We write $N_G(\{v\})$ as $N_G(v)$ and $N_G[\{v\}]$ as $N_G[v]$ for short.
If $u$ and $v$ are two vertices in $G$, then $G+uv$ is the graph obtained from $G$ by adding edge $uv$, and $G/uv$ is the graph obtained from $G$ by identifying $u,v$ and then deleting resulting loops and parallel edges.
We define $G[S]$ to be the subgraph induced by the set $S$ when $S$ is a subset of $V(G)$.
The {\it degree} of a vertex $v$ in $G$, denoted by $\deg_G(v)$, is the number of edges incident with $v$.
$G$ is {\it cubic} if every vertex in $G$ is of degree $3$.
A {\it leaf} is a vertex of degree one, and a {\it support vertex} is a vertex that is adjacent to a leaf.
A {\it matching} $M$ is a subset of edges such that no two edges in $M$ have a common end; we say that $M$ {\it saturates} a vertex $v$ if $v$ is an end of some edge in $M$.
Given any set $S$, we define the function $1_S$ by letting $1_S(x)=1$ if $x \in S$, and $1_S(x)=0$ otherwise.

A {\it digraph} $D$ is a graph equipped with an orientation of edges.
We say that a vertex $u$ is {\it pointed} by $v$ (or $v$ {\it points to} $u$) if there exists an edge with the head $u$ and the tail $v$.
For each vertex $v$ of $D$, we define $N^+_D(v)$ (and $N^-_D(v)$, respectively) to be the set of vertices that are pointed by $v$ (point to $v$, respectively).
The {\it out-degree} $\deg_D^+(v)$ ({\it in-degree} $\deg_D^-(v)$, respectively) of $v$ is $\lvert N^+_D(v) \rvert$ ($\lvert N^-_D(v) \rvert$, respectively).

Let $H$ be a subgraph of a graph $G$.
We say a function $f$ defined on $V(H)$ can be {\it extended} to a function $g$ defined on $V(G)$ if $g(v)=f(v)$ for all $v \in V(H)$.
We also say that $f$ can be extended to $G$ in this case.
Note that $f$ can always be extended to $g$ if there is no requirement for $g$.
But in our application, we require that $g$ satisfies some extra conditions.

\section{Proof of Theorem \ref{main}}

Given a graph $G$ and a function $F: V(G) \rightarrow 2^{[14]}$, we say that $f:V(G) \rightarrow 2^{[14]}$ is an {\it $F$-avoiding coloring} if $f(v)$ is disjoint from $F(v) \cup f(u)$ for every pair of adjacent vertices $u$ and $v$.

Let $f$ be a proper $3$-coloring of a subcubic graph $G$.
Given integers $i,j$ such that $1 \leq i < j \leq 3$, define $G_{i,j}$ to be the subgraph of $G$ induced by $f^{(-1)}(i) \cup f^{(-1)}(j)$, and let $C$ be a union of some components of the subgraph $G_{i,j}$.
For every vertex $v$ in $N_G(C)$, we define $n_C(v)$ to be the number of vertices in $C$ adjacent to $v$.
We say that $B$ is a {\it boundary-graph} of $C$ if $B$ is a graph such that $V(B)=N_G(C)$ and $B$ does not contain a triangle $x_1x_2x_3$ such that $n_C(x_i)=1$ for $1 \leq i \leq 3$.
We say that $(B, {\mathcal J})$ is a {\it boundary-pair} of $C$ if $B$ is a boundary-graph of $C$ and ${\mathcal J}$ is a collection of pairwise disjoint subsets $J$ of $V(B)$ such that every vertex $v$ in $J$ satisfies that $n_C(v) \geq 2$.
Notice that if $(B_1,\J_1)$ is a boundary-pair of $C_1$ and $(B_2,\J_2)$ is a boundary-pair of $C_2$, where $C_1$ is disjoint from $C_2$, then $(B_1 \cup B_2, \J_1 \cup \J_2)$ is a boundary-pair of $C_1 \cup C_2$.

Given a boundary-pair $(B, {\mathcal J})$, we say that a function $F: V(C) \rightarrow 2^{[14]}$ is {\it $(B, {\mathcal J})$-compatible} if there is a $(14:2)$-coloring $h$ of $B$ such that each set $S \in {\mathcal J}$ contains two vertices $x_S,y_S$ with $h(x_S) = h(y_S)$, and $F(v) = \bigcup_{u \in N_G(v) \cap V(B)} h(u)$ for every vertex $v$ of $C$.
Observe that $\lvert F(v) \rvert \leq 2 \lvert N_G(v) \cap V(B) \rvert \leq 6-2\deg_C(v)$ for every vertex $v$ in $C$ and every $(B, {\mathcal J})$-compatible $F$.
Note that when every vertex in $N(C)$ is pre-colored by $h$, then for every $v$, $F(v)$ is the set of colors that $v$ cannot use.

Once the vertices in $N_G(C)$ are colored, we want to extend the coloring to $C$.
We say that a boundary-pair $(B, {\mathcal J})$ of $C$ {\it penetrates} $C$ if there is an independent set $I \subseteq \{v \in V(C): \deg_C(v)=3\}$ in $C$ such that for every $(B, {\mathcal J})$-compatible $F$, there are $F$-avoiding colorings $g_1, g_2$ of $C$ such that $\lvert g_1(v) \rvert + \lvert g_2(v) \rvert = 16 - 2\deg_C(v) - 2 \cdot 1_I(v)$ for every vertex $v$ of $C$.
We say that a boundary-pair $(B_1, {\mathcal J}_1)$ of $C$ {\it cooperates} with a boundary-pair $(B_2, {\mathcal J}_2)$ of $C$ if there exists an independent set $I \subseteq \{v \in V(C): \deg_C(v) =3\}$ of $C$ such that given $F_1,F_2:V(G) \rightarrow 2^{[14]}$, where $F_1$ is $(B_1,\J_1)$-compatible and $F_2$ is $(B_2,\J_2)$-compatible, there exist an $F_1$-avoiding coloring $g_1$ of $C$ and an $F_2$-avoiding coloring $g_2$ of $C$ such that $\lvert g_1(v) \rvert + \lvert g_2(v) \rvert = 16-2\deg_C(v)-2 \cdot 1_I(v)$ for every vertex $v$ of $C$.

Given two positive integers $a$, $b$, a set $S$, and two functions $f:S \rightarrow 2^{[a]}$ and $g:S \rightarrow 2^{[b]}$, define $f \uplus g: S \rightarrow 2^{[a+b]}$ by letting $f \uplus g (v) = \{x, y+a: x \in f(v), y \in g(v)\}$.

All but the following lemma given in this section are laborious to prove, so we postpone the proofs to other sections.
On the other hand, the proof of the following lemma shows how the notion of penetrations and cooperations help us define an $(172:60)$-coloring, so we include the proof here.

\begin{lemma} \label{pene and coop}
Let $f$ be a proper $3$-coloring of $G$.
Given integers $s,t$ such that $1 \leq s <t \leq 3$, define $G_{s,t}$ to be the subgraph of $G$ induced by $f^{(-1)}(s) \cup f^{(-1)}(t)$.
Let $s,t,p$ be positive integers such that $1 \leq s < t \leq 3$, and let $C_1, C_2, ..., C_{r_{s,t}}$ be the components of $G_{s,t}$.
Assume that there are boundary-pairs $(B_i, {\mathcal J}_i)$ that penetrate $C_i$ for $1 \leq i \leq p$, and there are boundary-pairs $(B_{i,1}, {\mathcal J}_{i,1})$ and $(B_{i,2}, {\mathcal J}_{i,2})$ for $p+1 \leq i \leq r_{s,t}$ such that $(B_{i,1}, {\mathcal J}_{i,1})$ cooperates with $(B_{i,2}, {\mathcal J}_{i,2})$ for each $i$ such that $p+1 \leq i \leq r_{s,t}$.
Given $k=1,2$, let $D_k = \bigcup_{1 \leq i \leq p} B_i \cup \bigcup_{p+1 \leq i \leq r_{s,t}} B_{i,k}$ and ${\mathcal S}_k = \bigcup_{1 \leq i \leq p} {\mathcal J}_i \cup \bigcup_{p+1 \leq i \leq r_{s,t}} {\mathcal J}_{i,k}$.
If there are a $(D_1, {\mathcal S}_1)$-compatible function and a $(D_2, {\mathcal S}_2)$-compatible function, then there are an independent set $I_{s,t} \subseteq \{v \in V(G_{s,t}): \deg_{G_{s,t}}(v)=3\}$ and a function $g_{s,t}: V(G) \rightarrow 2^{[56]}$ such that $\lvert g_{s,t}(v) \rvert = 8$ if $f(v) \not \in \{s,t\}$, and $\lvert g_{s,t}(v) \rvert = 32-4 \deg_{G_{s,t}}(v) - 4 \cdot 1_{I_{s,t}}(v)$ if $f(v) \in \{s,t\}$, and $g_{s,t}(x) \cap g_{s,t}(y) = \emptyset$ for every pair of adjacent vertices $x$ and $y$.
Furthermore, if such $g_{1,2}, g_{1,3}$ and $g_{2,3}$ exist, then $G$ has an $(172: 60)$-coloring, and $\chi_f(G) \leq 43/15$.
\end{lemma}

\begin{pf}
Let $F_1$ and $F_2$ be a $(D_1, {\mathcal S}_1)$-compatible function and a $(D_2, {\mathcal S}_2)$-compatible function, respectively.
Observe that for every integer $i$ such that $1 \leq i \leq p$, the function $F_1$ and $F_2$ restricting on $B_i$ is $(B_i,{\mathcal J}_i)$-compatible, respectively.
Similarly, for every $i,j$ such that $p+1 \leq i \leq r_{s,t}$ and $1 \leq j \leq 2$, $F_j$ restricting on $B_{i,j}$ is $(B_{i,j}, {\mathcal J}_{i,j})$-compatible.
Since $(B_i, {\mathcal J}_i)$ penetrates $C_i$ for $1 \leq i \leq p$, there are independent sets $I_{s,t,i} \subseteq \{v \in V(C_i): \deg_{G_{s,t}}(v)=3\}$, $F_1$-avoiding colorings $h_{i,1}$ and $h'_{i,1}$, and $F_2$-avoiding colorings $h_{i,2}$ and $h'_{i,2}$ such that $\lvert h_{i,1}(v) \rvert + \lvert h'_{i,1}(v) \rvert = \lvert h_{i,2}(v) \rvert + \lvert h'_{i,2}(v) \rvert = 16 - 2\deg_{G_{s,t}}(v) - 2 \cdot 1_{I_{s,t,i}}(v)$ for vertex $v$ in $C_i$ and for $1 \leq i \leq p$.
Similarly, since $(B_{i,1}, {\mathcal J_{i,1}})$ cooperates with $(B_{i,2}, {\mathcal J_{i,2}})$ for $p+1 \leq i \leq r_{s,t}$, there are independent sets $I_{s,t,i} \subseteq \{v \in V(C_i): \deg_{G_{s,t}}(v)=3\}$ and $F_k$-avoiding colorings $h''_{i,k}$, for $k=1,2$, such that $\lvert h''_{i,1}(v) \rvert + \lvert h''_{i,2}(v) \rvert = 16 - 2\deg_{G_{s,t}}(v) - 2 \cdot 1_{I_{s,t,i}}(v)$ for vertex $v$ in $C_i$ and for $p+1 \leq i \leq r_{s,t}$.

Let $I_{s,t} = \bigcup_{1 \leq i \leq r_{s,t}} I_{s,t,i}$, and let $h_1,h_2$ be a $(14:2)$-coloring of $(D_1, {\mathcal S}_1)$ and $(D_2, {\mathcal S}_2)$, respectively, such that for each $S \in {\mathcal S}_i$, there exist two vertices in $S$ with the same $h_i$-value and $F_i(v) = \bigcup_{u \in N_G(v) \cap V(D_i)} h_i(u)$, for each vertex $v$ of $G_{s,t}$ and $i=1,2$.
Define $g_{s,t}(v)=(h_1 \uplus h_1 \uplus h_2 \uplus h_2)(v)$ for $v \in V(G)-V(G_{s,t})$, 
$g_{s,t}(v) = (h_{i,1} \uplus h'_{i,1} \uplus h_{i,2} \uplus h'_{i,2})(v)$ for $v \in V(C_i)$ and $1 \leq i \leq p$, and $g_{s,t}(v) = (h''_{i,1} \uplus h''_{i,1} \uplus h''_{i,2} \uplus h''_{i,2})(v)$ for $v \in V(C_i)$ and $p+1 \leq i \leq r_{s,t}$.
It is clear that $I_{s,t}$ is an independent set, and $g_{s,t}$ maps every pair of adjacent vertices to two disjoint sets in $2^{[56]}$.
In addition, $\lvert g_{s,t}(v) \rvert=8$ for $v \in V(G)-V(G_{s,t})$ and $\lvert g_{s,t}(v) \rvert = 32 - 4\deg_{G_{s,t}}(v) - 4 \cdot 1_{I_{s,t}}(v)$ for every vertex $v$ in $G_{s,t}$.

If $g_{1,2}, g_{1,3}$ and $g_{2,3}$ exist, then define $g = g_{1,2} \uplus g_{1,3} \uplus g_{2,3}$.
Hence, $g(v) \subseteq [168]$ and $\lvert g(v) \rvert = 72 -4\deg_G(v) - 4(1_{I_{1,2}}(v) + 1_{I_{1,3}}(v) + 1_{I_{2,3}}(v))$ for every vertex $v$ in $G$.
Let $I = I_{1,2} \cup I_{1,3} \cup I_{2,3}$.
Note that every vertex in $I_{1,2}$, $I_{1,3}$ or $I_{2,3}$ is adjacent to three vertices of the same $f$-value, so these three sets are pairwise disjoint, and $I$ is an independent set as well.
Define $g':V(G) \rightarrow 2^{[172]}$ by assigning $g'(v)=g(v)$ if $v \not \in I$, and $g'(v)=g(v) \cup \{169,170,171,172\}$ if $v \in I$.
Consequently, $g'$ is an $(172:60)$-coloring since the maximum degree of $G$ is at most three, and hence $\chi_f(G) \leq 172/60 = 43/15$.
\end{pf}

\bigskip

The following lemma shows the existence of boundary-pairs of a path on an odd number of vertices that cooperate with each other or penetrate the path.

\begin{lemma} \label{admit odd path}
Let $G$ be a triangle-free subcubic graph, and let $f$ be a proper $3$-coloring of $G$.
Let $k$ be an odd number, and let $C$ be a component of a subgraph of $G$ induced by two color classes of $f$.
Assume that $C$ is a path on $k$ vertices.
If $k \neq 3$, then there is a boundary-graph $B$ of $C$ such that $\deg_B(v) \leq 2n_C(v)$ for every vertex $v$ in $B$, and $(B, \emptyset)$ penetrates $C$.
If $k=3$, then there are boundary-graphs $B_1$ and $B_2$ of $C$ such that $\deg_{B_i}(v) \leq 2n_C(v)$ for every vertex $v$ of $B_i$ and for every $i=1,2$, and $(B_1, \emptyset)$ cooperates with $(B_2, \emptyset)$.
\end{lemma}

Let $L_0$ be the graph obtained by adding two vertices to a path on four vertices such that each of the new vertices is adjacent to the both ends of the path.
We say that the tuple $(v_1, v_2, v_3, v_4, u_1, u_2)$ of six distinct vertices of $G$ {\it forms a copy of $L_0$} in $G$ if these six vertices induce a graph isomorphic to $L_0$ and $u_1v_1v_2v_3v_4u_1$ and $u_2v_1v_2v_3v_4u_2$ are the two 5-cycles in $L_0$.
Notice that two different copies of $L_0$ in $G$ are not necessary vertex-disjoint.

Given any proper 3-coloring $f$ of $G$, we say that $(v_1, v_2, v_3, v_4, u_1, u_2)$ is a {\it rainbow} copy of $L_0$ with respect to $f$ (or {\it $f$-rainbow}) if $(v_1, v_2, v_3, v_4, u_1, u_2)$ forms a copy of $L_0$, and $f(u_1)=f(u_2)=\pi(1)$, $f(v_1)=f(v_3)=\pi(2)$, $f(v_2)=f(v_4)=\pi(3)$ for some permutation $\pi$ of $\{1,2,3\}$.

For even numbers $i,j$, where $i, j \geq 4$, we denote the path on $i$ and $j$ vertices by $P_i$ and $P_j$, respectively, and we define $H_{i,j}$ to be the graph that is obtained from two disjoint paths $P_i$ and $P_j$ by adding an edge incident with one support vertex in $P_i$ and one support vertex in $P_j$.
The family ${\mathcal H}$ consists of graphs $H_{i,j}$ for all even numbers $i,j$ with $i,j \geq 4$.

The following lemma shows the existence of boundary-pairs of a path on even number of vertices or a graph in ${\mathcal H}$ that penetrate this graph.
Furthermore, in order to prove Theorem \ref{main}, we need some specific structures for those boundary-pairs.

\begin{lemma} \label{admit even path and H}
Let $G$ be a triangle-free subcubic graph, and let $f$ be a proper $3$-coloring of $G$ such that there are no rainbow copies of $L_0$ with respect to $f$.
Let $C$ be a component of a subgraph of $G$ induced by two color classes of $f$.
If $C$ is a path on even number of vertices or it is a graph in ${\mathcal H}$, then there are boundary-pairs $(B_1, {\mathcal J_1})$ and $(B_2, {\mathcal J_2})$ of $C$ such that the following hold.
	\begin{enumerate}
		\item Both $(B_1, {\mathcal J_1})$ and $(B_2, {\mathcal J_2})$ penetrate $C$. 
		\item Either for $i=1,2$, every vertex in $B_i$ has $\deg_{B_i}(v) \leq 2n_C(v)$, or there exist $w_1 \in V(B_1)$, $w_2 \in V(B_2)$ such that for $i=1,2$, $\deg_{B_i}(w_i)=2n_C(w_i)+1$ and $\deg_{B_i}(v) \leq 2n_C(v)$ for every $v \in V(B_i) - \{w_i\}$.
		\item If $w_1$ and $w_2$ exist, then ${\mathcal J_1}={\mathcal J_2} = \emptyset$, and either $w_1=w_2$ and $n_C(w_1)=n_C(w_2)=3$, or $\deg_{B_{3-i}}(w_i) \leq 2n_C(w_i)-1$ and $n_C(w_i) \leq 2$ for each $i=1,2$.
		\item If $w_1$ and $w_2$ exist, then for $i=1,2$, there exists a vertex $w'_i$ in the component of $B_i$ containing $w_i$ such that $\deg_{B_i}(w'_i) \leq 2 n_C(w'_i)-1$.
		\item Each ${\mathcal J_i}$ contains at most one set.
		\item If ${\mathcal J}_i \neq \emptyset$, then let $J_i$ be the set in ${\mathcal J}_i$, and the following hold:
			\begin{enumerate}
				\item $\lvert \bigcup_{v \in J_i, n_C(v) = 2}N_{B_i}(v) \rvert \leq 2$.
				\item Either $J_i$ is an independent set of size $2$ in $B_i$ such that one vertex $x$ in $J_i$ has $n_C(x)=3$, or $J_i$ is an independent set of size $3$ in $B_i$ such that every vertex $x$ in $J_i$ has $n_C(x)=2$.
				\item If $\lvert J_i \rvert=3$, then either $\lvert \bigcup_{v \in J_i}N_{B_i}(v) \rvert \leq 1$, or there is a vertex $v \in V(B_i)$ such that $\lvert N_{B_i}(v) \cap J_i \rvert = 2$.
				\item If $\lvert J_i \rvert=2$, then $\lvert \bigcup_{v \in J_i}N_{B_i}(v) \rvert \leq 5$; if $\lvert J_i \rvert=2$ and there exists a vertex $x \in J_i$ with $n_C(x)=2$, then $\lvert \bigcup_{v \in J_i}N_{B_i}(v) \rvert \leq 3$.
			\end{enumerate}
	\end{enumerate}
\end{lemma}

We say that a graph $G$ is {\it good} if $G$ is bipartite, subcubic, and the following hold:
\begin{enumerate}
	\item[(G1)] No vertex of degree three is adjacent to three vertices of degree three in $G$.
	\item[(G2)] For every pair of adjacent vertices $x,y$ of degree three in $G$, either at least one of $x,y$ is adjacent to a leaf, or each of $x,y$ is adjacent to a support vertex of degree two.	
	\item[(G3)] If $x,y,z$ are three vertices of degree three such that $xyz$ is a path in $G$, and each $x$ and $z$ is adjacent to a leaf, and the neighbor $y'$ of $y$ other than $x$ and $z$ has degree two, then the neighbor $y'_1$ of $y'$  other than $y$ is a leaf.
	\item[(G4)] If $x,y,z$ are three vertices of degree three such that $xyz$ is a path in $G$,
$\deg_G(z_1) \geq 2$, and $\deg_G(z_2) \geq 2$, where $N_G(z) = \{z_1,z_2,y\}$, then $z_1$ or $z_2$ is a support vertex of degree two.
\end{enumerate}

We say that a $3$-coloring $f$ of $G$ is {\it good} if it is proper, $G$ contains no $f$-rainbow copies of $L_0$, and every pair of color classes of $f$ induces a good graph.
The following lemma shows the existence of a boundary-pair that penetrates a good graph.

\begin{lemma}  \label{admit non-path non-H}
Let $G$ be a triangle-free subcubic graph, and let $f$ be a good $3$-coloring of $G$.
Let $C$ be a component of a subgraph of $G$ induced by two color classes of $f$.
If $C$ is neither a path nor a graph in ${\mathcal H}$, then there exists a boundary-graph $B$ of $C$ such that $\deg_B(v) \leq 2n_C(v)$ and $(B, \emptyset)$ penetrates $C$.
\end{lemma}

$G$ is {\it fractionally $t$-critical} if $\chi_f(G) > t$ but $\chi_f(H) \leq t$ for every proper subgraph $H$ of $G$.
Readers who are familiar with the notion of coloring critical graphs might notice that the equality in the above definition is on the side of subgraphs instead of on the side of $G$ as in the definition of coloring critical graphs.
Even though our definition of fractionally critical graphs is not consistent with coloring critical graphs, we think that it is more reasonable to define critical graph in this way when we deal with fractional coloring.

\begin{lemma} \label{good 3-coloring}
If $G$ is a fractionally $t$-critical triangle-free subcubic graph with $t \geq 8/3$, then $G$ has a good $3$-coloring.
\end{lemma}

The following lemma is proved in \cite{aeltw} and probably elsewhere.

\begin{lemma} \cite{aeltw} \label{orientation}
If $G$ is a graph (allowing parallel edges), then there is an orientation of $G$ such that the in-degree and the out-degree of every vertex $v$ is at least $\lfloor \deg_G(v)/2 \rfloor$.
\end{lemma}

Theorem \ref{main} is an immediate consequence of Lemma \ref{good 3-coloring} and the following theorem.

\begin{theorem} \label{strong main}
Let $G$ be a triangle-free subcubic graph.
If there exists a good $3$-coloring of $G$, then $G$ has an $(172:60)$-coloring.
\end{theorem}

\begin{pf}
Let $s,t$ be integers such that $1 \leq s < t \leq 3$ and let $G_{s,t}$ be the subgraph of $G$ induced by $f^{(-1)}(s) \cup f^{(-1)}(t)$.
For every component $C$ of $G_{s,t}$ isomorphic to a path on an even number of vertices or a graph in ${\mathcal H}$, there exist boundary-pairs $(B_{C,1},{\mathcal J}_{C,1})$ and $(B_{C,2},{\mathcal J}_{C,2})$ of $C$ satisfying Statements 1-6 of Lemma \ref{admit even path and H}.
So for $i=1,2$, $B_{C,i}$ contains at most one vertex, denoted by $w_{C,i}$ if it exists, satisfying that $n_C(w_{C,i}) \leq 2$ and $\deg_{B_{C,i}}(w_{C,i})=2n_C(w_{C,i})+1$.
Observe that either both $w_{C,1}$ and $w_{C,2}$ exist, or none of them exists.
And when both of them exist, $\deg_{B_{C,3-i}}(w_{C,i}) \leq 2n_C(w_{C,i})-1$.
We say that a component $C'$ of $G_{s,t}$ isomorphic to a path on even number of vertices or a graph in ${\mathcal H}$ is {\it dangerous} if $w_{C',1}$ and $w_{C',2}$ exist.
Note that $w_{C,i}$ could be equal to $w_{C',i'}$ for some different components $C',C$ of $G_{s,t}$ and for some $i'$.
Construct a graph $A$ (allowing parallel edges), where $V(A) = \{w_{C',1}, w_{C',2}:C'$ is a dangerous component$\}$ and $E(A) = \{w_{C',1}w_{C',2}:C'$ is a dangerous component$\}$ (viewed as a multiset).
By Lemma \ref{orientation}, there is an orientation of $A$ such that every vertex $v$ in $A$ has in-degree and out-degree at most $\lceil \deg_A(v)/2 \rceil$.
Since we can exchange the name of $B_{C,1}$ and $B_{C,2}$, we may assume that every directed edge of $A$ in the orientation is from $w_{C,1}$ to $w_{C,2}$.
Note that every $w_{C,i}$ is adjacent to at most three components of $G_{s,t}$.
So the maximum degree of $A$ is at most three, and $\lvert \deg^+_A(v) - \deg^-_A(v) \rvert \leq 1$ for every vertex $v$ in $A$.

On the other hand, by Lemmas \ref{admit odd path} and \ref{admit non-path non-H}, for every component $C$ of $G_{s,t}$ not isomorphic to a path on an even number of vertices or a graph in ${\mathcal H}$, if $C$ is a path on three vertices, then there are boundary-pairs $(B_{C,1},\emptyset)$ and $(B_{C,2}, \emptyset)$ of $C$ such that $(B_{C,1}, \emptyset)$ cooperates with $(B_{C,2}, \emptyset)$ and $\deg_{B_{C,i}}(v) \leq 2 n_C(v)$ for every vertex $v$ in $B_{C,i}$ and $i=1,2$; otherwise, there is a boundary-pair $(B_C, \emptyset)$ penetrating $C$ and $\deg_{B_C}(v) \leq 2n_C(v)$ for every vertex $v$ in $B_C$.
For each $i=1,2$, construct a graph $H_i$ by defining $V(H_i)=N_G(G_{s,t})$ and $E(H_i)= \bigcup E(B_C) \cup \bigcup E(B_{C,i})$, where the second union runs through all components $C$ of $G_{s,t}$ which are isomorphic to paths on three or even numbers of vertices or graphs in ${\mathcal H}$, and the first union runs through all other components $C$ of $G_{s,t}$.
Note that for every vertex $v \in V(H_i)$, if $\deg_{B_{C,1}}(v)=\deg_{B_{C,2}}(v)=2n_C(v)+1$ for some component $C$, then $n_C(v)=3$ and $\deg_{H_i}(v)=7$, otherwise, 
$$\begin{array}{lcl}
\deg_{H_i}(v) & \leq & 2 \sum_C n_C(v) + \lvert \{C:v=w_{C,i}\} \rvert - \lvert \{C: v=w_{C,3-i}\} \rvert \\ 
 & \leq & 2\deg_G(v) + \lvert \deg^+_A(v) - \deg^-_A(v) \rvert \leq 7.
\end{array}$$
Similarly, we define ${\mathcal J}_i = \bigcup {\mathcal J_{C,i}}$ for each $i=1,2$, where the union runs through all components $C$ of $G_{s,t}$ which are isomorphic to paths on even numbers of vertices or graphs in ${\mathcal H}$.

Note that it is sufficient to show that for each $i=1,2$, there is a proper $7$-coloring $h_i:V(H_i) \rightarrow [7]$ of $H_i$ such that every set in ${\mathcal J}_i$ has two vertices of the same $h_i$-value.
If such colorings $h_1,h_2$ exist, then for $i=1,2$, define $F_i: V(G_{s,t}) \rightarrow 2^{[14]}$ by setting $F_i(v) = \bigcup_{u \in N_G(v) \cap V(H_i)} h_i'(u)$, where $h'_i$ is the $(14:2)$-coloring of $H_i$ obtained by setting $h_i'(v) = \{2h_i(v)-1, 2h_i(v)\}$ for every vertex $v$ in $V(H_i)$.
It is clear that $F_i$ is $(H_i, {\mathcal J}_i)$-compatible in this case.
Since it is true for every $1 \leq s < t \leq 3$, Lemma \ref{pene and coop} implies that $G$ has an $(172:60)$-coloring.

Now we show that there is a proper $7$-coloring $h_i$ of $H_i$ such that every set in ${\mathcal J}_i$ contains two vertices having the same $h_i$-value.
The proofs for $H_1$ and $H_2$ are the same, so we write $H$ instead of $H_1,H_2$ for convenience.
Also let ${\mathcal J_C}$, $B_C$ and $w_C$ be the corresponding ${\mathcal J_{C,i}}$, $B_{C,i}$ and $w_{C,i}$, respectively, for every component $C$ of $G_{s,t}$.

First, we show that $H$ is properly $7$-colorable.
By Brooks' Theorem, it is sufficient to show that no component of $H$ is isomorphic to $K_8$, since the maximum degree of $H$ is at most $7$.
Suppose to the contrary, and let $X$ be a component of $H$ isomorphic to $K_8$.
The following two claims are clear.

\noindent {\bf Claim 1:} If $v \in V(X)$, then $v$ is adjacent to at least one component $C$ of $G_{s,t}$ such that $v=w_C$.

\noindent {\bf Claim 2:} If $v \in V(X)$ and $\deg_{B_C}(v) \leq 2n_C(v)-1$ for some $C$, then $n_C(v)=1$ and there are two other components $C',C''$ such that $v=w_{C'}=w_{C''}$.

By Statement 4 of Lemma \ref{admit even path and H}, for every vertex $x$ in $X$ such that $x=w_C$ for some component $C$ of $G_{s,t}$, there exists a vertex $w_C'$ that is in the component of $B_C$ containing $w_C$ such that $\deg_{B_C}(w_C') \leq 2n_C(w'_C)-1$.
In particular, $w'_C$ is in $X$.
Construct a directed graph $X'$ (allowing parallel edges here) by setting $V(X')=V(X)$ and $(w_C,w'_C) \in E(X')$ for every component $C$ of $G_{s,t}$.
Note that every vertex in $X'$ has out-degree at least $1$ (by Claim 1) and in-degree at most $1$ (by Claim 2).
Furthermore, if the vertex has in-degree $1$, then it has out-degree $2$ by Claim 2.
So the sum of the out-degrees of vertices in $X'$ is greater than the sum of the in-degrees of vertices in $X'$, a contradiction.
Hence, $H$ is $7$-colorable.

Denote ${\mathcal J_{C}}$ by $\{J_C\}$ for every nonempty ${\mathcal J_{C}}$.
Note that for every components $C,C'$ of $G_{s,t}$, we have that every vertex $v$ in $J_C$ satisfies that $n_C(v) \geq 2$, and no triangle in $B_{C'}$ contains three vertices whose $n_{C'}$-values are $1$.
In addition, if $\lvert J_C \rvert=2$, then some vertex $v$ in $J_C$ has $n_C(v)=3$.
Hence, $J_C$ is not a clique in $H$.
Moreover, every vertex $v$ in $J_C$ has degree at most $6$ in $H$ by Statements 6(a) and (d) of Lemma \ref{admit even path and H}.

\noindent {\bf Claim 3:} Every proper $7$-coloring of $H-J_C$ can be extended to a proper $7$-coloring of $H$ such that a pair of vertices in $J_C$ receive the same color if one of the following is satisfied:
	\begin{enumerate}
		\item $\lvert J_C \rvert=2$;
		\item there are two edges in $H$ with the both ends in $J_C$;
		\item there is a vertex in $J_C$ incident with at most two edges in $E(H)-E(B_C)$, and there is at most one edge in $H$ with the both ends in $J_C$;
		\item $\lvert \bigcup_{u \in J_C} N_{B_C}(u) \rvert \leq 1$.
	\end{enumerate}

\noindent {\bf Proof of Claim 3:}
	\begin{enumerate}
		\item If $\lvert J_C \rvert=2$, then let $J_C = \{x,y\}$.
By Statement 6 (b) of Lemma \ref{admit even path and H}, we know that one vertex in $J_C$, say $x$, satisfies that $n_C(x)=3$.
By Statement 6 (d) of Lemma \ref{admit even path and H}, either $n_C(y)=3$ and $\lvert N_{B_C}(x) \cup N_{B_C}(y) \rvert \leq 5$, or $n_C(y)=2$, $\lvert N_{B_C}(x) \cup N_{B_C}(y) \rvert \leq 3$ and $\deg_{H-B_C}(y) \leq 3$.
So $\lvert N_H(x) \cup N_H(y) - J_C \rvert \leq 6$.
Therefore, we can extend any proper $7$-coloring of $H-J_C$ to $H$ such that $x$ and $y$ receive the same color.

		\item We may assume that the previous case does not happen, so $\lvert J_C \rvert =3$ for the rest of the proof of Claim 3.
Let $J_C=\{x,y,z\}$ and $xy, yz \in E(H)$.
Since $n_C(v) \geq 2$ for every $v \in J_C$, there is at most one component $C' \neq C$ of $G_{s,t}$ such that $J_C \cap V(B_{C'}) \neq \emptyset$.
Also, $J_C$ is an independent set in $B_C$ and $x,y,z$ cannot form a triangle in $B_{C'}$, so $xz \not \in E(H)$.
In addition, since $\lvert N_{B_C}(x) \cup N_{B_C}(z) \rvert \leq 2$ by Statement 6(a) of Lemma \ref{admit even path and H} and both $\deg_{B_{C'}-\{xy\}}(x)$ and $\deg_{B_{C'}-\{yz\}}(y)$ are at most $2$ for the component $C' \neq C$ of $G_{s,t}$ containing $x,y$, we know that $\lvert N_H(x) \cup N_H(z) - J_C \rvert \leq 6$, so every proper $7$-coloring $h$ of $H-J_C$ can be extended to $H$ such that $h(x)=h(z)$ by first defining $h(x)=h(z)$ to be a color in $[7]-( N_H(x) \cup N_H(z) - J_C)$ and then defining $h(y)$ to be any feasible color.

		\item Assume that there is a vertex in $J_C$, denoted by $a$, incident with at most two edges in $E(H)-E(B_C)$.
If there is exactly one edge in $H$ with the both ends in $J_C$, then by Statement 6(a) of Lemma \ref{admit even path and H}, there is a pair of non-adjacent vertices $u,v$ in $J_C$ such that $\lvert N_H(u) \cup N_H(v) - J_C \rvert \leq 6$, where $a \in \{u,v\}$, so every proper $7$-coloring of $H$ can be extended to a proper $7$-coloring of $H$ such that two vertices in $J_C$ receive the same color.

So we may assume that no edge in $H$ has both ends in $J_C$.
Denote the other two vertices in $J_C$ other than $a$ by $b$ and $c$.
We shall prove that there exists a pair of non-adjacent vertices $u,v$ in $J_C$ such that $\lvert h((N_H(u) \cup N_H(v)) - J_C) \rvert \leq 6$, and hence it can be extended to a proper $7$-coloring of $H$ such that $h(u)=h(v)$.
Suppose to the contrary.
Note that $\lvert (N_H(a) \cup N_H(b)) - J_C \rvert$ and $\lvert (N_H(a) \cup N_H(c)) - J_C \rvert$ are at most $7$.
So they are equal to $7$.
This implies that $(N_H(a) \cap N_H(b)) - (J_C \cup N_{B_C}(a) \cup N_{B_C}(b) \cup N_{B_C}(c)) = \emptyset = (N_H(a) \cap N_H(c)) - (J_C \cup N_{B_C}(a) \cup N_{B_C}(b) \cup N_{B_C}(c))$.
In addition, $\lvert h((N_H(a) \cup N_H(b))-J_C) \rvert = \lvert h((N_H(a) \cup N_H(c))-J_C) \rvert=7$.
So $h(N_H(b)-(J_C \cup N_{B_C}(a) \cup N_{B_C}(b) \cup N_{B_C}(c))) = h(N_H(c)-(J_C \cup N_{B_C}(a) \cup N_{B_C}(b) \cup N_{B_C}(c)))$.
Hence, $\lvert h((N_H(b) \cup N_H(c))-J_C) \rvert \leq 6$, a contradiction.
		\item Let $\lvert \bigcup_{u \in J_C} N_{B_C}(u) \rvert \leq 1$.
We may assume that the previous three cases do not happen.
Suppose the conclusion does not hold.
Then there does not exist a pair of non-adjacent vertices $u,v$ in $J_C$ such that $\lvert h((N_H(u) \cup N_H(v)-J_C) \rvert \leq 6$.
This implies that $J_C$ is an independent set in $H$, and there exist $x,y \in J_C$ such that $h(N_H(x)-(J_C \cup \bigcup_{u \in J_C}N_{B_C}(u)) = h(N_H(y)-(J_C \cup \bigcup_{u \in J_C} N_{B_C}(u)))$.
So $\lvert h((N_H(x) \cup N_H(y))-J_C) \rvert \leq 6$, a contradiction.
This completes the proof of Claim 3.
$\Box$
	\end{enumerate}

For each component $C$ of $G_{s,t}$, if any hypothesis of Claim 3 applies on $J_C$, then define $S_C = J_C$; otherwise, define $S_C = J_C \cup \{u_C\}$, where $u_C$ is a vertex in $\bigcup_{u \in J_C} N_{B_C}(u)$ such that $u_C$ is adjacent in $B_C$ to exactly two vertices in $J_C$.
Note that $u_C$ exists and $u_C \not \in J_C$ by Statements 6(b) and (c) of Lemma \ref{admit even path and H}.
By Claim 3, if $S_C=J_C$ and $h$ is a proper $7$-coloring of $H-J_C$, then $h$ can be extended to a proper $7$-coloring of $H$ such that at least two vertices in $J_C$ receive the same color.
On the other hand, if $u_C \in S_C \cap J_{C'}$ for some two distinct components $C$ and $C'$, then $n_C(u_C)=1$, $n_{C'}(u_C)=2$, $w_C$ does not exist, and there are at most two edges in $E(H)-E(B_{C'})$ incident with $u_C$, so $S_{C'} =J_{C'}$.

Now, we construct a proper $7$-coloring $h$ of $H$ such that each $J_C$ has two vertices getting the same $h$-value.
Note that $J_C \cap J_{C'} = \emptyset$ unless $C=C'$.
Set $H'= H - \bigcup_C S_C$, where the union runs through all components $C$ of $G_{s,t}$ such that ${\mathcal J}_C \neq \emptyset$.
Let $h$ be a proper $7$-coloring of $H - \bigcup_C S_C$.
Recall that $H$ is $7$-colorable, so such $h$ exists.
Pick a component $C^*$, where ${\mathcal J}_{C^*} = \{J_{C^*}\} \neq \emptyset$ and $S_{C^*} \neq J_{C^*}$ and $J_{C^*} \cap H' = \emptyset$, then extend $h$ to $H' \cup \{x_{C^*}, y_{C^*}\}$ such that $h(x_{C^*})=h(y_{C^*})$ for some distinct vertices $x_{C^*}$ and $y_{C^*}$ in $J_{C^*}$, and set $H'$ to be the subgraph of $H$ induced by $V(H') \cup \{x_{C^*},y_{C^*}\}$.
Note that since $u_{C^*}$ has not been assigned a color, such coloring extension exists by the argument for proving the fourth case of Claim 3.
Repeat this process until no such $C^*$ exists.
Note that $J_{C^*} \cap S_C = \emptyset$ for other component $C$ of $G_{s,t}$, since $J_{C^*} \neq S_{C^*}$.
So, at this moment, the vertices in $H$ that have not received $h$-values are $u_{C'}$ and an uncolored vertex in $J_{C'}$ for each component $C'$ with $J_{C'} \neq S_{C'}$, and all vertices in $J_{C''}$ for each component $C''$ with $J_{C''}=S_{C''}$.
Note that for every $C$ such that $J_C \neq S_C$, since $u_C$ is adjacent to at least two vertices in $J_C$, either two neighbors of $u_C$ have the same $h$-value, or one neighbor of $u_C$ has not been assigned an $h$-value.
Furthermore, every vertex in $J_C$ has degree at most six, so we can extend $h$ to $H-\bigcup_C J_C$, where the union runs through all components $C$ of $G_{s,t}$ such that $J_C=S_C$.
Then finally we can extend $h$ to $H$ by Claim 3.
This completes the proof of the theorem.
\end{pf}

\section{Structure of fractionally critical graphs}
The objective of this section is to prove Lemma \ref{good 3-coloring}.
Recall that $G$ is {\it fractionally $t$-critical} if $\chi_f(G) > t$ but $\chi_f(H) \leq t$ for every proper subgraph $H$ of $G$.
The first step is to give a list of forbidden subgraphs for every fractionally $t$-critical graph for $t \geq 8/3$.

\begin{lemma} \label{basic structure}
Let $G_1$ and $G_2$ be two induced subgraphs of a graph $G$ such that $G_1 \cup G_2 = G$ and $G_1 \cap G_2$ is a clique.
Let $a,b,c,d$ be positive integers such that $a/b \geq c/d$.
If $G_1$ has an $(a:b)$-coloring and $G_2$ has a $(c:d)$-coloring, then $G$ has an $(as/b:s)$-coloring, where $s$ is the least common multiple of $b$ and $d$.
As a result, given $t >0$, every fractionally $t$-critical graph is $2$-connected, and it has no vertex-cut which induces a clique.
\end{lemma}

\begin{pf}
Since $G_1$ has an $(a:b)$-coloring $f_1$ and $G_2$ has a $(c:d)$-coloring $f_2$, there exist an $(as/b:s)$-coloring of $G_1$ and a $(cs/d:s)$-coloring of $G_2$.
As $G_1 \cap G_2$ is a clique, every pair of vertices in $G_1 \cap G_2$ receive disjoint sets of size $s$ by $f_1$ and $f_2$, respectively.
By swapping colors, we may assume that $f_1(v)=f_2(v)$ for every vertex $v$ in $G_1 \cap G_2$.
Therefore, there exists an $(as/b:s)$-coloring $f$ of $G$ such that $f(v)=f_1(v)$ for every $v \in V(G_1)$, and $f(v)=f_2(v)$ otherwise.
\end{pf}

\bigskip

Lu and Peng \cite{lp} gave the following result for graphs that have a vertex-cut with size $2$.

\begin{lemma} \cite{lp} \label{2 cut set}
Let $G_1$ and $G_2$ be two induced subgraph of a graph $G$ such that $V(G_1) \cup V(G_2) = V(G)$ and $V(G_1) \cap V(G_2) = \{u,v\}$, where $\{u,v\}$ is a vertex-cut.
	\begin{enumerate}
		\item If $uv$ is an edge in $G$, then $\chi_f(G) = \max\{\chi_f(G_1), \chi_f(G_2)\}$.
		\item If $uv$ is not an edge in $G$, then $\chi_f(G) \leq \max\{\chi_f(G_1), \chi_f(G_2+uv), \chi_f(G_2/uv)\}$.
	\end{enumerate}
\end{lemma}

The following lemma is a simple but useful observation to extend an $(8:3)$-coloring, and we leave the proof to the reader.

\begin{lemma} \label{extend 8:3}
Let $G$ be a graph, and let $f$ be an $(8:3)$-coloring of $G$.
Let $u,v$ be two distinct vertices of $G$.
Let $H$ be the graph obtained from $G$ by adding a new vertex $x$ and two edges $ux,vx$, and let $H'$ be the graph obtained from $G$ by adding two new vertices $y,z$ and three edges $uy,yz,zv$.
If $f(u)$ is not disjoint from $f(v)$, then $f$ can be extended to an $(8:3)$-coloring of $H$.
If $f(u) \neq f(v)$, then $f$ can be extended to an $(8:3)$-coloring of $H'$.
\end{lemma}

\begin{figure}
\unitlength=1mm
\begin{picture}(120, 150)

\put(0,175){$R_0$}
\multiput(0,155)(12.5,0){4}{\circle*{2}}
\multiput(18.25,175)(0,-40){2}{\circle*{2}}
\drawline[1000](0,155)(37,155)
\drawline[1000](0,155)(18.25,175)(37,155)(18.25,135)(0,155)
\put(21.3625,165){\circle*{2}}
\drawline[1000](18.25,175)(24.75,155)
\put(-5,150){$123$}
\put(7.5,150){$456$}
\put(20,150){$178$}
\put(35,150){$234$}
\put(18,177){$567$}
\put(18,130){$567$}
\put(14,165){$234$}

\put(50,175){$R_1$}
\multiput(50,155)(12.5,0){4}{\circle*{2}}
\multiput(68.25,175)(0,-40){2}{\circle*{2}}
\drawline[1000](50,155)(87,155)
\drawline[1000](50,155)(68.25,175)(87,155)(68.25,135)(50,155)
\put(71.5,165){\circle*{2}}
\drawline[1000](68.25,175)(74.75,155)
\put(65.5,165){\circle*{2}}
\drawline[1000](65.5,165)(62.25,155)
\drawline[1000](65.5,165)(71.5,165)
\put(45,150){$123$}
\put(59.5,150){$456$}
\put(72,150){$178$}
\put(85,150){$234$}
\put(68,177){$567$}
\put(68,130){$567$}
\put(73,162){$234$}
\put(58.5,162){$178$}

\put(100,175){$R_2$}
\multiput(100,155)(12.5,0){4}{\circle*{2}}
\multiput(118.25,175)(0,-40){2}{\circle*{2}}
\drawline[1000](100,155)(137,155)
\drawline[1000](100,155)(118.25,175)(137,155)(118.25,135)(100,155)
\put(121.5,165){\circle*{2}}
\drawline[1000](118.25,175)(124.75,155)
\put(115.5,165){\circle*{2}}
\drawline[1000](118.25,135)(115.5,165)(121.5,165)
\put(95,150){$123$}
\put(109.5,150){$456$}
\put(122,150){$178$}
\put(135,150){$234$}
\put(118,177){$678$}
\put(118,130){$578$}
\put(123,162){$235$}
\put(108.5,162){$146$}

\put(0,115){$R_3$}
\multiput(0,95)(12.5,0){4}{\circle*{2}}
\multiput(18.25,115)(0,-40){2}{\circle*{2}}
\drawline[1000](0,95)(37,95)
\drawline[1000](0,95)(18.25,115)(37,95)(18.25,75)(0,95)
\put(21.3625,105){\circle*{2}}
\drawline[1000](21.3625,105)(18.25,75)
\drawline[1000](18.25,115)(24.75,95)
\put(-5,90){$123$}
\put(9.5,90){$456$}
\put(22,90){$178$}
\put(35,90){$234$}
\put(18,117){$567$}
\put(18,70){$567$}
\put(14,105){$234$}

\put(50,115){$R_4$}
\multiput(50,95)(12.5,0){4}{\circle*{2}}
\multiput(68.25,115)(0,-40){2}{\circle*{2}}
\drawline[1000](50,95)(87,95)
\drawline[1000](50,95)(68.25,115)(87,95)(68.25,75)(50,95)
\put(71.5,105){\circle*{2}}
\drawline[1000](68.25,75)(71.5,105)(68.25,115)
\put(65.5,105){\circle*{2}}
\drawline[1000](62,95)(65.5,105)(71.5,105)
\put(45,90){$123$}
\put(59.5,90){$456$}
\put(72,90){$178$}
\put(85,90){$234$}
\put(68,117){$567$}
\put(68,70){$567$}
\put(73,102){$348$}
\put(58.5,102){$127$}

\put(100,115){$R_5$}
\multiput(100,95)(12.5,0){4}{\circle*{2}}
\multiput(118.25,115)(0,-40){2}{\circle*{2}}
\drawline[1000](100,95)(137,95)
\drawline[1000](100,95)(118.25,115)(137,95)(118.25,75)(100,95)
\put(121.5,105){\circle*{2}}
\drawline[1000](118.25,75)(115.5,105)
\put(115.5,105){\circle*{2}}
\drawline[1000](112,95)(115.5,105)(121.5,105)
\drawline[1000](118.25,115)(125,95)
\put(95,90){$123$}
\put(109.5,90){$456$}
\put(122,90){$378$}
\put(135,90){$124$}
\put(118,117){$578$}
\put(118,70){$568$}
\put(123,102){$246$}
\put(108.5,102){$137$}

\put(0,55){$R_6$}
\multiput(0,35)(12.5,0){4}{\circle*{2}}
\multiput(18.25,55)(0,-40){2}{\circle*{2}}
\drawline[1000](0,35)(37,35)
\drawline[1000](0,35)(18.25,55)(37,35)(18.25,15)(0,35)
\put(21.3625,45){\circle*{2}}
\drawline[1000](12.75,35)(18.25,15)
\put(15.1375,25){\circle*{2}}
\drawline[1000](18.25,55)(24.75,35)
\drawline[1000](15.1375,25)(21.3625,45)
\put(-5,30){$123$}
\put(7.5,30){$456$}
\put(22,30){$378$}
\put(35,30){$124$}
\put(18,57){$578$}
\put(18,10){$568$}
\put(14,45){$246$}
\put(16.5,22){$137$}

\put(50,55){$R_7$}
\multiput(50,35)(12.5,0){4}{\circle*{2}}
\multiput(68.25,55)(0,-40){2}{\circle*{2}}
\drawline[1000](50,35)(87,35)
\drawline[1000](50,35)(68.25,55)(87,35)(68.25,15)(50,35)
\put(71.3625,45){\circle*{2}}
\drawline[1000](62.75,35)(68.25,15)
\put(65.1375,25){\circle*{2}}
\drawline[1000](68.25,55)(74.75,35)
\put(45,30){$123$}
\put(57.5,30){$456$}
\put(72,30){$378$}
\put(85,30){$124$}
\put(68,57){$567$}
\put(68,10){$568$}
\put(64,45){$124$}
\put(66.5,22){$123$}

\end{picture}
\caption{Graphs $R_i$ for $0 \leq i \leq 7$ and $(8:3)$-colorings for them.}    \label{R's}
\end{figure}
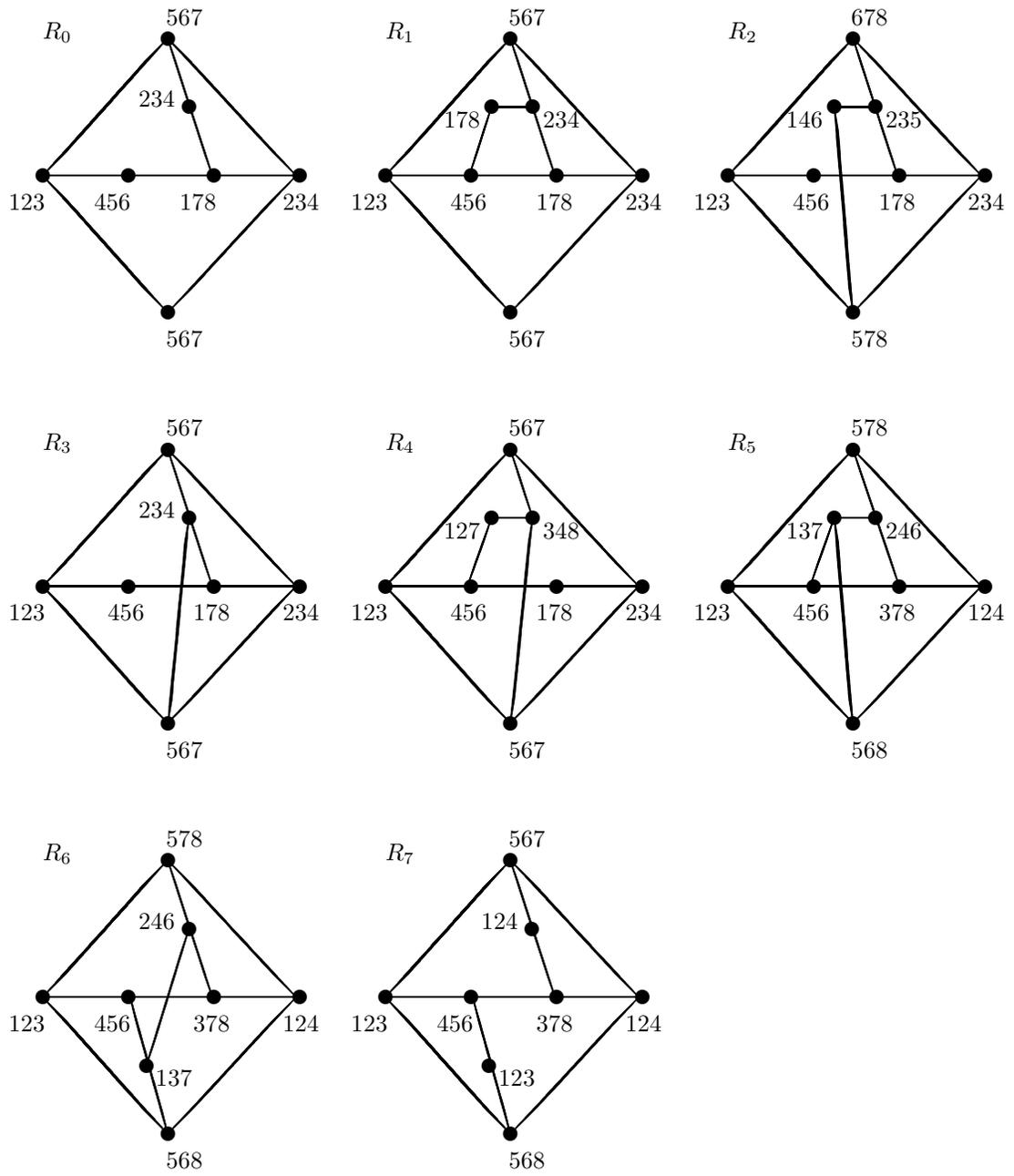

Now, we are ready to find a partial list of forbidden subgraphs.
We define $R_i$ to be the graphs in Figure \ref{R's} for $0 \leq i \leq 7$.
The following lemma excludes $R_1, R_2, ..., R_7$ as subgraphs from a fractionally $t$-critical graph.
But we cannot exclude $R_0$.

\begin{lemma}  \label{R_i-free}
	Every fractionally $t$-critical subcubic graph with $t \geq 8/3$ is $\{R_i: 1 \leq i \leq 7\}$-free.
\end{lemma}

\begin{pf}
Suppose that $H$ is an induced subgraph of a fractionally $t$-critical graph $G$ with $t \geq 8/3$ isomorphic to $R_i$ for some $1 \leq i \leq 7$.
As shown in Figure \ref{R's}, $R_i$ has an $(8:3)$-coloring for every $1 \leq i \leq 7$, so $\lvert N(H) \rvert \geq 2$ by Lemma \ref{basic structure}.
Since $R_5$ and $R_6$ are cubic, and $R_3$ contains only one vertex of degree at most two, so $i \neq 3,5,6$.
Given $i=1,2,4,7$, let $u_i$ and $v_i$ be the two vertices in $N(H)$, and let $H'$ be the subgraph of $G$ induced by $V(H) \cup \{u_i,v_i\}$.
Note that the pair of vertices of degree two in $H$ receive different and non-disjoint sets of colors in the $(8:3)$-colorings of $H$ shown in Figure \ref{R's}.
By Lemma \ref{extend 8:3}, $H'+u_iv_i$ and $H'/u_iv_i$ have $(8:3)$-colorings for every $i=1,2,4,7$.
Hence, every fractionally $t$-critical graph is $R_i$-free by Lemma \ref{2 cut set}, for $i=1,2,4,7$.
\end{pf}

\bigskip

Next, we investigate structure of fractionally $t$-critical graphs that contain $L_0$ as a subgraph, where $t \geq 8/3$, to obtain another list of forbidden subgraphs.
Note that two different copies of $L_0$ are not necessarily vertex-disjoint.

We denote the path on two vertices by $P_2$ and denote the cycle on four vertices by $C_4$.
For $i \in {\mathbb N}$, define $\L_0 = \{L_0\}$ and define $\L_i$ to be the collection of triangle-free subcubic graphs $H$ that can be obtained from some graph $H'$ in $\L_{i-1}$ by one of the following operations:
\begin{itemize}
	\item ({\bf Operation 1}) adding a disjoint $P_2$ to $H'$ and two edges, where one edge is incident with one end of the $P_2$ and one vertex in a $C_4$ in $H'$, and the other edge is incident with the other end of the $P_2$ and the diagonal vertex of the same $C_4$ in $H'$, or
	\item ({\bf Operation 2}) adding a disjoint $C_4$ to $H'$ and two edges, where one edge is incident with one vertex in the $C_4$ and one end of a $P_2$ in $H'$, and the other edge is incident with the diagonal vertex in the $C_4$ and the other end of the $P_2$.
\end{itemize}

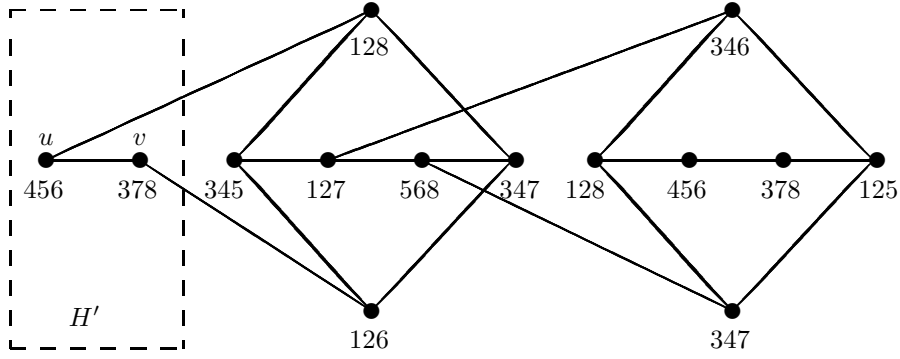
\begin{figure}
\unitlength=1mm
\begin{picture}(120, 30)(7,0)

\multiput(17,25)(12.5,0){2}{\circle*{2}}
\put(14,20){$456$}
\put(26.5,20){$378$}
\put(16,27){$u$}
\put(28.5,27){$v$}

\put(12.25,0){\dashbox{2}(23,45)}
\put(20,3){$H'$}

\multiput(42,25)(12.5,0){4}{\circle*{2}}
\multiput(60.25,45)(0,-40){2}{\circle*{2}}
\drawline[1000](42,25)(79,25)
\drawline[1000](42,25)(60.25,45)(79,25)(60.25,5)(42,25)
\put(38,20){$345$}
\put(51.5,20){$127$}
\put(64,20){$568$}
\put(77.3,20){$347$}
\put(57.25,39){$128$}
\put(57.25,0){$126$}

\multiput(90,25)(12.5,0){4}{\circle*{2}}
\multiput(108.25,45)(0,-40){2}{\circle*{2}}
\drawline[1000](90,25)(127,25)
\drawline[1000](90,25)(108.25,45)(127,25)(108.25,5)(90,25)
\put(86,20){$128$}
\put(99.5,20){$456$}
\put(112,20){$378$}
\put(125,20){$125$}
\put(105.25,39){$346$}
\put(105.25,0){$347$}

\drawline[1000](60.25,45)(17,25)(29.5,25)(60.25,5)
\drawline[1000](108.25,45)(54.5,25)(67,25)(108.25,5)

\end{picture}
\caption{The box surrounded by the dotted line denotes a graph $H'$ in ${\mathcal L}_{k-4}$, where $u,v$ is a pair of adjacent vertices in $H'$. The whole graph $H$ is in ${\mathcal L}_k$ obtained from $H'$ by consecutively doing Operations 2,1,2,1, where the first Operation 2 adding a new $C_4$ and two edges incident with a $P_2$ $uv$ in $H'$. An $(8:3)$-coloring of a subgraph of $H$ is presented.}     \label{4 operations}
\end{figure}

Figure \ref{4 operations} shows an example for consecutive applying the above operations four times.
Observe that every graph in ${\mathcal L}_i$, where $i \geq 0$, contains exactly four vertices of degree two, and these four vertices can be paired such that each pair either induces a path on two vertices or consists of two diagonal vertices of a $4$-cycle.
For every finite sequence $(a_1, a_2, ..., a_k)$ with $a_i \in \{1,2\}$ for all $1 \leq i \leq k$, we denote $L_{a_1, a_2 ,..., a_k}$ as the graph that is obtained from $L_0$ by doing Operations $a_1, a_2, ..., a_k$ consecutively.
Note that given a positive integer $t$ and $j=1,2$, for every graph $H \in \L_t$ and for each Operation $j$, there is at most one way (up to isomorphism) to add edges to link $H$ and the new $P_2$ or $C_4$, so $L_{a_1, a_2 ,..., a_k}$ is well-defined.
However, $L_{a_1, a_2 ,..., a_k}$ may be equal to $L_{b_1, b_2 ,..., b_m}$ for two different sequences $(a_1, a_2 ,..., a_k)$ and $(b_1, b_2, ..., b_m)$.
The following lemma ensures that every graph in $\L_t$ can be generated by Operations 1 and 2 alternately.

\begin{lemma} \label{normalize L sequence graph}
Every graph $H$ in $\L_t$, there is a sequence $(a_1, a_2, ..., a_t)$ with $a_1 \in \{1,2\}$ and $a_{i+1} = 3-a_i$ for $1 \leq i \leq t-1$ such that $H = L_{a_1, a_2,...,a_t}$.
\end{lemma}

\begin{pf}
Let $(b_1,b_2,...,b_t)$ be an $1$-$2$ sequence such that $H=L_{b_1,b_2,...,b_t}$.
Suppose that $i$ is the smallest index such that $b_i=b_{i+1}$.
Clearly, $i \neq 1$.
Then it is easy to see that $H=L_{b_{i+1},b_1,b_2,...,b_i, b_{i+2},...,b_t}$.
And the lemma follows from repeating this process.
\end{pf}

\bigskip

We denote $L_{a_1, a_2,...,a_t}$ by $L^{(a_1)}_t$ for each sequence $(a_1, a_2,...,a_t)$ with $a_i \in \{1,2\}$ and $a_{i+1} = 3-a_i$.
By Lemma \ref{normalize L sequence graph}, $\L_t = \{L^{(1)}_t, L^{(2)}_t\}$ for each positive integer $t$.
Let $\L'$ be the set of triangle-free subcubic graphs $H'$ for which there exist an integer $i$ and a graph $H$ in $\L_i$ with $\lvert E(H') \rvert > \lvert E(H) \rvert$ such that $H'$ contains $H$ as a spanning subgraph.
Observe that every graph in $\L'$ either is a cubic graph or contains exactly two vertices of degree two.

\begin{lemma}  \label{L'-free}
Every graph in $\L'$ has an $(8:3)$-coloring such that if it is not cubic, then the pair of the degree-two vertices receive different non-disjoint sets of colors.
Furthermore, every fractionally $t$-critical triangle-free subcubic graph with $t \geq 8/3$ is $\L'$-free.
\end{lemma}

\begin{figure}
\unitlength=0.8mm
\begin{picture}(120, 200)(15,-60)

\put(0,175){$L_1$}
\multiput(0,155)(12.5,0){4}{\circle*{2}}
\multiput(18.25,175)(0,-40){2}{\circle*{2}}
\drawline[1000](0,155)(37,155)
\drawline[1000](0,155)(18.25,175)(37,155)(18.25,135)(0,155)
\put(11,157){{\scriptsize$w$}}
\put(23.5,157){{\scriptsize$x$}}
\put(-3,150){{\scriptsize$123$}}
\put(9.5,150){{\scriptsize$456$}}
\put(22,150){{\scriptsize$378$}}
\put(37,150){{\scriptsize$124$}}
\put(16,177){{\scriptsize$567$}}
\put(16,130){{\scriptsize$678(568)$}}

\multiput(45,165)(0,-20){2}{\circle*{2}}
\drawline[1000](18.25,175)(45,165)(45,145)(18.25,135)
\put(47,165){{\scriptsize$y$}}
\put(47,145){{\scriptsize$z$}}
\put(41,170){{\scriptsize$238(128)$}}
\put(41,140){{\scriptsize$145(347)$}}

\put(70,175){$L_2$}
\multiput(70,155)(12.5,0){4}{\circle*{2}}
\multiput(88.25,175)(0,-40){2}{\circle*{2}}
\drawline[1000](70,155)(107,155)
\drawline[1000](70,155)(88.25,175)(107,155)(88.25,135)(70,155)
\put(87,171){{\scriptsize$w$}}
\put(87,137){{\scriptsize$x$}}
\put(67,150){{\scriptsize$123$}}
\put(79.5,150){{\scriptsize$456$}}
\put(92,150){{\scriptsize$378$}}
\put(107,150){{\scriptsize$124$}}
\put(86,177){{\scriptsize$567$}}
\put(86,130){{\scriptsize$578$}}

\multiput(120,155)(37.5,0){2}{\circle*{2}}
\multiput(138.75,175)(0,-40){2}{\circle*{2}}
\drawline[1000](120,155)(138.75,175)(157.5,155)(138.75,135)(120.55,155)
\drawline[1000](82.5,155)(138.75,175)
\drawline[1000](95,155)(138.75,135)
\put(118.5,157){{\scriptsize$y$}}
\put(156,157){{\scriptsize$z$}}
\put(136,177){{\scriptsize$127$}}
\put(136,130){{\scriptsize$125$}}
\put(117,150){{\scriptsize$348$}}
\put(154.75,150){{\scriptsize$346(348)$}}

\put(0,121){$L_{1,2} \cong L_{2,1}$}
\multiput(0,95)(12.5,0){4}{\circle*{2}}
\multiput(18.25,115)(0,-40){2}{\circle*{2}}
\drawline[1000](0,95)(37,95)
\drawline[1000](0,95)(18.25,115)(37,95)(18.25,75)(0,95)
\put(11,97){{\scriptsize$w$}}
\put(23.5,97){{\scriptsize$x$}}
\put(-3,90){{\scriptsize$123$}}
\put(9.5,90){{\scriptsize$456$}}
\put(22,90){{\scriptsize$378$}}
\put(37,90){{\scriptsize$124$}}
\put(16,117){{\scriptsize$567$}}
\put(16,70){{\scriptsize$678(568)$}}

\multiput(45,105)(0,-20){2}{\circle*{2}}
\drawline[1000](18.25,115)(45,105)(45,85)(18.25,75)
\put(39,109){{\scriptsize$238(128)$}}
\put(39,79.5){{\scriptsize$145(347)$}}

\multiput(52.5,95)(37.5,0){2}{\circle*{2}}
\multiput(70.75,115)(0,-40){2}{\circle*{2}}
\drawline[1000](52.5,95)(70.75,115)(90,95)(70.75,75)(52.5,95)
\drawline[1000](45,105)(70.75,115)
\drawline[1000](45,85)(70.75,75)
\put(51,98){{\scriptsize$y$}}
\put(89,98){{\scriptsize$z$}}
\put(49.5,90.5){{\scriptsize$123$}}
\put(48.5,87.5){{\scriptsize$(138)$}}
\put(89,90.5){{\scriptsize$124$}}
\put(87.5,87.5){{\scriptsize$(148)$}}
\put(68,117){{\scriptsize$567$}}
\put(68,70){{\scriptsize$678(256)$}}

\put(100,121){$L_{1,2,1}$}
\multiput(100,95)(12.5,0){4}{\circle*{2}}
\multiput(118.25,115)(0,-40){2}{\circle*{2}}
\drawline[1000](100,95)(137,95)
\drawline[1000](100,95)(118.25,115)(137,95)(118.25,75)(100,95)
\put(111,97){{\scriptsize$w$}}
\put(123.5,97){{\scriptsize$x$}}
\put(97,90){{\scriptsize$123$}}
\put(109.5,90){{\scriptsize$456$}}
\put(122,90){{\scriptsize$378$}}
\put(137,90){{\scriptsize$124$}}
\put(116,117){{\scriptsize$567$}}
\put(116,70){{\scriptsize$568$}}

\multiput(145,105)(0,-20){2}{\circle*{2}}
\drawline[1000](118.25,115)(145,105)(145,85)(118.25,75)
\put(142,109){{\scriptsize$128$}}
\put(142,80){{\scriptsize$347$}}

\multiput(152.5,95)(12.5,0){4}{\circle*{2}}
\multiput(170.75,115)(0,-40){2}{\circle*{2}}
\drawline[1000](152.5,95)(190,95)
\drawline[1000](152.5,95)(170.75,115)(190,95)(170.75,75)(152.5,95)
\drawline[1000](145,105)(170.75,115)
\drawline[1000](145,85)(170.75,75)
\put(163.5,98){{\scriptsize$y$}}
\put(176.5,98){{\scriptsize$z$}}
\put(149.5,90.5){{\scriptsize$124$}}
\put(148.5,87.5){{\scriptsize$(247)$}}
\put(162,90.5){{\scriptsize$378$}}
\put(161,87.5){{\scriptsize$(138)$}}
\put(174.5,90.5){{\scriptsize$456$}}
\put(173.5,87.5){{\scriptsize$(567)$}}
\put(189,90.5){{\scriptsize$123$}}
\put(188,87.5){{\scriptsize$(248)$}}
\put(168,117){{\scriptsize$567(356)$}}
\put(168,70){{\scriptsize$568(156)$}}

\put(0,60){$L_{2,1,2}$}
\multiput(0,35)(12.5,0){4}{\circle*{2}}
\multiput(18.25,55)(0,-40){2}{\circle*{2}}
\drawline[1000](0,35)(37,35)
\drawline[1000](0,35)(18.25,55)(37,35)(18.25,15)(0,35)
\put(17,51){{\scriptsize$w$}}
\put(17,17){{\scriptsize$x$}}
\put(-3,30){{\scriptsize$123$}}
\put(9.5,30){{\scriptsize$456$}}
\put(22,30){{\scriptsize$378$}}
\put(37,30){{\scriptsize$124$}}
\put(16,57){{\scriptsize$567$}}
\put(16,10){{\scriptsize$578$}}

\multiput(50,35)(12.5,0){4}{\circle*{2}}
\multiput(68.75,55)(0,-40){2}{\circle*{2}}
\drawline[1000](50,35)(87.5,35)
\drawline[1000](50,35)(68.75,55)(87.5,35)(68.75,15)(50.55,35)
\drawline[1000](12.5,35)(68.75,55)
\drawline[1000](25,35)(68.75,15)
\put(66,57){{\scriptsize$127$}}
\put(66,10){{\scriptsize$125$}}
\put(47,30){{\scriptsize$348$}}
\put(59.5,30){{\scriptsize$167$}}
\put(72,30){{\scriptsize$258$}}
\put(86.5,30){{\scriptsize$346$}}

\multiput(100,35)(37.5,0){2}{\circle*{2}}
\multiput(118.75,55)(0,-40){2}{\circle*{2}}
\drawline[1000](100,35)(118.75,55)(137.5,35)(118.75,15)(100.55,35)
\drawline[1000](62.5,35)(118.75,55)
\drawline[1000](75,35)(118.75,15)
\put(98.5,37){{\scriptsize$y$}}
\put(136,37){{\scriptsize$z$}}
\put(116,57){{\scriptsize$345$}}
\put(116,10){{\scriptsize$347$}}
\put(97,30){{\scriptsize$128$}}
\put(136,30){{\scriptsize$126(128)$}}

\put(0,1){$L_{1,2,1,2} \cong L_{2,1,2,1}$}
\multiput(0,-25)(12.5,0){4}{\circle*{2}}
\multiput(18.25,-5)(0,-40){2}{\circle*{2}}
\drawline[1000](0,-25)(37.5,-25)
\drawline[1000](0,-25)(18.25,-5)(37,-25)(18.25,-45)(0,-25)
\put(11,-23){{\scriptsize$w$}}
\put(23.5,-23){{\scriptsize$x$}}
\put(-3,-30){{\scriptsize$123$}}
\put(9.5,-30){{\scriptsize$456$}}
\put(22,-30){{\scriptsize$378$}}
\put(37,-30){{\scriptsize$124$}}
\put(16,-3){{\scriptsize$567$}}
\put(16,-50){{\scriptsize$678(568)$}}

\multiput(45,-15)(0,-20){2}{\circle*{2}}
\drawline[1000](18.25,-5)(45,-15)(45,-35)(18.25,-45)
\put(39,-11){{\scriptsize$138(128)$}}
\put(39,-41){{\scriptsize$245(347)$}}

\multiput(52.5,-25)(12.5,0){4}{\circle*{2}}
\multiput(70.75,-5)(0,-40){2}{\circle*{2}}
\drawline[1000](52.5,-25)(90,-25)
\drawline[1000](52.5,-25)(70.75,-5)(90,-25)(70.75,-45)(52.5,-25)
\drawline[1000](45,-15)(70.75,-5)
\drawline[1000](45,-35)(70.75,-45)
\put(49.5,-29.5){{\scriptsize$148$}}
\put(48.5,-32.5){{\scriptsize$(124)$}}
\put(62,-29.5){{\scriptsize$237$}}
\put(61,-32.5){{\scriptsize$(356)$}}
\put(74.5,-29.5){{\scriptsize$156$}}
\put(73.5,-32.5){{\scriptsize$(478)$}}
\put(91.75,-27.5){{\scriptsize$248$}}
\put(90.75,-30.5){{\scriptsize$(123)$}}
\put(68,-3){{\scriptsize$567$}}
\put(68,-50){{\scriptsize$367(568)$}}

\multiput(102.5,-25)(37.5,0){2}{\circle*{2}}
\multiput(120.75,-5)(0,-40){2}{\circle*{2}}
\drawline[1000](102.5,-25)(120.75,-5)(140,-25)(120.75,-45)(102.5,-25)
\drawline[1000](65,-25)(120.75,-5)
\drawline[1000](77.5,-25)(120.75,-45)
\put(101.5,-22){{\scriptsize$y$}}
\put(139,-22){{\scriptsize$z$}}
\put(106.75,-25){{\scriptsize$123$}}
\put(105.75,-28){{\scriptsize$(378)$}}
\put(143,-25){{\scriptsize$125$}}
\put(142,-28){{\scriptsize$(357)$}}
\put(118,-3){{\scriptsize$468(124)$}}
\put(118,-50){{\scriptsize$478(126)$}}

\end{picture}
\caption{$(8:3)$-colorings $f_1$ of graphs belonging to $\L_k$, for $1 \leq k \leq 4$, such that $f_1(w) \cap f_1(y) = f_1(x) \cap f_1(z) = \emptyset$. Let $f_2$ be the $(8:3)$-coloring obtained from $f_1$ by replacing the colors of each vertex by the colors in the brackets. Then $f_2(w) \cap f_2(y) = \emptyset$ and $1 \leq \lvert f_2(x) \cap f_2(z) \rvert \leq 2$.}    \label{(8:3)-coloring of L's}
\end{figure}
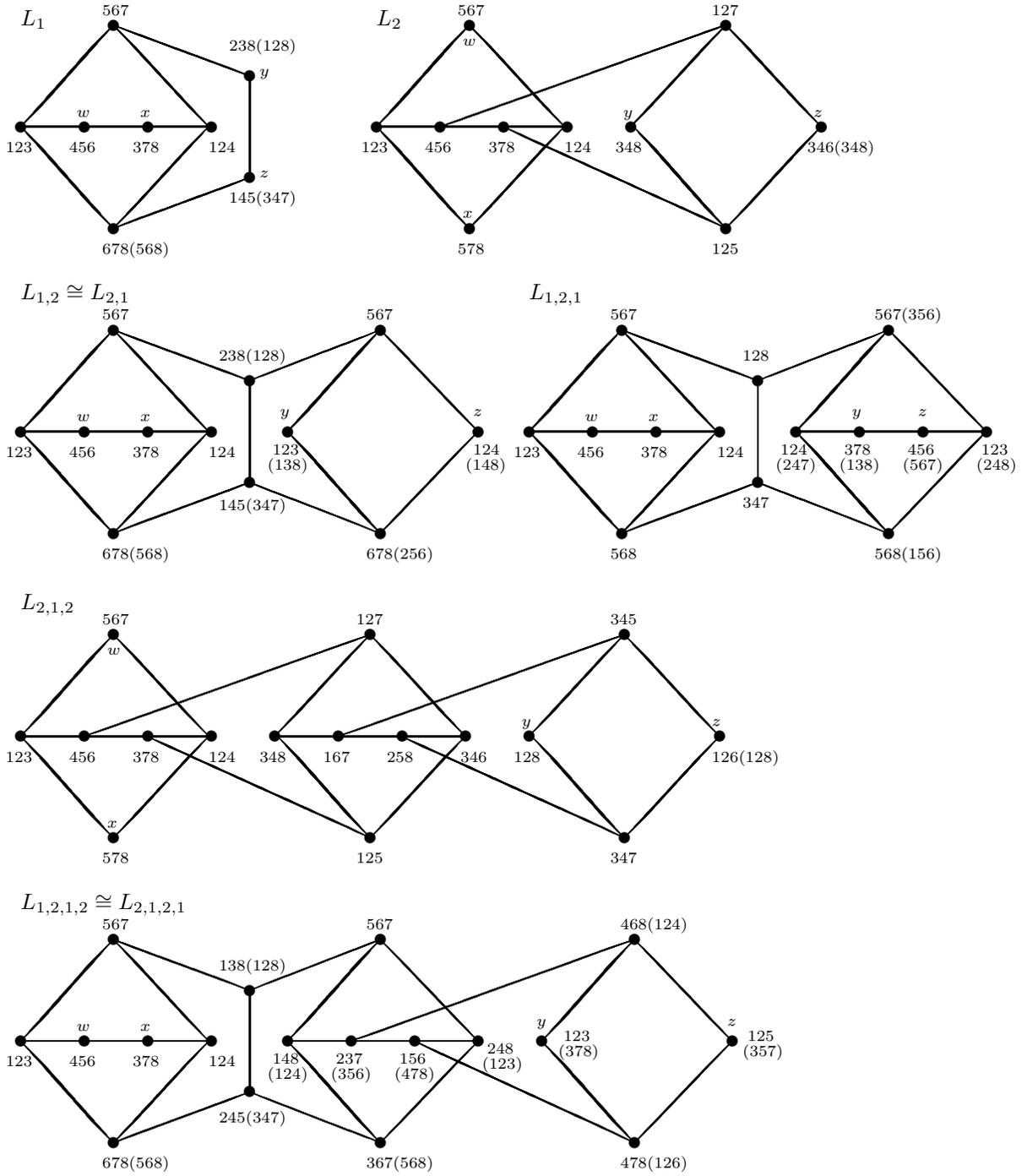

\begin{pf}
Let $k$ be a nonnegative integer, and let $H$ be a graph in $\L_k$.
Let $H'$ be in $\L'$ containing $H$ as a spanning subgraph with $\lvert E(H') \rvert > \lvert E(H) \rvert$.
Denote the four vertices of degree two in $H$ by $\{w,x,y,z\}$.
By symmetry, we may assume one of the following holds: 
	\begin{enumerate}
		\item Both $wx$ and $yz$ are edges in $H$.
		\item $wx$ is an edge in $H$ and $y,z$ are diagonal vertices in a $4$-cycle.
		\item $w,x$ are diagonal vertices in a $4$-cycle and $y,z$ are diagonal vertices in another $4$-cycle.
	\end{enumerate}
To show that $H'$ has an $(8:3)$-coloring, by Lemma \ref{extend 8:3}, it is sufficient to show the following.
The coloring $f_1$ deals with the case that $H'$ is cubic, and $f_2$ deals with the case that $H'$ is not cubic. 
	\begin{enumerate}
		\item If $wx$ and $yz$ are edges, then there are $(8:3)$-colorings $f_1, f_2$ of $H$ such that $f_1(w) \cap f_1(y) = \emptyset$, $f_1(x) \cap f_1(z) = \emptyset$, $f_2(w) \cap f_2(y) = \emptyset$, and $1 \leq \lvert f_2(x) \cap f_2(z) \rvert \leq 2$.
		\item If $wx$ is an edge and $y,z$ are diagonal vertices in a $4$-cycle, then there are $(8:3)$-colorings $f_1, f_2$ of $H$ such that $f_1(w) \cap f_1(y) = \emptyset$, $f_1(x) \cap f_1(z) = \emptyset$, $f_2(w) \cap f_2(y) = \emptyset$, $1 \leq \lvert f_2(x) \cap f_2(z) \rvert \leq 2$.
		\item If $w,x$ are diagonal vertices in a $4$-cycle and $y,z$ are diagonal vertices in another $4$-cycle, then there are $(8:3)$-colorings $f_1, f_2$ of $H$ such that $f_1(w) \cap f_1(y) = \emptyset$, $f_1(x) \cap f_1(z) = \emptyset$, $f_2(w) \cap f_1(y) = \emptyset$, $1 \leq \lvert f_2(w) \cap f_2(y) \rvert \leq 2$.
	\end{enumerate}

The coloring given in Figure \ref{4 operations} shows that the patterns of the colors on the pair of adjacent degree-two vertices can be kept if we do Operations 2,1,2,1 consecutively.
Hence, it is sufficient to check the above three cases for $1 \leq k \leq 4$.
Desired $f_1$ and $f_2$ are shown in Figure \ref{(8:3)-coloring of L's}.
Note that the case $k=0$ is excluded by the triangle-freeness of $H'$.

Let $G$ be a fractionally $t$-critical triangle-free subcubic graph with $t \geq 8/3$.
Suppose that $G$ contains $H'$ as an induced subgraph.
If $H'$ is cubic, then $G=H'$ has an $(8:3)$-coloring, a contradiction.
If $H'$ is not cubic, then there exist a pair of vertices $u,v$ such that $H'$ is a component of $G-\{u,v\}$.
By the nice property of $f_2$ and Lemmas \ref{2 cut set} and \ref{extend 8:3}, $\chi(G)\leq 8/3$, a contradiction.
\end{pf}

\bigskip

The rest of this section is dedicated to proving Lemma \ref{good 3-coloring}.

\begin{lemma} \label{intersecting L_0}
Let $G$ be a $\{K_3, R_i: 1 \leq i \leq 7\} \cup \L'$-free subcubic graph.
Let $t$ be a positive integer, and for each $i$ with $1 \leq i \leq t$, let $H_i$ be a subgraph of $G$ such that $H_i$ is isomorphic to $L_0$.
If $V(H_{i+1}) \cap (\bigcup_{j=1}^i V(H_j)) \neq \emptyset$ for all $i \geq 1$, then $G[\bigcup_{i=1}^t V(H_i)]$ is isomorphic to $R_0$ or a graph in $\L_k$ for some $k \geq 0$.
\end{lemma}

\begin{pf}
For each $1 \leq i \leq t$, let $V(H_i) = \{v_{i,1}, v_{i,2}, v_{i,3}, v_{i,4}, u_{j,1}, u_{j,2}\}$, where $(v_{i,1}, v_{i,2}, v_{i,3}, v_{i,4}, u_{i,1}, u_{i,2})$ forms a copy of $L_0$.
We shall prove this lemma by induction on $t$.
The lemma is clear if $t=1$.
Now, we prove the case that $t=2$.

Suppose this lemma is not true when $t=2$.
If $(v_{2,1},v_{2,2}) = (u_{1,1}, v_{1,1})$, then $G[V(H_1) \cup V(H_2)]$ is isomorphic to $R_0$ by the triangle-freeness.
So we may assume that $(v_{2,1},v_{2,2}) \neq (u_{1,1}, v_{1,1})$.
If $v_{2,1} = u_{1,1}$, then $v_{2,2} \not \in V(H_1)$ by symmetry, and hence $G[V(H_1) \cup V(H_2)]$ is isomorphic to the graph that can be obtained from $L_0$ by Operation 1.
So we may assume that $\{v_{2,1},v_{2,4}\} \cap \{u_{1,1},u_{1,2}\} = \emptyset$.
If $v_{2,2} = u_{1,1}$ and $v_{2,1} \in V(H_1)$, then $G[V(H_1) \cup V(H_2)]$ is isomorphic to the graph that can be obtained from $L_0$ by Operation 1 and then adding an edge, so it is in ${\mathcal L}'$.
If $v_{2,2} = u_{1,1}$ and $v_{2,1} \not \in V(H_1)$, then $G[V(H_1) \cup V(H_2)]$ is isomorphic to $R_1$, a contradiction.
Hence, we may assume that $\{u_{1,1}, u_{1,2}\} \cap \{v_{2,j}: 1 \leq j \leq 4\} = \{u_{2,1}, u_{2,2}\} \cap \{v_{1,j}: 1 \leq j \leq 4\} = \emptyset$.
Consequently, if $v_{2,1}$ or $v_{2,4}$ is in $V(H_1)$, say $v_{2,1} \in V(H_1)$, then $v_{2,1}$ is equal to $v_{1,1}$ or $v_{1,4}$, and hence $G[V(H_1) \cup V(H_2)]$ is isomorphic to $L_0$, $R_3$ or $R_4$.
So $v_{2,1}$ and $v_{2,4}$ are not in $V(H_1)$.
Furthermore, if $v_{2,2}$ or $v_{2,3}$ is in $V(H_1)$, say $v_{2,2} \in V(H_1)$, then $v_{2,2}$ is $v_{1,2}$ or $v_{1,3}$, and hence $G[V(H_1) \cup V(H_2)]$ is isomorphic the graph that can be obtained from $L_0$ by Operation 2.
As a result, $\{v_{2,j}: 1 \leq j \leq 4\}$ is disjoint from $V(H_1)$, and it implies that $\{u_{2,1},u_{2,2}\} \cap V(H_1) = \emptyset$, so $V(H_1) \cap V(H_2) = \emptyset$, a contradiction.
This proves the case that $t=2$.

Let $t \geq 3$.
First, suppose that $G[\cup_{i=1}^{t-1}V(H_i)] = R_0$.
If $\{v_{t,1},v_{t,4}\} \cap (\bigcup_{i=1}^{t-1}V(H_i)) \neq \emptyset$, then $G[\cup_{i=1}^t V(H_i)]$ is isomorphic to $R_0$, $R_1$, $R_2$, $R_3$, $R_5$, $R_6$ or $R_7$, a contradiction.
So $\{v_{t,1},v_{t,4}\} \cap (\bigcup_{i=1}^{t-1}V(H_i)) = \emptyset$, but it implies that $\{v_{t,2},v_{t,3}, u_{t,1}, u_{t,2}\} \cap (\bigcup_{i=1}^{t-1}V(H_i)) = \emptyset$ as well, a contradiction.
Therefore, $G[\cup_{i=1}^{t-1}H_i]$ is isomorphic to a graph in $\L_k$ for some $k \geq 0$ by the induction hypothesis.
Note that $k \geq 1$ since we may assume that $V(H_{i+1}) - (\bigcup_{j=1}^{i} V(H_j)) \neq \emptyset$ for $1 \leq i \leq t-1$, without loss of generality.
In this case, we have the following two claims since $G$ is $\L'$-free.

\noindent {\bf Claim 1:} Each vertex of degree two in $G[\bigcup_{j=1}^{t-1}V(H_j)]$ is adjacent to a vertex of degree three in $G[\bigcup_{j=1}^{t-1}V(H_j)]$.

\noindent {\bf Claim 2:} Each pair of non-adjacent vertices of degree two in $G[\bigcup_{j=1}^{t-1}V(H_j)]$ either are in a $4$-cycle or have distance at least three.
Furthermore, if $v$ is a vertex of degree two in $G[\bigcup_{j=1}^{t-1}V(H_j)]$, and $v$ is not adjacent to any vertex which has degree two in $G[\bigcup_{j=1}^{t-1}V(H_j)]$, then there exists uniquely a vertex $u$ having degree two in $G[\bigcup_{j=1}^{t-1}V(H_j)]$ such that the distance between $u$ and $v$ is two, and all other degree two vertices in $G[\bigcup_{j=1}^{t-1}V(H_j)]$ have distance at least four from $v$ in $G[\bigcup_{j=1}^{t-1}V(H_j)]$.

Observe that $\{v_{t,j}: 1 \leq j \leq 4\} \cap \bigcup_{j=1}^{t-1}V(H_j) \neq \emptyset$, since every graph in $\L_k$ has minimum degree two.
So it is impossible that $v_{t,1}, v_{t,2} \not \in \bigcup_{j=1}^{t-1}V(H_j)$ by Claims 1 and 2.
Assume that $v_{t,1} \not \in \bigcup_{j=1}^{t-1}V(H_j)$, so $v_{t,2} \in \bigcup_{j=1}^{t-1}V(H_j)$, and $u_{t,1}$ and $u_{t,2}$ do not have degree three in $G[\bigcup_{j=1}^{t-1}V(H_j)]$.
By Claims 1 and 2, $v_{t,4} \not \in \bigcup_{j=1}^{t-1}V(H_j)$, but this implies that $G[\bigcup_{j=1}^t V(H_j)] \in \L_{k+1}$.
Hence, we may assume that $v_{t,1}, v_{t,4} \in \bigcup_{j=1}^{t-1}V(H_j)$ by symmetry.
Therefore, one of $u_{t,1}$ and $u_{t,2}$ is in $G[\bigcup_{j=1}^{t-1}V(H_j)]$ as well, so $v_{t,1}$ and $v_{t,4}$ have distance two in $G[\bigcup_{j=1}^{t-1}V(H_j)]$.
It implies that either $v_{t,1}$, $v_{t,4}$ are diagonal vertices of a $4$-cycle, or one of $v_{t,1}$ and $v_{t,4}$ has degree three in $\bigcup_{j=1}^{t-1}V(H_j)$ by Claim 2.
In the latter case, $u_{t,1}$ and $u_{t,2}$ are in $G[\bigcup_{j=1}^{t-1}V(H_j)]$, and $v_{t,1}$ and $v_{t,2}$ are diagonal vertices of a $4$-cycle.
Consequently, $G[\bigcup_{j=1}^{t}V(H_j)]$ is isomorphic to a graph in $\L_{k+1}$ if $H_t \not \subseteq G[\bigcup_{j=1}^t V(H_j)]$.
This completes the proof.
\end{pf}

\bigskip

Let $G$ be a graph and $H$ a subgraph of $G$.
We say that a proper $3$-coloring $f$ of $G$ is {\it $H$-rainbow-free} if $H$ contains no rainbow copy of $L_0$ with respect to $f$.

\begin{lemma} \label{extend coloring to L_0's}
Let $G$ be a subcubic $\{K_3\} \cup {\mathcal L}'$-free graph and $H$ a subgraph of $G$ isomorphic to a graph in $\L_t$ for some $t \geq 0$.
Let $f$ be a proper 3-coloring of $G-H$.
If either there exists a vertex $v$ of $H$ such that $\deg_G(v) < 3$, or $\lvert f(N_G(H)) \rvert \geq 2$, then $f$ can be extended to a proper $H$-rainbow-free 3-coloring of $G$.
\end{lemma}

\begin{pf}
We shall prove this lemma by induction on $t$.
We may assume that $\deg_G(x)=3$ for all $x \in V(H)$; otherwise, we construct a graph $G'$ from $G$ by adding a vertex $u_v$ and an edge $vu_v$ for each vertex $v$ of $H$ of degree $2$ in $G$, and then define $f(u_v)$ to be an element in $\{1,2,3\}$ such that $\lvert f(N_{G'}(H)) \rvert$ is as large as possible.
This ensures that $\lvert f(N_{G'}(H)) \rvert \geq 2$, and any proper $H$-rainbow-free $3$-coloring $g$ of $G'$ extended from $f$ is a proper $H$-rainbow-free $3$-coloring of $G$ as well.
So it is sufficient to deal with the case that $\lvert f(N_G(H)) \rvert \geq 2$.

When $t=0$, $H$ is an induced subgraph as $G$ is triangle-free.
Let $V(H) = \{u_i, v_j: i=1,2, j=1,2,3,4\}$ and $(v_1, v_2, v_3, v_4, u_1, u_2)$ form a copy of $L_0$.
Let $x_i \in N_G(u_i) - V(H)$, and $y_j \in N_G(v_j)-V(H)$, for $i=1,2$ and $j=2,3$. 
Then $f$ can be extended to a proper $H$-rainbow-free 3-coloring $g$ of $G$ by first defining $g(v_1)=g(v_4)=f(y_2)$ and then defining $g$ on the remaining vertices in $H$.

Now, we assume that $t \geq 1$, so $H$ is isomorphic to $L^{(k)}_t$ for some $k=1,2$.
First, we assume that $k=2$.
So there is a cycle $abcda$ such that $H - \{a,b,c,d\}$ is isomorphic to a graph in $\L_{t-1}$, where $\deg_H(a) = \deg_H(c)=2$.
Let $a',c'$ be the neighbor of $a$ and $c$ not in $\{b,d\}$, respectively.
Let $\{x,y\} = N_G(H) - \{a',c'\}$.
Note that $a',c',x,y$ may not be pairwise distinct.
If $f(a') \neq f(c')$, say $f(a')=1$ and $f(c')=2$, then without loss of generality, we may assume that $(f(x), f(y)) \neq (1,1)$ (since we can swap colors $1$ and $2$), and hence we can extend $f$ to a proper $H$-rainbow-free 3-coloring $g$ of $G$ by first defining $(g(a),g(b),g(c),g(d)) = (3,1,3,1)$, and then defining $g$ on $V(H) - \{a,b,c,d\}$ by the induction hypothesis.
If $f(a')=f(c')$, say $f(a')=f(c')=1$, then without loss of generality, we may assume that $(f(x), f(y)) \neq (3,3)$ (since we can swap colors $2$ and $3$), and hence we can extend $f$ to a proper $H$-rainbow-free 3-coloring $g$ of $G$ by first defining $(g(a),g(b),g(c),g(d)) = (2,3,2,3)$, and then defining $g$ on $V(H) - \{a,b,c,d\}$ by the induction hypothesis.
This proves the case that $k=2$.

Assume that $k=1$.
So there is a tuple $(p_1, p_2, p_3, p_4, q_1, q_2)$ of vertices in $H$ such that $(p_1, p_2, p_3, p_4, q_1, q_2)$ forms a copy of $L_0$ with $\deg_H(p_2) \allowbreak =\deg_H(p_3)=2$, and $H - \{p_1,p_2,p_3,p_4,q_1,q_2\}$ is in $\L_{t-2}$ if $t \geq 2$, and a $P_2$ if $t=1$.
Let $p'_2, p'_3$ be the neighbor of $p_2,p_3$ not in $\{p_1,p_2,p_3,p_4\}$, respectively.
Let $\{x,y\} = N_G(H) - \{p_2',p_3'\}$.
Note that $p'_2, p'_3,x,y$ may not be pairwise distinct.

Suppose that $t=1$.
Let $\{w_x,w_y\}=V(H)-\{p_1,p_2,p_3,p_4,q_1,q_2\}$, where $w_x$ is adjacent to $q_1$ and $x$, and $v$ is adjacent to $q_2$ and $y$.
First, assume that $f(p_2') \neq f(p_3')$ and $f(x) \neq f(y)$.
By symmetry, we assume that $f(p_2')=f(x)=1$, $f(p_3')=2$, and $f(y)=3$.
Then $f$ can be extended to a proper $H$-rainbow-free $3$-coloring $g$ by defining $(g(p_1),g(p_2),g(p_3),g(p_4),g(q_1),g(q_2),g(w_x),g(w_y))$ to be $(1,3,1,2,3,3,2,1)$.
Second, we assume that $f(p_2') \neq f(p_3')$ and $f(x)=f(y)$.
If $f(x) \in \{f(p_2'),f(p_3')\}$, then we may assume that $f(p_2')=1$ and $f(x)=f(y)=f(p_3')=2$.
In this case, $f$ can be extended to a proper $H$-rainbow-free $3$-coloring $g$ by defining $(g(p_1),g(p_2),g(p_3),g(p_4),g(q_1),g(q_2),g(w_x),g(w_y))$ to be $(1,2,1,3,2,2,1,3)$.
If $f(x) \not \in \{f(p_2'),f(p_3')\}$, then we may assume that $f(p_2')=1$, $f(p_3')=2$ and $f(x)=f(y)=3$.
Then $f$ can be extended to a proper $H$-rainbow-free $3$-coloring $g$ by defining $(g(p_1),g(p_2),g(p_3), \allowbreak g(p_4),g(q_1),g(q_2),g(w_x),g(w_y))$ to be $(1,2,3,2,3,3,2,1)$.
Therefore, by symmetry, it is sufficient to consider the case that $f(p'_2)=f(p'_3)$ and $f(x)=f(y)$.
Without loss of generality, we may assume that $f(p'_2)=f(p'_3)=1$ and $f(x)=f(y)=2$.
Then it is clear that we can extend $f$ to a proper $H$-rainbow-free 3-coloring $g$ by defining $(g(p_1),g(p_2), g(p_3), \allowbreak g(p_4), g(q_1), \allowbreak g(q_2), g(w_x), g(w_y)) = (1,2, \allowbreak 3,1,2,2,3, \allowbreak 1)$.
This proves the case that $k=1$ and $t=1$.

It remains to deal with the case that $k=1$ and $t \geq 2$.
If $f(p'_2) \neq f(p'_3)$, say $f(p'_2)=1, f(p'_3)=2$, then without loss of generality, we may assume $(f(x),f(y)) \neq (1,1)$ (since we can swap colors $1$ and $2$), and hence we can extend $f$ to a proper $H$-rainbow-free 3-coloring $g$ of $G$ by first defining $(g(p_1),g(p_2), g(p_3), g(p_4), g(q_1), g(q_2)) = (2,3,1,3,1,1)$ and then defining $g$ on $V(H) - \{p_1,p_2,p_3,p_4,q_1,q_2\}$ by the induction hypothesis.
If $f(p'_2) = f(p'_3)$, say $f(p'_2)=f(p'_3)=1$, then without loss of generality, we may assume $(f(x),f(y)) \neq (2,2)$ (since we can swap colors $2$ and $3$), and hence we can extend $f$ to a proper $H$-rainbow-free 3-coloring $g$ of $G$ by first defining $(g(p_1),g(p_2), g(p_3), g(p_4), g(q_1), g(q_2)) = (3,2,3,1,2,2)$ and then defining $g$ on $V(H) - \{p_1,p_2,p_3,p_4,q_1,q_2\}$ by the induction hypothesis.
This completes the proof.
\end{pf}

\begin{lemma} \label{extend to R_0}
Let $G$ be a $\{K_3,R_i: 1 \leq i \leq 7\} \cup \L'$-free subcubic graph and $H$ a subgraph of $G$.
Let $f$ be a proper $(G-H)$-rainbow-free 3-coloring of $G-H$.
If $H$ is isomorphic to $R_0$, then $f$ can be extended to a proper $G$-rainbow-free 3-coloring $g$ of $G$.
\end{lemma}

\begin{pf}
By Lemma \ref{intersecting L_0}, it is sufficient to show that $f$ can be extend to a proper 3-coloring $g$ of $G$ such that $H$ contains no rainbow copy of $L_0$ with respect to $g$.
Let $V(H) = \{v_i, u_j, w: 1 \leq i \leq 4, j=1,2\}$ such that $(v_1,v_2,v_3,v_4,u_1,u_2)$ forms a copy of $L_0$, and $w$ is adjacent to $u_1$ and $v_3$.
Note that if $h$ is a proper 3-coloring of $H$ such that $\{h(v_1), h(v_2), h(v_3)\} = \{1,2,3\}$, then $H$ contains no rainbow copy of $L_0$ with respect to $h$.
Let $N_G(H) = \{x,y,z\}$ such that $xv_2, wy,u_2z \in E(G)$.
Without loss of generality, we may assume that $f(z)=1$.
Then we extend $f$ to $g$ by first defining $g(v_1)=1, g(v_2) \in \{2,3\}-\{f(x)\}, g(v_3)=\{2,3\}-\{g(v_2)\}$, and then it is easy to further define $g$ on the remaining vertices such that $g$ is a proper $G$-rainbow-free $3$-coloring of $G$.
\end{pf}

\begin{lemma}  \label{avoid rainbow}
Let $G$ be a $K_4$-free subcubic graph and $S$ a subset of $E(G)$.
If $G - S$ is $\{K_3, R_i: 1 \leq i \leq 7\} \cup \L'$-free, and every subgraph $H$ of $G-S$ belonging to $\L_t$ for some $t \geq 0$ is an induced subgraph in $G$, then there is a proper $(G-S)$-rainbow-free 3-coloring $f$ of $G$.
\end{lemma}

\begin{pf}
We shall prove this lemma by induction on $\lvert V(G) \rvert$.
If $G-S$ contains no copy of $L_0$ as a subgraph, then a proper $(G-S)$-rainbow-free 3-coloring exists by Brooks' Theorem.
This proves the base case, and we assume that $G-S$ contains a copy of $L_0$ as a subgraph.
Let $H$ be a maximal subgraph of $G-S$ induced by a set $T$ which is of the form $\bigcup_{i=1}^kH_i$ for some positive integer $k$, where each $H_i$ is isomorphic to $L_0$ and $H_j \cap (\bigcup_{i=1}^{j-1}H_i) \neq \emptyset$ for every $1 \leq j \leq k$.
That is, if $T'$ is a subset of $V(G)$ with $(G-S)[T']$ isomorphic to $L_0$, then $T' \cap T = \emptyset$ or $T' \subseteq T$.
By Lemma \ref{intersecting L_0}, $H$ is isomorphic to $R_0$ or a graph in $\L_t$ for some nonnegative integer $t$.
Note that if $H$ is isomorphic to a graph in $\L_t$, then $H$ is an induced subgraph of $G$. 

If $H$ is isomorphic to $R_0$ and $H$ is an induced subgraph of $G$, or $H$ is in $\L_t$ such that $N_G(H)$ is not an independent set, or $\deg_G(v)=2$ for some $v \in V(H)$, then there exists a proper $(G-S)$-rainbow-free 3-coloring of $G$ by extending a proper $(G'-S')$-rainbow-free 3-coloring of $G-V(H)$ that is obtained from applying induction to the graph $G'=G - V(H)$ and the set $S'=S \cap E(G-V(H))$, by Lemmas \ref{extend to R_0} and \ref{extend coloring to L_0's}.
Similarly, if $H$ is isomorphic to $R_0$ but not an induced subgraph of $G$, then $\lvert N_G(H) \rvert \leq 1$, and it is easy to extend a coloring obtained from applying induction on $G-V(H)$ and $S \cap E(G-V(H))$ to a proper $(G-S)$-rainbow-free $3$-coloring.
So we may assume that $H$ is in $\L_t$ such that $N_G(H)$ is an independent set and $\deg_G(v)=3$ for every $v \in V(H)$.

Assume that $\lvert N_G(H) \rvert \leq 3$, then some vertex $x$ in $N_G(H)$ has degree at most one in $G-V(H)$.
Let $G'' = G-V(H)-\{x\}$ and $S'' = S \cap E(G'')$.
Applying induction to $G''$ and $S''$, we obtain a proper $(G''-S'')$-rainbow-free $3$-coloring $f$ of $G''$.
Then we can extend $f$ to a proper $(G-S)$-rainbow-free $3$-coloring of $G$ by first assigning $f(x)$ a color which is different from $f(y)$ for some $y \in N_G(H)-\{x\}$, and then coloring the remaining vertices by Lemma \ref{extend coloring to L_0's}.
Hence, we assume that $\lvert N_G(H) \rvert =4$.

If not all vertices in $N_G(H)$ are contained in a subgraph of $G-V(H)$ isomorphic to $\L_t$ for some $t \geq 0$, then there is a pair of vertices $u,v$ in $N_G(H)$ such that the following hold.
Let $G'''=G - V(H)+uv$ and $S'''=(S \cap E(G-V(H))) \cup \{uv\}$.
Then $G'''$ is $K_4$-free and subcubic, and $G'''-S'''$ is $\{K_3, R_i: 1 \leq i \leq 7\} \cup \L'$-free, and every subgraph $H'$ of $G'''-S'''$ belonging to $\L_t$ for some $t \geq 0$ is an induced subgraph in $G'''$.
In this case, we are done by extending a proper $(G'''-S''')$-rainbow-free $3$-coloring obtained from applying induction to $G'''$ and $S'''$.

Therefore, we assume that all vertices in $N_G(H)$ are in a subgraph $J$ of $G-H$ isomorphic to $\L_t$ for some $t \geq 0$.
But $N_G(H)$ is an independent set, so $J=L_r^{(2)}$ for some odd number $r$.
Note that there is a $4$-cycle $C$ such that $J-C$ is isomorphic to $L^{(2)}_{r-1}$ or $L_0$.
Let $a,c$ be the two diagonal vertices in $C$ such that $a,c \in N_G(H)$.
Denote $(G-V(H)-(V(J)-V(C)))+ac$ by $G^*$ and denote $(S \cap E(G-V(H)-(V(J)-V(C)))) \cup \{ac\}$ by $S^*$.
Apply induction to $G^*$ and $S^*$ to obtain a proper $(G^*-S^*)$-rainbow-free $3$-coloring $f^*$ of $G^*$.
By Lemma \ref{extend coloring to L_0's}, since some vertex in $J-C$ has degree at most two in $G-V(H)$, $f^*$ can be extended to a $(G-(V(H) \cup V(C) \cup (S \cap E(G-V(H)))))$-rainbow-free $3$-coloring $g^*$ of $G-V(H)$ such that $g^*(a) \neq g^*(c)$.
Note that $g^*(a) \neq g^*(c)$ implies that $g^*$ is $(G-V(H)-S^*)$-rainbow-free.
Finally, since $g^*(a) \neq g^*(c)$, $g^*$ can be extended to a proper $(G-S)$-rainbow-free $3$-coloring of $G$ by Lemma \ref{extend coloring to L_0's}.
\end{pf}

\bigskip

The following lemma shows the existence of good colorings for $\{K_3, R_i: 1 \leq i \leq 7\} \cup \L'$-free subcubic graphs.
Hence, Lemma \ref{good 3-coloring} follows from Lemmas \ref{basic structure}, \ref{R_i-free}, \ref{L'-free}, and the following lemma.

\begin{lemma} \label{stronger form for good 3-coloring}
Every $2$-connected $\{K_3, R_i: 1 \leq i \leq 7\} \cup \L'$-free subcubic graph $G$ has a good $3$-coloring.
\end{lemma}

\begin{pf}
Given any proper 3-coloring $f$ of $G$, we define $N(f)$ to be the number of vertices of degree three whose three neighbors have the same color.
Applying Lemma \ref{avoid rainbow} by choosing $S=\emptyset$, there exists a proper $G$-rainbow-free $3$-coloring of $G$.
Let $f$ be a proper $G$-rainbow-free $3$-coloring of $G$ such that $N(f)$ is as small as possible.
We shall show that $f$ is a good $3$-coloring, and it is sufficient to show that every pair of color classes induce a good graph.
Given integers $i,j$ such that $1 \leq i < j \leq 3$, we define $G_{i,j}$ to be the subgraph of $G$ induced by $f^{(-1)}(i) \cup f^{(-1)}(j)$.

\noindent {\bf Claim 1:} For every $1 \leq i <j \leq 3$, $G_{i,j}$ satisfies (G1).

\noindent {\bf Proof of Claim 1:}
Suppose to the contrary, and assume that $G_{2,3}$ does not satisfy (G1).
Then there is a vertex $x$ such that $f(x)=2$, $f(z)=3$ and $\deg_{G_{2,3}}(z) = 3$ for all $z \in N_{G_{2,3}}[x]$.
However, replacing $f(x)$ by 1 reduces $N(f)$, and the coloring is still proper $G$-rainbow-free, a contradiction.
$\Box$

\noindent {\bf Claim 2:} For every $1 \leq i <j \leq 3$, $G_{i,j}$ satisfies (G2).

\noindent {\bf Proof of Claim 2:}
Suppose to the contrary, and assume that $G_{2,3}$ does not satisfy (G2).
Then let $x$, $y$ be two adjacent vertices such that $\deg_{G_{2,3}}(x) = \deg_{G_{2,3}}(y) \allowbreak = 3$ and $\deg_{G_{2,3}}(z) \geq 2$ for all $z \in N_{G_{2,3}}(x) \cup N_{G_{2,3}}(y)$, and $f(x)=2$.
Let $x_1$, $x_2$ be the two neighbors of $x$ in $G$ other than $y$, and let $y_1$, $y_2$ be the two neighbors of $y$ in $G$ other than $x$.
Let $x_{i,j}$ and $y_{i,j}$ be the neighbors of $x_i$ and $y_i$ other than $x$ and $y$, respectively, for $i=1,2$ and $j=1,2$.
Note that if $x_i$ (and $y_i$, respectively) has degree at most two in $G$, then we assume that $x_{i,2}$ (and $y_{i,2}$, respectively) does not exist, and we just ignore it when we mention it in the rest of the proof.
Let $f_x$ (and $f_y$, respectively) be the coloring that is obtained from $f$ by changing $f(x)$ (and $f(y)$, respectively) to 1.
Notice that $f_x$ and $f_y$ are proper 3-colorings with $N(f_x)<N(f), N(f_y) < N(f)$ by the degree condition of neighbors of $x$ and $y$, so there exists a rainbow copy of $L_0$ with respect to each of $f_x$ and $f_y$.
Furthermore, the rainbow copy of $L_0$ with respect to $f_x$ (or $f_y$, respectively) contains $x$ (or $y$, respectively).

First, we assume that $x_{1,1}=x_{2,1}$.
Then by the symmetry between $x_1$ and $x_2$, we may assume that the rainbow copy of $L_0$ with respect to $f_x$ is formed by $(y, x, x_1, x_{1,2}, y_1, y_2)$.
However, no $f_y$-rainbow copy of $L_0$ exists in this case.
Hence $N_{G_{2,3}}[x_1] \cap N_{G_{2,3}}[x_2] = \{x\}$ and $N_{G_{2,3}}[y_1] \cap N_{G_{2,3}}[y_2]=\{y\}$ by symmetry.

Second, we assume that $y_1 = x_{1,1}$ and $y_2 = x_{2,1}$.
By the minimality of $N(f)$ and the triangle-freeness, $(x_{i,2},x_i,x,x_{3-i}, x_{3-i,1}, x_{3-i,2})$ forms a rainbow copy of $L_0$ with respect to $f_x$ for some $i=1,2$.
Then it is easy to check that there is no rainbow copy of $L_0$ with respect to $f_y$.
Similarly, if $y_1=x_{1,1}$, $y_2=x_{1,2}$, then either $G$ has a cut-edge or $N(f)>N(f_x)$ and there is no $f_x$-rainbow copy of $L_0$, a contradiction.
Hence, $\lvert \{x_{1,1}, x_{1,2}, x_{2,1}, x_{2,2}\} \cap \{y_1,y_2\} \rvert \leq 1$ by symmetry.
Similarly, $\lvert \{y_{1,1},y_{1,2},y_{2,1},y_{2,2}\} \cap \{x_1,x_2\} \rvert \leq 1$.

Third, we assume that $\lvert \{x_{1,1},x_{1,2},x_{2,1},x_{2,2}\} \cap \{y_1,y_2\} \rvert = 1$, say $x_{2,1} = y_2$ and $y_{2,1}=x_2$.
By symmetry, every $f_x$-rainbow copy of $L_0$ is formed by $(x_{1,1}, x_1, x, x_2, x_{2,1}, x_{2,2})$, $(x_1, x,x_2, x_{2,2}, x_{1,1}, x_{1,2})$ or $(x_{1,1}, x_1,x,y,y_1, \allowbreak y_2)$ by the triangle-freeness.
If $(x_{1,1}, x_1, x, x_2, x_{2,1}, x_{2,2})$ forms an $f_x$-rainbow copy of $L_0$, then no $f_y$-rainbow copy of $L_0$ exists.
If $(x_1, x,x_2, x_{2,2}, x_{1,1}, x_{1,2})$ forms a rainbow copy of $L_0$ with respect to $f_x$, then $x_{1,2} \neq y_1$ by the $R_1$-freeness of $G$, and the unique rainbow copy of $L_0$ with respect to $f_y$ is formed by $(y_{2,2}, y_2, y, y_1, y_{1,1}, y_{1,2})$, where $y_{2,2}$ is the neighbor of $y_2$ other than $x_2$ and $y$.
Let $f_{y_1, y_2, y_{2,2}}$ be the 3-coloring obtained from $f$ by redefining $(f(y_1), f(y_2), f(y_{2,2}))$ to be $(1,1,2)$.
Then $f_{y_1,y_2,y_{2,2}}$ is a proper 3-coloring with $N(f_{y_1,y_2,y_{2,2}}) < N(f)$ without any rainbow copy of $L_0$.
If $(x_{1,1}, x_1,x,y,y_1,y_2)$ forms an $f_x$-rainbow copy of $L_0$, then no $f_y$-rainbow copy of $L_0$ exists.
Consequently, $\{x_{i,j}: i=1,2, j=1,2\}$ is disjoint from $\{y_1,y_2\}$.
Similarly, $\{x_1,x_2\} \cap \{y_{i,j}: i=1,2,j=1,2\} = \emptyset$.

Therefore, since $f_x$-rainbow copy of $L_0$ exists, then by symmetry, we may assume that any $f_x$-rainbow $L_0$ is formed by $(x_{1,1}, x_1, x, x_2, x_{2,1}, x_{2,2})$ or $(x_{1,1}, x_1, x, y, y_1,y_2)$.
For the latter case, no $f_y$-rainbow copy of $L_0$ exists.
So $(x_{1,1}, x_1, x, x_2, x_{2,1}, x_{2,2})$ forms an $f_x$-rainbow copy of $L_0$, then we may assume that the $f_y$-rainbow copy of $L_0$ is formed by $(y_{1,1}, y_1, y, y_2, y_{2,1}, \allowbreak y_{2,2})$ by symmetry.
Let $f_{y_1,y_2,y_{1,1}}$ be the $3$-coloring obtained from $f$ by changing $(f(y_1), f(y_2), f(y_{1,1}))$ to $(1,1,2)$.
Then $f_{y_1,y_2,y_{1,1}}$ is a proper 3-coloring of $G$ with $N(f_{y_1,y_2,y_{1,1}}) \leq N(f)$, but no $f_{y_1,y_2,y_{1,1}}$-rainbow copy of $L_0$ exists.
So $N(f_{y_1,y_2,y_{1,1}}) = N(f)$ and $f(y_{1,2,1}) = f(y_{1,2,2})=1$, where $y_{1,2,1}$ and $y_{1,2,2}$ are the two neighbors of $y_{1,2}$ other than $y_1$.
In other words, $\deg_{G_{2,3}}(y_1)=2$ and $\deg_{G_{2,3}}(y_{1,2})=1$.
Similarly, we can define $f_{x_1, x_2, x_{1,1}}$ by changing $(f(x_1), f(x_2), f(x_{1,1}))$ to $(1,1,3)$, and then we obtain that $\deg_{G_{2,3}}(x_1)=2$ and $\deg_{G_{2,3}}(x_{1,2})=1$.
So $G_{2,3}$ satisfies (G2), a contradiction.
$\Box$

\noindent {\bf Claim 3:} For every $1 \leq i <j \leq 3$, $G_{i,j}$ satisfies (G3).

\noindent {\bf Proof of Claim 3:}
Suppose to the contrary, and assume that $G_{2,3}$ does not satisfy (G3).
Let $x,y,z \in V(G_{2,3})$, $x_1 \in N_{G_{2,3}}(x) \setminus \{y\}$, $y' \in N_{G_{2,3}}(y) \setminus \{x,z\}$, $y'_1 \in N_{G_{2,3}}(y')-\{y\}$, $z_1 \in N_{G_{2,3}}(z) \setminus \{y\}$ such that $xyz$ is a path in $G_{2,3}$, $\deg_{G_{2,3}}(x)= \deg_{G_{2,3}}(y)= \deg_{G_{2,3}}(z)=3$, $\deg_{G_{2,3}}(x_1)=\deg_{G_{2,3}}(z_1)=1$, $\deg_{G_{2,3}}(y')=2$, and $\deg_{G_{2,3}}(y'_1) \geq 2$, and $f(y)=2$.
Let $r_y$ be the 3-coloring obtained from $f$ by changing $f(y)$ to $1$.
Since $N(r_y)<N(f)$, there is a rainbow copy of $L_0$ with respect to $r_y$.
Without loss of generality, we may assume that the $r_y$-rainbow $L_0$ is formed by $(z,y,y',z_{1,1},z_1,z_2)$, where $N_{G_{2,3}}(z) = \{z_1,z_2,y\}$ and $z_{1,1}$ is a neighbor of $z_1$ distinct from $z$.
Let $r_{z,y',z_{1,1}}$ be the 3-coloring obtained from $f$ by changing $(f(z),f(y'),f(z_{1,1}))$ to $(1,1,3)$.
We derive a contradiction since no $r_{z,y',z_{1,1}}$-rainbow copy of $L_0$ exists and $N(r_{z,y',z_{1,1}})<N(f)$.
$\Box$

\noindent {\bf Claim 4:} For every $1 \leq i <j \leq 3$, $G_{i,j}$ satisfies (G4).

\noindent {\bf Proof of Claim 4:}
Suppose to the contrary, and assume that $G_{2,3}$ does not satisfy (G4).
Let $x,y,z \in V(G_{2,3})$ such that $xyz$ is a path in $G_{2,3}$, $\deg_{G_{2,3}}(x)= \deg_{G_{2,3}}(y)= \deg_{G_{2,3}}(z)=3$, $\deg_{G_{2,3}}(z_1) \geq 2$, and $\deg_{G_{2,3}}(z_2) \geq 2$, where $N_G(z) = \{z_1,z_2,y\}$, and $f(z)=2$, but none of $z_1$, $z_2$ is a support vertex of degree two in $G_{2,3}$.
By Claim 2, $y'$ is not adjacent to $z_1$ and $z_2$, where $y'$ is the neighbor of $y$ other than $x$ and $z$.
Let $p_z$ be the $3$-coloring obtained from $f$ by changing $f(z)$ to $1$.
Since $\deg_{G_{2,3}}(z_1) \geq 2$ and $\deg_{G_{2,3}}(z_2) \geq 2$, we have that $N(p_z) < N(f)$, so there is a rainbow copy of $L_0$ with respect to $p_z$.
By symmetry between $z_1$ and $z_2$, we may assume that the rainbow copy of $L_0$ is formed by $(z_1, z, z_2, z_{2,1}, z_{1,1}, z_{1,2})$, where $N_{G_{2,3}}(z_1)=\{z, z_{1,1}, z_{1,2}\}$ and $N_{G}(z_2) = \{z, z_{2,1}, z_{2,2}\}$.
Notice that it implies that $\deg_{G_{2,3}}(z_1)=3$ and $\deg_{G_{2,3}}(z_2)=2$.
Let $p_{z_1, z_2, z_{2,1}}$ be the $3$-coloring obtained from $f$ by changing $(f(z_1), f(z_2), f(z_{2,1}))$ to $(1,1,3)$.
Since there is no rainbow $L_0$ with respect to $p_{z_1, z_2,z_{2,1}}$, we have $N(p_{z_1, z_2, z_{2,1}})=N(f)$, and this implies that $\deg_{G_{2,3}}(z_{2,2})=1$.
That is, $z_2$ is a support vertex of degree two in $G_{2,3}$, a contradiction.
$\Box$

Together with Claims 1 through 4, $f$ is a good $3$-coloring.
This completes the proof of Lemma \ref{good 3-coloring}.
\end{pf}

\section{Avoiding coloring}
In this section, we investigate sufficient conditions that allow us to extend a coloring of a subgraph to the whole graph, and these lemmas will be used in Section 6.

Recall that given $F:V(G) \rightarrow 2^{[14]}$, we say that $f: V(G) \rightarrow 2^{[14]}$ is an {\it $F$-avoiding coloring} of a graph $G$ if $f(v)$ is disjoint from $F(v) \cup f(u)$ for every pair of adjacent vertices $u$ and $v$.

\begin{lemma}  \label{odd path}
Let $0 \leq r_1 \leq 2$ and $0 \leq r_{2k+1} \leq 2$.
Let $P$ be a path $v_1v_2...v_{2k+1}$ on an odd number of vertices, where $k \geq 1$, and let $F:V(P) \rightarrow 2^{[14]}$ be a function.
If $\lvert F(v_j) \rvert \leq 3+r_j$ for $j=1, 2k+1$ and $\lvert F(v_i) \rvert =2$ for all $2 \leq i \leq 2k$, $\lvert F(v_1) \cap F(v_2) \rvert \leq r_1$, and $\lvert F(v_{2k}) \cap F(v_{2k+1}) \rvert \leq r_{2k+1}$, then there is an $F$-avoiding coloring $f$ of $P$ such that $\lvert f(v_j) \rvert=6$ for every $2 \leq j \leq 2k$, and $\lvert f(v_j) \rvert = 8-r_j$ for $j=1,2k+1$.
\end{lemma}

\begin{pf}
We shall prove this lemma by induction on $k$.
Note that we may assume that $\lvert F(v_j) \rvert = 3+r_j$ for $j=1,2k+1$.
Let $S_1$ be a subset of $F(v_1)-F(v_2)$ of size $3$, and $S_{2k+1}$ be a subset of $F(v_{2k+1})-F(v_{2k})$ of size $3$.
When $k=1$, define $f:V(P) \rightarrow 2^{[14]}$ by letting $f(v_2)=S_1 \cup S_{2k+1} \cup S$, where $S$ is a subset of $[14]-(S_1 \cup S_{2k+1} \cup F(v_2))$ of size $6- \lvert S_1 \cup S_{2k+1}\rvert$, and letting $f(v_1)$ and $f(v_{2k+1})$ be a subset of $[14] - (F(v_1) \cup f(v_2))$ of size $8-r_1$ and a subset of $[14] - (F(v_{2k+1}) \cup f(v_{2k}))$ of size $8-r_{2k+1}$, respectively.
Notice that $\lvert F(v_1) \cup f(v_2) \rvert \leq 6+r_1$ and $\lvert F(v_{2k+1}) \cup f(v_{2k}) \rvert \leq 6+r_{2k+1}$, so $f(v_1)$ and $f(v_{2k+1})$ are well-defined.
This proves the base case.

Now, we assume that $2k+1 \geq 5$ and the lemma holds for the path on $2k-1$ vertices.
Define $F': \{v_i: 3 \leq i \leq 2k+1\} \rightarrow 2^{[14]}$ by $F'(v_3) = F(v_3) \cup S_1$ and $F'(v_i) = F(v_i)$ for all $3<i \leq 2k+1$.
Hence, by the induction hypothesis, there is an $F'$-avoiding coloring $f'$ of $P -\{v_1,v_2\}$ such that $\lvert f'(v_{2k+1}) \rvert = 8-r_{2k+1}$ and $\lvert f'(v_i) \rvert =6$ for all $3 \leq i \leq 2k+1$.
Then define $f:V(P) \rightarrow 2^{[14]}$ by letting $f(v_i)=f'(v_i)$ for all $3 \leq i \leq 2k+1$, $f(v_2) = S_1 \cup T_2$, where $T_2$ is a subset of $[14] - (F(v_2) \cup f'(v_3) \cup S_1)$ of size $3$, and letting $f(v_1)$ be a subset of $[14] - (F_1 \cup f(v_2))$ of size $8-r_1$.
Notice that $f$ is an $F$-avoiding coloring of $P$, so this proves the lemma.
\end{pf}

\begin{lemma}  \label{even path}
Let $P$ be a path $v_1v_2...v_{2k}$ on an even number of vertices, where $k \geq 1$.
Let $F:V(P) \rightarrow 2^{[14]}$ be a function such that $\lvert F(v_j) \rvert \leq 4$ for $j=1, 2k$ and $\lvert F(v_i) \rvert \leq 2$ for all $2 \leq i \leq 2k-1$.
Then there is an $F$-avoiding coloring $f$ of $P$ such that $\lvert f(v_j) \rvert = 8-\deg_P(v_j)$ for $1 \leq j \leq 2k$ if one of the following holds:
	\begin{enumerate}
		\item $k=1$ and $F(v_1) \cap F(v_2) = \emptyset$;

		\item $k=2$, $F(v_1) \cap F(v_2)= \emptyset$, $F(v_3) \cap F(v_4)=\emptyset$, and $\lvert F(v_1) \cap F(v_4) \lvert \leq 2$;

		\item $k \geq 3$, $F(v_{2i+4}) \subseteq F(v_{2i+1})$ for $0 \leq i \leq  k-3$, $F(v_{2k-3}) \subseteq F(v_{2k})$, $(F(v_1)-F(v_4)) \cap (F(v_{2k})-F(v_{2k-3})) = \emptyset$, and $F(v_{2j-1}) \cap F(v_{2j}) = \emptyset$ for $1 \leq j \leq k$;
		\item $k \geq 3$, $F(v_1) \cap F(v_2) =\emptyset$, $F(v_1) \cap F(v_4) = \emptyset$, $F(v_3) \cap F(v_4) = \emptyset$, and $F(v_{2k-1}) \cap F(v_{2k}) = \emptyset$;

		\item $k \geq 3$, $F(v_{2i+4}) \subseteq F(v_{2i+1})$ for $0 \leq i \leq  k-3$, $F(v_{2j-1}) \cap F(v_{2j}) = \emptyset$ for $1 \leq j \leq k-2$, $F(v_{2k-3}) \cap F(v_{2k}) = \emptyset$, and $F(v_{2k-1}) \cap F(v_{2k}) = \emptyset$;

		\item $k \geq 4$, $F(v_1) \cap F(v_2)=\emptyset$, and there exists $0 \leq s \leq k-4$ such that $F(v_{2i+4}) \subseteq F(v_{2i+1})$ for $0 \leq i \leq s$, $F(v_{2i+3}) \cap F(v_{2i+4})=\emptyset$ for $0 \leq i \leq s-1$, $F(v_{2s+3}) \cap F(v_{2s+6})=\emptyset$, $F(v_{2s+5}) \cap F(v_{2s+6})=\emptyset$, $(F(v_1)-F(v_4)) \cap F(v_{2s+6}) = \emptyset$, and $F(v_{2k-1}) \cap F(v_{2k}) = \emptyset$;

		\item $k \geq 2$, $\lvert F(v_1) \rvert \leq 2$, $\lvert F(v_1) \cap F(v_2) \rvert \leq 1$, and $F(v_{2k-1}) \cap F(v_{2k}) = \emptyset$;

		\item $k \geq 3$, $F(v_3) \subseteq F(v_1)$, $F(v_1) \cap F(v_2)=\emptyset$, and $F(v_{2k-1}) \cap F(v_{2k}) = \emptyset$;

		\item $k=4$, $F(v_1)=F(v_4) \cup F(v_6)$, $F(v_1) \cap F(v_2)=\emptyset$, $(F(v_5) \cup F(v_7)) \cap F(v_8) = \emptyset$;

		\item $k \geq 5$, $F(v_1)=F(v_4) \cup F(v_6)$, and $F(v_1) \cap F(v_2) = (F(v_3) \cup F(v_5) \cup F(v_7)) \cap F(v_8) = F(v_{2k-1}) \cap F(v_{2k})= \emptyset$;
		\item $k \geq 4$, $F(v_1)=F(v_4) \cup F(v_6)$, and there exists $0 \leq t \leq k-4$ such that $F(v_{2i+5}) = F(v_{2i+8})$ for $0 \leq i \leq t-1$,
			\begin{enumerate}
				\item $F(v_{2i+5}) \subseteq F(v_{2i+8})$ if $t=k-4$,
				\item $F(v_{2i+5}) = F(v_{2i+8})$ if $t \neq k-4$,
			\end{enumerate}
$F(v_1) \cap F(v_2) = F(v_5) \cap F(v_6) = F(v_{2k-1}) \cap F(v_{2k})= \emptyset$, $F(v_{2j+7}) \cap F(v_{2j+8}) = \emptyset$ for $0 \leq j \leq t-1$, and one of the following holds:
			\begin{enumerate}
				\item $t=k-4$, and $F(v_3) \cap (F(v_{2k})-F(v_{2k-3})) = \emptyset$.
				\item $k \geq 5$, $t = k-5$, $F(v_{2k-3}) \cap F(v_{2k}) = \emptyset$.
				\item $k \geq 6$, $t \leq k-6$, and $(F(v_3) \cup F(v_{2t+7}) \cup F(v_{2t+9})) \cap F(v_{2t+10}) = \emptyset$.
			\end{enumerate}

		\item $k \geq 2$, $\lvert F(v_i) \rvert \leq 4 - 2 \deg_P(v_i)$ for some $1 \leq i \leq 2k$, $i \neq 2,2k-1$, and $F(v_1) \cap F(v_2) = F(v_{2k-1}) \cap F(v_{2k}) = \emptyset$.
	\end{enumerate}
\end{lemma}

\begin{pf}
	\begin{enumerate}
		\item Observe that, since $F(v_1) \cap F(v_2) = \emptyset$, there is a partition $X_1 \cup X_2$ of $[14]$ such that $\lvert X_1 \rvert =\lvert X_2 \rvert=7$, $X_1 \cap F(v_1)= \emptyset$, and $X_2 \cap F(v_2)= \emptyset$.
In that case, $f(v_i)=X_i$ will do.

		\item We may assume that $\lvert F(v_1) \rvert = \lvert F(v_4) \rvert =4$.
Since $\lvert F(v_1) \cap F(v_4) \lvert \leq 2$, there are two disjoint sets $X$ and $Y$, each of size $3$, such that $X \subseteq F(v_1)$ and $Y \subseteq F(v_4)$.
Pick a subset $Y'$ of $[14] - (X \cup Y \cup F(v_3))$ of size $3$, and then assign $f(v_3)=Y \cup Y'$.
Define $f(v_2)$ to be a subset of $[14] - (f(v_3) \cup F(v_2))$ containing $X$ of size $6$.
Finally, pick a subset $W$ of $[14]-(f(v_2) \cup F(v_1))$ of size $7$.
Notice that $\lvert [14]-(f(v_2) \cup F(v_1)) \rvert \geq 7$ since $f(v_2) \cap F(v_1)$ contains $X$ which has size $3$.
Similarly, there is a subset $Z$ of $[14]-(f(v_3) \cup F(v_4))$ of size $7$.
Then assigning $f(v_1) = W$ and $f(v_4)=Z$ gives a desired $F$-avoiding coloring.

		\item We shall do induction on $k$.
Without loss of generality, we may assume that $\lvert F(v_1) \rvert=4$, $\lvert F(v_3) \rvert = \lvert F(v_4) \rvert=2$, and if $k \geq 4$, then $\lvert F(v_7) \rvert =2$.
Let $S$ be a subset of $F(v_1)$ containing $F(v_1)-F(v_4)$ of size 3.
Define $F': \{v_i: 3 \leq i \leq 2k\} \rightarrow 2^{[14]}$ by letting $F'(v_3) = F(v_3) \cup (F(v_1)-F(v_4))$, and $F'(v_i) = F(v_i)$ for all $i > 3$.
Then there is an $F'$-avoiding coloring $f'$ of the subpath $v_3v_4...v_{2k}$ by Statement 2 (when $k=3$) and by the induction hypothesis (when $k >3$).
Note that when $k > 3$, $F(v_3)=F(v_7)$ as $\lvert F(v_3) \rvert = \lvert F(v_7) \rvert=2$, so $F'(v_3)-F'(v_7) \subseteq F(v_1)-F(v_4)$ is disjoint from $F'(v_{2k})-F'(v_{2k-3})$, so the induction applies.
Define $f(v_i)=f'(v_i)$ for $i \geq 4$, and define $f(v_3)$ to be a subset of $f'(v_3)-S$ of size $6$, $f(v_2)$ to be a subset of $[14]-(F(v_2) \cup f(v_3))$ containing $S$ of size $6$, and $f(v_1)$ to be a subset of $[14]-(f(v_2) \cup F(v_1))$ of size $7$.
Note that $\lvert f'(v_3) \cap S \rvert \leq 1$, so $f$ is well-defined.
		
		\item We may assume that $\lvert F(v_1) \rvert = \lvert F(v_{2k}) \rvert =4$ and $\lvert F(v_4) \rvert=2$.
Let $S$ be a subset of $[14]-(F(v_1) \cup F(v_3) \cup F(v_4))$ of size $3$.
Define $F': \{v_i: 4 \leq i \leq 2k\} \rightarrow 2^{[14]}$ by $F'(v_4) = F(v_4) \cup S$, and $F'(v_i) = F(v_i)$ for all $i > 4$.
By Lemma \ref{odd path}, there is an $F'$-avoiding coloring $f'$ with $\lvert f'(v_{2k})\rvert =7$ and $\lvert f'(v_i) \rvert =6$ for all $4 \leq i <2k$.
Pick a subset $S'$ of $[14] - (f'(v_4) \cup F(v_3) \cup F(v_4) \cup S)$ of size $1$.
Assign $f(v_3) = F(v_4) \cup S \cup S'$ and $f(v_i) = f'(v_i)$ for all $4 \leq i \leq 2k$.
Let $T$ be a subset of $[14] - (F(v_2) \cup f(v_3) \cup (F(v_1)-S'))$ with size $6-\lvert F(v_1) - S' \rvert$, and then assign $f(v_2) = (F(v_1) - S') \cup T$.
Finally, assign $f(v_1)$ by any subset of $[14]-(F(v_1) \cup f(v_2))$ with size $7$.
It is clear that $f$ is a desired $F$-avoiding coloring.

		\item Without loss of generality, we may assume that $\lvert F(v_1) \rvert=4$ and $\lvert F(v_3) \rvert = \lvert F(v_4) \rvert =2$.
We shall do induction on $k$.
Assume that $k=3$.
If $F(v_3)=F(v_4)$, then define $S=F(v_4) \cup \{s\}$, where $s \in F(v_1)-F(v_4)$, and $S'=\{s\}$; otherwise, let $S = (F(v_1)-F(v_4)) \cup \{s\}$, where $s \in F(v_4)$ such that $\lvert S \cap F(v_3) \rvert$ is as large as possible, and let $S'$ be a subset of $(S \cup F(v_3))-F(v_4)$ of size $\lvert S \cup F(v_3) \rvert - 2$ such that $\lvert S' \cap F(v_3) \rvert$ is as large as possible.
In fact, $\lvert S' \cap F(v_3) \rvert \geq \lvert S' \rvert-1$.
Then there is a subset $X$ of $F(v_6)$ of size $\lvert F(v_6) \rvert -1$ such that $X \cap S'=\emptyset$.
Define an $F$-avoiding coloring $f: V(G) \rightarrow 2^{[14]}$ by assigning $f(v_5)$ any subset of $[14]-(F(v_5) \cup S')$ containing $X$ of size $6$, $f(v_6)$ any subset of $[14]-(f(v_5) \cup F(v_6))$ of size $7$, $f(v_4)$ any subset of $[14]-(F(v_4) \cup f(v_5))$ containing $S'$ of size $6$, $f(v_3)$ any subset of $[14]-(F(v_3) \cup f(v_4) \cup S)$ of size $6$, $f(v_2)$ any subset of $[14]-(F(v_2) \cup f(v_3))$ containing $S$ of size $6$, and $f(v_1)$ any subset of $[14]-(F(v_1) \cup f(v_2))$ of size $7$.
This prove the base case of the induction.

When $k >3$, let $S_2$ be a subset of $F(v_1)$ of size $3$ such that $\lvert S_2 \cap F(v_4) \rvert \leq 1$.
Define $F':\{v_i: 3 \leq i \leq 2k\} \rightarrow 2^{[14]}$ by letting $F'(v_3)$ be a subset of $(F(v_3) \cup S_2)-F(v_4)$ of size $\lvert F(v_3) \cup S_2 \rvert-1$, and $F'(v_i)=F(v_i)$ for $i \geq 4$.
By induction, there is an $F'$-avoiding coloring $f'$ of the subpath $v_3v_4...v_{2k}$ such that $\lvert f'(v_3) \rvert=\lvert f'(v_{2k}) \rvert=7$ and $\lvert f'(v_i) \rvert=6$ for $4 \leq i \leq 2k-1$.
Then define $f(v_i)=f'(v_i)$ for $i \geq 4$, $f(v_3)$ to be a subset of $f'(v_3)-S_2$ of size $6$, and $f(v_2)$ to be a superset of $S_2$ disjoint from $f(v_3) \cup F(v_2)$ of size $6$, and $f(v_1)$ to be a subset of $[14]-(f(v_2) \cup F(v_1))$ of size $7$.
It is easy to check that $f$ is desired.

		\item Induction on $k$.
Without loss of generality, we assume that $\lvert F(v_1) \rvert=4$, $\lvert F(v_3) \rvert = \lvert F(v_4) \rvert=\lvert F(v_5) \rvert = \lvert F(v_6) \rvert=2$.
Let us first assume that $s=0$.
Let $S_5$ be a superset of $F(v_6)$ disjoint from $(F(v_1)-F(v_4)) \cup F(v_3) \cup F(v_5)$ with size $5$.
By Lemma \ref{odd path}, there is an $F$-avoiding coloring $f$ of the subpath $v_6v_7...v_{2k}$ such that $f(v_6) \cap S_5 = \emptyset$ and $\lvert f(v_i) \rvert = 8-\deg_P(v_i)$ for $6 \leq i \leq 2k$.
Now, we extend $f$ to $V(P)$.
First, define $f(v_5)$ to be a superset of $S_5$ disjoint from $F(v_5) \cup f(v_6)$ of size $6$.
If $F(v_3) \neq F(v_4)$, then let $S_2$ be a subset of $F(v_1)$ of size $3$ such that $\lvert (S_2 \cup F(v_3)) \cap F(v_4) \rvert \leq 1$, and let $S_4$ be a subset of $S_2 \cup F(v_3) - (F(v_4) \cup f(v_5))$ of size $\lvert S_2 \cup F(v_3) \rvert-2$.
If $F(v_3)=F(v_4)$, then let $S_2$ be a subset of $F(v_1)$ containing $F(v_4)$ of size $3$ such that $(S_2-F(v_4)) \cap f(v_5) =\emptyset$, and let $S_4=S_2-F(v_4)$.
Define $f(v_4)$ to be a superset of $S_4$ disjoint from $F(v_4) \cup f(v_5)$ of size $6$, $f(v_3)$ to be a subset of $[14]-(f(v_4) \cup F(v_3))$ of size $6$, $f(v_2)$ to be a superset of $S_2$ disjoint from $F(v_2) \cup f(v_3)$ of size $6$, and $f(v_1)$ to be a subset of $[14]-(F(v_1) \cup f(v_2))$ of size $7$.
Therefore, the claim holds when $s=0$, and thus in particular in the base case $k=4$.

Now assume that the statement holds for every smaller $k \geq 4$, and that $s>0$.
In this case, $F(v_3)=F(v_6)$ since $F(v_6) \subseteq F(v_3)$ and they have the same size.
Let $S$ be a subset of $F(v_1)$ containing $F(v_1)-F(v_4)$ of size $3$.
Define $F': \{v_i: 3 \leq i \leq 2k\} \rightarrow 2^{[14]}$ by letting $F'(v_3)$ be $(F(v_1)-F(v_4)) \cup F(v_3)$, and $F'(v)=F(v)$ for every other $v$.
Note that $F'(v_3)-F'(v_6) = (F(v_1)-F(v_4)) \cup F(v_3)-F(v_6) \subseteq F(v_1)-F(v_4)$, so $(F'(v_3)-F'(v_6)) \cap F'(v_{2s+6}) = \emptyset$.
By induction, there is an $F'$-avoiding coloring $f'$ of the subpath $v_3v_4...v_{2k}$.
Extend $f'$ to $f$ by defining $f(v_i)=f'(v_i)$ for $i \geq 4$, $f(v_3)$ to be a subset of $f'(v_3)-S$ of size $6$, and $f(v_2)$ to be a superset of $S$ disjoint from $F(v_2) \cup f(v_3)$ of size $6$, and $f(v_1)$ to be a subset of $[14]-F(v_1)-f(v_2)$ of size $7$.

		\item Induction on $k$.
We may assume that $\lvert F(v_1) \rvert=2$.
The case that $k=2$ is easy.
When $k \geq 3$, let $S$ be a subset of $F(v_1)-F(v_2)$ of size $1$.
Define $F':\{v_i: 3 \leq i \leq 2k\} \rightarrow 2^{[14]}$ by letting $F'(v_3)$ be a subset of $S \cup F(v_3)$ of size $\lvert S \cup F(v_3) \rvert -1$ such that $\lvert F'(v_3) \cap F(v_4) \rvert \leq 1$, and define $F'(v_i)=F(v_i)$ for $i \geq 4$.
By induction, there is an $F'$-avoiding coloring $f'$ of the subpath.
Extend $f'$ to $f$ by defining $f(v_i)=f'(v_i)$ for $4 \leq i \leq 2k$, $f(v_3)$ to be a subset of $f'(v_3)-(S \cup F(v_3))$ of size $6$, $f(v_2)$ to be a superset of $S$ disjoint from $F(v_2) \cap f(v_3)$ of size $6$, and define $f(v_1)$ to be a subset of $[14]-(f(v_2) \cup F(v_1))$ of size $7$.

		\item We may assume that $\lvert F(v_1) \rvert \geq 3$, otherwise, we are done by Statement 7.
Let $S$ be a subset of $F(v_1)$ containing $F(v_3)$ of size $3$.
Define $F':\{v_i: 3 \leq i \leq 2k\} \rightarrow 2^{[14]}$ by letting $F'(v_3)$ be a subset of $S$ of size $\lvert S \rvert-1$ such that $\lvert F'(v_3) \cap F(v_4) \rvert \leq 1$, and define $F'(v_i)=F(v_i)$ for $i \geq 4$.
By statement 7, there is an $F'$-avoiding coloring $f'$ of the subpath.
And it is easy to use the same argument as in the proof of statement 7 to extend $f'$ to an $F$-avoiding coloring $f$ of $P$.

		\item We may assume that $\lvert F(v_1) \rvert \geq 3$ by Statement 7.
Let $S_2$ be a subset of $F(v_1)$ containing $F(v_6)$ of size $\lvert F(v_1) \rvert -1$, $S_4$ a subset of $(S_2 \cup F(v_3))- F(v_4)$ of size $\lvert S_2 \cup F(v_3) \rvert-2$ such that $\lvert S_4 \cap F(v_6) \rvert$ is as large as possible, $S_6$ a subset of $(S_4 \cup F(v_5))-F(v_6)$ of size $\lvert S_4 \cup F(v_5) \rvert-2$ such that $\lvert S_6 \cap F(v_5) \rvert$ is as large as possible.
It is not hard to see that $\lvert S_6 \cap F(v_8) \rvert \leq 1$, so there is a subset $S_7$ of $F(v_8)- S_6$ of size $\lvert F(v_8) \rvert-1$.
Therefore, we can define a desired $F$-avoiding coloring $f$ of $P$ by defining $f(v_7)$ to be a subset of $[14]-S_6$ of size $6$ containing $S_7$, and then defining $f(v_8)$ to be a subset of $[14]-(F(v_8) \cup f(v_7))$ of size $7$, and $f(v_j)$ from $j=6$ down to $1$ to be a set such that $f(v_j)$ is disjoint from $F(v_j) \cup f(v_{j+1})$ and contains $S_j$ whenever $S_j$ is defined.

		\item Let $S_7$ be a subset of $[14]-(F(v_3) \cup F(v_5) \cup F(v_7))$ containing $F(v_8)$ with size $5$.
By Lemma \ref{odd path}, there is an $F'$-avoiding coloring $f'$ of the subpath $v_8v_9...v_{2k}$, where $F'(v_8)=S_7$ and $F'(v_i)=F(v_i)$ for $9 \leq i \leq 2k$, such that $\lvert f'(v_i) \rvert =6$ for $8 \leq i \leq 2k-1$ and $\lvert f'(v_{2k}) \rvert =7$.
Define $f(v_i)=f'(v_i)$ for $8 \leq i \leq 2k$, and $f(v_7)$ to be a superset of $S_7$ of size $6$ disjoint from $F(v_7) \cup f(v_8)$.
Let $S_5$ be the subset of $(f(v_7) \cup F(v_6))-F(v_5)$ of size $\lvert f(v_7) \cup F(v_6) \rvert-2$ such that $\lvert S_5 \cap S_7 \rvert$ is as large as possible (so $\lvert S_5 \cap F(v_6) \rvert \leq 1$), and $S_3$ be the subset of $(S_5 \cup F(v_4))-F(v_3)$ of size $\lvert S_5 \cup F(v_4)\rvert-2$ such that $\lvert S_3 \cap F(v_1) \rvert = \lvert S_3 \cap (F(v_4) \cup F(v_6)) \rvert \leq 1$.
Define $f(v_2)$ to be a subset of $[14]-(F(v_2) \cup S_3)$ of size $6$ such that $\lvert f(v_2) \cap F(v_1) \rvert \geq \lvert F(v_1) \rvert-1$.
And then define $f(v_1)$ to be a subset of $[14]-(F(v_1) \cup f(v_2))$ of size $7$, and define$f(v_j)$ from $j=3$ through $6$ such that $f(v_j)$ is disjoint from $F(v_j) \cup f(v_{j-1}) \cup S_{j-1} \cup S_{j+1}$ and contains $S_j$ whenever $S_j$ is defined.

	\item Define $S_2$ to be a subset of $F(v_1)$ of size $\lvert F(v_1) \rvert -1$ such that $\lvert S_2 \cap F(v_3) \rvert$ is as large as possible, and define $S_4$ to be a subset of $(S_2 \cup F(v_3)) - F(v_4)$ with size $\lvert S_2 \cup F(v_3) \rvert-2$ such that $\lvert S_4 \cap F(v_3) \rvert$ is as large as possible.
Note that $\lvert S_4 \cap F(v_6) \rvert \leq 1$.
Define $F': \{v_i: 5 \leq i \leq 2k\} \rightarrow 2^{[14]}$ by letting $F'(v_5)$ be a subset of $(S_4 \cup F(v_5)) - F(v_6)$ containing $F(v_5)$ with size $\lvert S_4 \cup F(v_5) \rvert -1$, and $F'(v_i)=F(v_i)$ for $6 \leq i \leq 2k$.
We claim that there is an $F'$-avoiding coloring $f'$ of the subpath $v_5v_6...v_{2k}$ such that $\lvert f'(v_5) \rvert = \lvert f'(v_{2k}) \rvert = 7$ and $\lvert f'(v_i) \rvert = 6$ for $6 \leq i \leq 2k-1$.
If $t=k-4$ and $F(v_3) \cap (F(v_{2k})-F(v_{2k-3})) = \emptyset$, then the claim follows from Statements 2 and 3.
If $k \geq 5$, $t=k-5$, $F(v_{2k-3}) \cap F(v_{2k})=\emptyset$, then the claim follows from Statement 5.
Similarly, if $k \geq 6$, $t \leq k-6$, and $(F(v_3) \cup F(v_{2t+7}) \cup F(v_{2t+9})) \cap F(v_{2t+10}) = \emptyset$, then the claim follows from Statement 6.
Define an $F$-avoiding coloring $f$ of $P$ by letting $f'(v_i)=f(v_i)$ for $6 \leq i \leq 2k$, and assigning $f(v_5)$ a subset of $f'(v_5) - (S_4 \cup F(v_5))$ with size $6$, and from $j=4$ down to $1$, define $f(v_j)$ to be a subset of $[14]-(F(v_j) \cup f(v_{j+1}))$ containing $S_j$ (if $S_j$ is defined) of size $6$ (size $7$ when $j=1$).
It is clear that $f$ is a desired $F$-avoiding coloring of $P$.

	\item Without loss of generality, we may assume that $i$ is odd, so $i \leq 2k-3$.
Define $S_2$ to be a subset of $F(v_1)$ of size $\lvert F(v_1) \rvert -1$, and define $S_{2j}$ to be a subset of $(S_{2j-2} \cup F(v_{2j-1})) - F(v_{2j})$ with size $\lvert S_{2j-2} \cup F(v_{2j-1}) \rvert-2$ for $2 \leq j \leq k-1$.
Since $i \leq 2k-3$, $\lvert S_{2k-2} \rvert \leq 1$.
Define $f(v_{2k-1})$ to be a subset of $[14]-(S_{2k-2} \cup F(v_{2k-1}))$ with size $6$ such that $\lvert f(v_{2k-1}) \cap F(v_{2k}) \rvert \geq \lvert F(v_{2k}) \rvert-1$, and then we define $f(v_i)$ for $i \neq 2k-1$ containing $S_i$ (if $S_i$ is defined) to make $f$ a desired $F$-avoiding coloring.
	\end{enumerate}
This completes the proof.
\end{pf}

\bigskip

Recall that given even numbers $a,b \geq 4$, $H_{a,b}$ is the graph obtained from two disjoint paths $P_a$ and $P_b$ by adding an edge incident with one support vertex in $P_a$ and one support vertex in $P_b$.

\begin{lemma} \label{H_{i,j}}
Let $a,b \geq 4$ be two even integers.
Denote $V(H_{a,b}) = \{v_i, u_j: 1 \leq i \leq a, 1 \leq j \leq b\}$ and $E(H_{a,b}) = \{v_iv_{i+1}, u_ju_{j+1}, v_2u_2: 1 \leq i \leq a-1, 1 \leq j \leq b-1\}$.
Let $F:V(H_{a,b}) \rightarrow 2^{[14]}$ be a function such that $\lvert F(v) \rvert \leq 6-2\deg_{H_{i,j}}(v)$ for every vertex $v$ of $H_{i,j}$.
If $F(u_{b-1}) \cap F(u_b)=\emptyset$, and any condition in Lemma \ref{even path} holds for the path $v_1v_2...v_a$ or the path $v'_1v'_2...v'_a$, where $v'_i = v_{a+1-i}$, then there is an $F$-avoiding coloring $f$ such that $\lvert f(v) \rvert = 8-\deg_{H_{a,b}}(v)-1_{\{u_2\}}(v)$.
\end{lemma}

\begin{pf}
We may assume that $F(u_1) \neq \emptyset$.
Let $S_2$ be a subset of $F(u_1)$ with size $1$.
Apply Lemma \ref{even path}, there is an $F$-avoiding coloring $f'$ of the path $v_1v_2...v_a$ such that $\lvert f'(v_2) \rvert=6$ and $\lvert f'(v_i) \rvert = 8-\deg_{H_{a,b}}(v_i)$ for every $i \neq 2$.
Define $f(v_2)$ to be a subset of $f'(v_2)-S_2$ of size $5$, and $f(v_i)=f'(v_i)$ for $i \neq 2$.
Let $S_3$ be a subset of $f(v_2)-F(u_3)$ of size $1$.
For every $j$ with $4 \leq j \leq b-2$, define $S_j$ to be a subset of $(S_{j-2} \cup F(u_{j-1})) - F(u_j)$ of size $\lvert S_{j-2} \cup F(u_{j-1}) \rvert -2$.
Define $f(u_{b-1})$ to be a superset of $S_{b-3} \cup F(u_{b-2}) \cup F(v_b)$ of size $6$ disjoint with $S_{b-2} \cup F(v_{b-1})$ such that $\lvert f(v_{b-1}) \cap (S_{b-3} \cup F(u_{b-2})) \rvert \geq \lvert S_{b-3} \cup F(u_{b-2}) \rvert - 2$ and $\lvert f(v_{b-1}) \cap F(v_b) \rvert \geq \lvert F(v_b) \rvert-1$.
Then, for $j=b-2$ down to $3$, it is easy to define $f(u_j)$ to be a superset of $S_j$ disjoint from $F(u_j) \cup f(u_{j+1}) \cup S_{j-1}$ such that $\lvert f(u_j) \rvert=6$, and define $f(u_b)$ to be a subset of $[14]-F(v_b)-f(v_{b-1})$ of size $7$.
Define $f(u_2)$ to be a subset of $[14]-(f(v_2) \cup f(u_3))$ containing $S_2$ of size $4$.
Note that such set exists since $f(v_2) \cap f(u_3) \neq \emptyset$, which implies that $\lvert f(v_2) \cup f(u_3) \rvert \leq 10$.
Finally, define $f(u_1)$ to be a subset of $[14]-f(u_2)$ of size $7$.
So it is clear that $f$ is an $F$-avoiding coloring such that $\lvert f(v) \rvert = 8-\deg_{H_{a,b}}(v)-1_{\{u_2\}}(v)$ for every vertex $v$.
\end{pf}

\bigskip

The following lemma is a restatement of Lemma 2.1 in \cite{g}.
Note that Lemma 2.1 in \cite{g} involves a notion called ``kernel", and 
Richardson \cite{r} proved that every digraph that does not contain an odd directed cycle has a nonempty kernel.

\begin{lemma}  \label{choosable}
Let $D$ be a digraph, and let $r_D: V(D) \rightarrow {\mathbb N} \cup \{0\}$ be a function.
For each $v \in V(D)$, let $S(v)$ be a set of size at least $r_D(v)+\sum_{u \in N_D^+(v)}r_D(u)$.
If $D$ contains no odd directed cycle, then there exist subsets $C(v) \subseteq S(v)$ of size $r_D(v)$ for all $v \in V(D)$ such that $C(u) \cap C(v)=\emptyset$ for every pair of adjacent vertices $u,v$ of $D$.
\end{lemma}

Recall that, by Lemma \ref{good 3-coloring}, every fractionally $t$-critical triangle-free subcubic graph, where $t \geq 8/3$, has a proper $3$-coloring of $G$ such that every pair of color classes induces a good graph.

\begin{lemma}  \label{keep good}
Let $G$ be a good subcubic graph.
If $X \subseteq V(G)$, then $G-X$ is good.
\end{lemma}

\begin{pf}
It is clear that $G-X$ satisfies (G1), (G2) and (G4) since $G$ satisfies them.
Furthermore, $G-X$ satisfies (G3) since $G$ satisfies (G2).
\end{pf}

\bigskip

Given a matching $M$ (not necessary induced) of $G$ such that $M$ does not saturate any support vertex, we say that a function $F:V(G) \rightarrow 2^{[14]}$ {\it obeys} $M$ if the following hold:
\begin{enumerate}
	\item $\lvert F(v) \rvert \leq 6-2\deg_G(v)$ for every vertex $v$ in $G$,
	\item $F(v)$ is disjoint from $F(u)$ for every leaf $v$ in $G$ and the neighbor $u$ of $v$,
	\item $\lvert F(x) \cap F(y) \rvert \leq 1$ for every $xy \in M$.
\end{enumerate}

We say that $G$ is {\it $(M,I)$-tractable} if $M$ is a matching of $G$ saturating no support vertices, and $I$ is an independent set of $G$ such that every vertex in $I$ has degree three in $G$, and for every $F$ that obeys $M$, there exist two $F$-avoiding colorings $f_1$ and $f_2: V(G) \rightarrow 2^{[14]}$ of $G$, where $f_1,f_2$ might not be distinct, such that $\lvert f_1(v) \rvert + \lvert f_2(v) \rvert = 16 - 2 \deg_G(v) - 2 \cdot 1_I(v)$.
In this case, we say that $(f_1,f_2)$ is an {\it $F$-avoiding pair} of $G$ with respect to $I$.
Note that every path on an odd number of vertices is $(\emptyset, \emptyset)$-tractable by Lemma \ref{odd path}.

\begin{lemma}  \label{delete P_2}
Let $G$ be a connected subcubic graph which is not a path, and let $u$ be a support vertex of degree two and $N_G(u)=\{v,w\}$, where $v$ is a leaf and $w$ is not a support vertex.
Let $G'=G-\{u,v\}$.
If $G'$ is $(M',I')$-tractable for some matching $M'$ of $G'$ and independent set $I'$ of $G'$, then $G$ is $(M,I')$-tractable for some matching $M$ of $G$.
\end{lemma}

\begin{pf}
If $w$ is a leaf in $G'$, or $w$ is saturated by $M'$, then define $x$ to be the support vertex adjacent to $w$ in $G'$ or the other end of the edge in $M'$ incident with $w$.
Observe that the degree of $x$ is at least two.
If $w$ is a leaf in $G'$ and $x$ has degree two, then $x$ is not a support vertex in $G$ since $G$ is not a path.
Define $M=M' \cup \{wx\}$ if $w$ is a leaf in $G'$ and $\deg_G(x)=2$, and define $M=M'$ otherwise.
So $M$ is a matching that saturates no support vertices of $G$.

Let $F$ defined on $V(G)$ obey $M$.
Define $S_u$ to be a subset of $F(v)$ of size $\lvert F(v) \rvert-1$ such that $\lvert (S_u \cup F(w)) \cap F(x) \rvert \leq 1$ (we assume that $F(x)$ is the empty set when $x$ is not defined).
Define $F':V(G') \rightarrow 2^{[14]}$ by letting $F'(w)$ be a subset of $(S_u \cup F(w))-F(x)$ of size $\lvert S_u \cup F(w) \rvert-1$, and $F'(z)=F(z)$ for every vertex $z$ in $G'$ other than $w$.
So $F'$ obeys $M'$, and there is an $F'$-avoiding pair $(f'_1,f'_2)$ of $G'$ with respect to $I'$.
It is easy to define an $F$-avoiding pair $(f_1,f_2)$ of $G$ with respect to $I'$ such that $f_i(w)$ is a subset of $f'_i(w)-S_u$ of size $\lvert f'_i(w) \rvert-1$, $f_i(u)$ contains $S_u$, and $f_i(y)=f'_i(y)$ for every $y \in V(G')-\{w,u\}$ and $i=1,2$.
\end{pf}

\bigskip

Given nonnegative integers $\alpha, \beta$, an {\it $(\alpha,\beta)$-star} is the graph obtained from $K_{1,\alpha+\beta}$ by subdividing $\beta$ edges.
A {\it spider} is a tree that has exactly one vertex, denoted by the {\it central vertex}, of degree at least three.

\begin{lemma}  \label{spider}
If $G$ is a subcubic spider with the central vertex $v$, then $G$ is $(M, \{v\})$-tractable for some matching $M$.
\end{lemma}

\begin{pf}
By Lemma \ref{delete P_2}, it suffices to show that every $(\alpha, \beta)$-star with central vertex $v$, where $\alpha+\beta=3$, is $(\emptyset, \{v\})$-tractable.
Let $F$ obey $\emptyset$.
Let $V(G)=\{v,x_i, y_j, y'_j: 1 \leq i \leq \alpha, 1 \leq j \leq \beta\}$ and $E(G) = \{vx_i, vy_j, y_jy'_j: 1 \leq i \leq \alpha, 1 \leq j \leq \beta\}$. 
We may assume that $\lvert F(x) \rvert = 6-2\deg_G(x)$ for every vertex $x$ in $G$.
If $(\alpha, \beta)=(3,0)$, then let $S_v$ be a subset of $[14]$ of size three such that $S_v$ intersects $F(x_i)$ for $1 \leq i \leq 3$.
If $(\alpha, \beta)=(2,1)$, then let $S_{y_1}$ be a subset of $F(y'_1)$ of size three, $S_v$ be a subset of $[14]-S_{y_1}$ of size two such that $S_v$ intersects $F(x_1)$ and $F(x_2)$.
If $(\alpha, \beta)=(1,2)$, then let $S_v$ be a subset of $F(x_1)$ of size one, and $S_{y_j}$ a subset of $F(y'_j)-S_v$ of size $3$ for $j=1,2$.
If $(\alpha, \beta)=(0,3)$, then let $S_{y_j}$ be a subset of $F(y'_j)$ of size $3$ for $j=1,2,3$.
Then it is easy to define an $F$-avoiding coloring $f$ of $G$ such that for every vertex $x$, $f(x)$ contains $S_x$ (if $S_x$ is defined) and $f(x)$ has size $8-\deg_G(x)-1_{\{v\}}(x)$ by first assigning $f(v)$ and then assigning $f(x_i)$ and $f(y_j)$ for $1 \leq i \leq \alpha$ and $1 \leq j \leq \beta$, and then assigning $f(y'_j)$ for $1 \leq j \leq \beta$.
Hence, $(f,f)$ is an $F$-avoiding pair of $G$ with respect to $\{v\}$.
\end{pf}

\begin{lemma}  \label{delete P_3 P_4}
Let $G$ be a connected subcubic graph which is not a path or a spider.
Assume that $G$ has a cut-edge $v_2w$ such that $G-v_2w$ contains a path $P=v_1v_2 ... v_{\lvert P \rvert}$ as a component with $\lvert P \rvert=3$ or $4$.
Let $G' = G-P$.
If $G'$ is $(M',I')$-tractable for some matching $M'$ of $G'$ and independent set $I'$ of $G'$, then $G$ is $(M,I' \cup \{v_2\})$-tractable for some matching $M$ of $G$.
\end{lemma}

\begin{pf}
When $w$ is a leaf in $G'$, or $w$ is saturated by $M'$, we define $x$ to be the support vertex adjacent to $w$ in $G'$ or the other end of the edge in $M'$ incident with $w$.
Since $G$ is not a path or a spider, and $M'$ does not saturate any support vertex in $G'$, we know that when $x$ has degree two and $w$ is a leaf in $G'$, $x$ is not a leaf and not a support vertex in $G$.
Define $M=M' \cup \{wx\}$ if $w$ is a leaf in $G'$ and $x$ has degree two, otherwise, define $M=M'$.
So $M$ is a matching that does not saturate any support vertex in $G$.

Let $Y$ be the set of leaves adjacent to $w$.
Note that $x \not \in Y$ if $x$ is defined.
Also, $F(w)=\emptyset$ when $\lvert Y \rvert = 1$, since $G$ is not a spider.
Let $F$ defined on $G$ obey $M$.
Without loss of generality, we may assume that $\lvert F(v) \rvert = 6-2\deg(v)$ for every $v \in V(G)$.
In the rest of the proof, let $F(x)$ be the empty set when $x$ is not defined.
If $\lvert P \rvert=4$, then define $S_{v_2}$ to be a subset of $F(v_1)-F(x)$ of size $\lvert F(v_1) \rvert-3$, and $S_{v_3}$ to be a subset of $F(v_4)-S_{v_2}$ of size $\lvert F(v_4) \rvert-1$; if $\lvert P \rvert=3$, then define $S_{v_2}$ to be a smallest subset of $F(v_1) \cup F(v_3)$ such that $\lvert S_{v_2} \cap F(v_1) \rvert \geq \lvert F(v_1) \rvert-3$, $\lvert S_{v_2} \cap F(v_3) \rvert \geq \lvert F(v_3) \rvert-3$, $\lvert (S_{v_2} \cup F(w)) \cap (\bigcup_{y \in Y}F(y) \cup F(x)) \rvert \leq 2$, and $\lvert (S_{v_2} \cup F(w)) \cap (\bigcup_{y \in Y}F(y) \cup F(x)) \rvert = 2$ only when $\lvert Y \rvert=2$.

Define $F'$ on $V(G')$ by assigning $F'(z)=F(z)$ for $z \neq w$, and assigning $F'(w)$ a subset of $S_{v_2} \cup F(w)$ of size $\lvert S_{v_2} \cup F(w) \rvert-1$ such that $\lvert (S_{v_2} \cup F(w)) \cap (\bigcup_{y \in Y}F(y) \cup F(x)) \rvert$ is as small as possible.
If $\lvert Y \rvert <2$, then $F'(w)$ is disjoint from $\bigcup_{y \in Y}F(y) \cup F(x)$, so $F'$ obeys $M'$.
In this case, there is an $F'$-avoiding pair $(f'_1,f'_2)$ of $G'$ with respect to $I'$.
And then it is easy to obtain an $F$-avoiding pair $(f_1,f_2)$ of $G$ with respect to $I' \cup \{v_2\}$ by modifying $f'_1$ and $f'_2$.
Note that $w$ has degree at most two in $G'$, so $w \not \in I'$ and $I' \cup \{v_2\}$ is an independent set in $G$.
If $\lvert Y \rvert=2$, then we define $f:V(G) \rightarrow 2^{[14]}$ by first setting $f(w)$ a subset of $[14]-S_{v_2}$ of size $5$ such that $\lvert F(w) \cap F(y) \rvert \geq 2$ for every $y \in Y$.
And then we extend $f$ to $G$ such that $\lvert f(v) \rvert = 8-\deg_G(v)-1_{\{v_2\}}(v)$ for every $v \in V(G)$.
So $(f,f)$ is an $F$-avoiding pair of $G$.
That is, $G$ is $(\emptyset, \{v_2\})$-tractable.
\end{pf}

\begin{lemma}  \label{super H}
Let $t$ be a nonnegative integer, and let $G$ be the graph obtained from a path $u_1u_2u_3u_4v_1v_2...v_{2t} \allowbreak  u_5u_6u_7u_8$ by attaching one leaf $u'_i$ on each vertex $u_i$ for $3 \leq i \leq 6$.
Then $G$ is $(\emptyset, \{u_3, u_6\})$-tractable.
\end{lemma}

\begin{pf}
Let $F$ obey the empty set, so $F(u_1)$ is disjoint from $F(u_2)$ and $F(u_7)$ is disjoint from $F(u_8)$.
Without loss of generality, we may assume that $\lvert F(v) \rvert = 6-2\deg(v)$ for every $v \in V(G)$.
Define $S_{u_4}$ to be a subset of $F(u'_4)$ of size $\lvert F(u_4) \rvert - 2$, and for each $1 \leq i \leq t$, define $S_{v_{2i}}$ to be a subset of $(S_{v_{2i-2}} \cup F(v_{2i-1})) - F(v_{2i})$ of size $\lvert S_{v_{2i-2}} \cup F(v_{2i-1}) \rvert - 2$, where $v_0$ is $u_4$.
Furthermore, define $S_{u_3}$ to be a subset of $F(u'_3)-S_{u_4}$ of size $\lvert F(u'_3) \rvert-3$, $S_{u_2}$ to be a subset of $F(u_1)-S_{u_3}$ of size $\lvert F(u_1) \rvert -1$, $S_{u_5}$ to be a subset of $F(u'_5) - S_{v_{2t}}$ of size $\lvert F(u'_5) \rvert -2$, $S_{u_6}$ to be a subset of $F(u'_6)-S_{u_5}$ of size $\lvert F(u'_6) \rvert-3$, and $S_{u_7}$ to be a subset of $F(u_8)-S_{u_6}$ of size $\lvert F(v_8) \rvert -1$.
Then it is easy to find an $F$-avoiding coloring $f$ of $G$ by first assigning $f(u_5)$ a superset of $S_{u_5}$ of size $5$, and then assigning $f(v_j)$ a subset of $[14]-(F(v_j) \cup f(v_{j+1}) \cup S_{v_{j-1}} \cup S_{v_{j+1}})$ containing $S_{v_j}$ (let $S_{v_k} = \emptyset$ for undefined $S_{v_k}$, where $-1 \leq k \leq 2t+1$) of size $8-\deg_G(v_j)$ from $j=2t$ down to $0$, where $v_{-1}=u_3$ and $v_{2t+1} = u_5$, and then assigning $f(z)$ by a superset of $S_z$ (if $S_z$ is defined) for other vertices $z$ of size $8-\deg_G(z)-1_{\{u_3,u_6\}}(z)$.
Hence, $(f,f)$ is an $F$-avoiding pair with respect to $\{u_3,u_6\}$.
\end{pf}

\begin{lemma}  \label{support deg 3 adj deg 2 and 2}
Let $G$ be a triangle-free subcubic graph such that no support vertex has degree two.
Let $v$ be a vertex of degree three in $G$, and let $N_G(v) = \{v', u,w\}$, where $v'$ is a leaf and $u,w$ are of degree two.
If $G'=G-\{v,v'\}$ is $(M',I')$-tractable for some matching $M'$ and independent set $I'$, then $G$ is $(M, I' \cup \{v\})$-tractable, for some matching $M$ in $G$.
\end{lemma}

\begin{pf}
Observe that $u$ is not adjacent to $w$ since $G$ is triangle-free.
Let $u',w'$ be the neighbor of $u,w$ distinct from $v$, respectively.
Let $M=M' \cup \{zz': (z,z')=(u,u')$ or $(w,w'), \deg_G(z')=2\}$.
Note that $u$ and $w$ are leaves in $G'$, so $M'$ does not saturate $u,u',w$ or $w'$.
And since no support vertex has degree two in $G$, $M$ is a matching of $G$ saturating no support vertices in $G$.

Let $F$ defined on $V(G)$ obey $M$.
We may assume that $F(u),F(w),F(v')$ are nonempty.
Define $S_v$ to be a subset of $[14]$ of size $4$ such that $S_v$ intersects $F(u)$, $F(w)$ and $F(v')$ such that $\lvert (S_v \cup F(u)) \cap F(u') \rvert \leq 1$ and $\lvert (S_v \cup F(w)) \cap F(w') \rvert \leq 1$.
Define $F'$ on $V(G')$ by assigning $F'(u)$ a subset of $(F(u) \cup S_v)-F(u')$ of size $\lvert F(u) \cup S_v \rvert -1$, and assigning $F'(w)$ a subset of $(F(w) \cup S_v)-F(w')$ of size $\lvert F(w) \cup S_v \rvert-1$, and define $F'(z)=F(z)$ for every other vertex $z$ in $G'$.
Hence $F'$ obeys $M'$, and there is an $F'$-avoiding pair $(f'_1,f'_2)$ of $G'$ with respect to $I'$, and then it is easy to obtain an $F$-avoiding pair $(f_1,f_2)$ of $G$ with respect to $I' \cup \{v\}$ from $f_1',f_2'$ by modifying $f_i'(u),f_i'(w)$, for $i=1,2$.
\end{pf}

\bigskip

Recall that ${\mathcal H}$ is the family of graphs consisting of $H_{a,b}$ for every pair of even integers $a,b \geq 4$.

\begin{lemma} \label{non-path non-H}
Let $G$ be a good graph such that no component of $G$ is a path or a graph in ${\mathcal H}$.
Then $G$ is $(M,I)$-tractable for some matching $M$ and independent set $I$.
\end{lemma}

\begin{pf}
We shall prove this lemma by induction on $\lvert V(G) \rvert$.
Clearly, $G$ is $(\emptyset,\emptyset)$-tractable if $\lvert V(G) \rvert \leq 3$.
Assume the lemma holds for every good proper subgraph of $G$ containing no path or a graph in ${\mathcal H}$ as a component, but $G$ is not $(M,I)$-tractable for any $M$ and $I$.
So $G$ is connected.
And by Lemma \ref{spider}, $G$ is not a spider.
The following claim is an immediate consequence of Lemmas \ref{keep good} and \ref{delete P_2}.

\noindent {\bf Claim 1:} If $v$ is a support vertex of degree two adjacent to a nonsupport vertex, then $G-\{v,v'\}$ is a graph in ${\mathcal H}$, where $v'$ is the leaf adjacent to $v$.

\noindent {\bf Claim 2:} No support vertex of degree two is adjacent to a support vertex of degree three.

\noindent {\bf Proof Claim 2:}
Suppose that there is a support vertex of degree two adjacent to a support vertex of degree three.
Then there is a cut-edge $e$ in $G$ such that one component $C$ of $G-e$ is a path of order $4$.
By Lemmas \ref{delete P_3 P_4}, \ref{keep good} and \ref{odd path}, $G$ is $(M,I)$-tractable for some $M$ and $I$ unless $G-C$ is a path on an even number of vertices or a graph in ${\mathcal H}$.
If $G-C$ is a path on an even number of vertices, then $G$ is either a spider or in ${\mathcal H}$ by Claim 1, a contradiction.
And if $G-C$ is a graph in ${\mathcal H}$, then by Claim 1, either there is another cut-edge $e'$ such that one component $C'$ of $G-e'$ is a path on four vertices but $G-C'$ is not a graph in ${\mathcal H}$, or $G$ is the graph mentioned in Lemma \ref{super H}.
So no support vertex of degree two adjacent to a support vertex of degree three.
$\Box$

\noindent {\bf Claim 3:} No support vertex is of degree two.

\noindent {\bf Proof of Claim 3:}
Let $v$ be a support vertex of degree $2$, and let $N_G(v) = \{v',w\}$, where $v'$ is a leaf.
Since $G$ is a not path, $w$ is not a support vertex by Claim 2.
Hence, every support vertex of degree $2$ is adjacent to a leaf and a nonsupport vertex.
Together with Claim 1, there is no support vertex of degree two.
$\Box$

Similarly, by taking advantage of Lemmas \ref{keep good}, \ref{delete P_3 P_4} and Claim 3, we have the following Claim.

\noindent {\bf Claim 4:} Every support vertex is adjacent to exactly one leaf.

In addition, the following claim follows from the bipartiteness of $G$, Claim 3 and Lemmas \ref{keep good} and \ref{support deg 3 adj deg 2 and 2}.

\noindent {\bf Claim 5:} No support vertex is adjacent to two vertices of degree two.
Consequently, every support vertex is of degree three and adjacent to exactly one leaf and at least one vertex of degree three.

Let $X$ be the subgraph of $G$ induced by the vertices of degree three.
By (G1), $X$ has maximum degree at most two, so $X$ is a disjoint union of paths and cycles.

\noindent {\bf Claim 6:} $X$ is a disjoint union of paths.

\noindent {\bf Proof of Claim 6:}
Suppose that some component of $X$ is a cycle.
Since $G$ is connected and bipartite, $G$ is obtained from an even cycle by attaching a leaf on each vertex, by (G2), (G3) and Claim 3.
We shall obtain a contradiction by showing that $G$ is $(\emptyset, \emptyset)$-tractable.
Denote $V(G)=\{v_i,u_i: 0 \leq i \leq 2k-1\}$ and $E(G)=\{v_iv_{i+1}, v_iu_i: 0 \leq i \leq 2k-1\}$, where the index is computed under modulo $2k$, for some positive integer $k$.
Let $F$ obey the empty set, and we may assume that $\lvert F(u_i) \rvert=4$ for $0 \leq i \leq 2k-1$.
Give an orientation on the edges of the cycle $v_0v_1...v_{2k-1}v_0$ such that every vertex in the cycle has in-degree and out-degree one.
By Lemma \ref{choosable}, for every $0 \leq i \leq 2k-1$, we can pick a subset $S_i$ of $F(u_i)$ with size $2$ such that $S_i$ is disjoint from $S_{i-1} \cup S_{i+1}$.
Again by Lemma \ref{choosable}, for every $0 \leq i \leq 2k-1$, we can pick a subset $T_i$ of $[14]-\bigcup_{j=i-1}^{i+1} S_j$ of size $3$ such that $T_i$ is disjoint from $T_{i-1} \cup T_{i+1}$.
Then it is easy to define an $F$-avoiding coloring $f$ of $G$ such that $\lvert F(z) \rvert = 8-\deg_G(z)$ for every vertex $z$ in $G$ and $f(v_i)=S_i \cup T_i$ for every $0 \leq i \leq 2k-1$.
So $(f,f)$ is an $F$-avoiding pair of $G$, and hence $G$ is $(\emptyset, \emptyset)$-tractable.
$\Box$

The following claim is an immediate consequence of (G2), (G3), (G4), Claims 3,4,5 and 6.

\noindent {\bf Claim 7:} If $P$ be a component in $X$, then $P$ is a path that satisfies the following.
	\begin{enumerate}
		\item If $\lvert V(P) \rvert \geq 3$, then every vertex in $P$ is a support vertex in $G$.
		\item If $\lvert V(P) \rvert = 2$, then $P$ contains a support vertex in $G$.
		\item If $\lvert V(P) \rvert = 1$, then the vertex in $P$ is a nonsupport vertex of degree $3$ in $G$.
	\end{enumerate}

Let $X_1, X_2, ...$ be the components of $X$.
If $X_i$ is a path on two vertices, then let $x_i \in X_i$ be a support vertex in $G$.
Let $I$ be the subset of $V(X)$ consisting of $x_i$, for each component $X_i$ with $\lvert V(X_i) \rvert=2$, and the ends of every other component of $X$.
So $I$ is an independent set in $G$, and every vertex in $I$ has degree three in $G$.
We shall show that $G$ is $(\emptyset, I)$-tractable.
Let $G'$ be the graph obtained from $G$ by deleting all leaves, so every support vertex in $G$ has degree two in $G'$ by Claims 3 and 4.
Note that $G'$ has minimum degree at least two.

\noindent {\bf Claim 8:} There exists an orientation $O_1$ of $G'$ such that every vertex has in-degree and out-degree at least one, and if $v \in V(G')$ with $\deg_{G'}(v)=3$ is a nonsupport vertex in $G$ adjacent to a support vertex in $G$, then $v$ is pointed by a nonsupport vertex in $G$, and $v$ points to a nonsupport vertex in $G$.

\noindent {\bf Proof of Claim 8:}
Let $Y$ be the graph obtained from $G'$ by deleting every edge that is incident with one support vertex in $G$ and one nonsupport vertex in $G$ of degree three.
By Lemma \ref{orientation}, there is an orientation $O$ of $Y$ such that every vertex $v$ has in-degree and out-degree at least $\lfloor \deg_Y(v)/2 \rfloor$.
Let $J = \{u \in V(G'): u$ is a support vertex in $G$, $N_G(u)$ contains a nonsupport vertex of degree three$\}$.
By (G4), $J$ is an independent set in $G$.
Note that if $u \in J$, then $u$ is incident with an edge in $E(G')-E(Y)$.
Define an orientation $O_1$ of $G'$ by assigning the direction of each edge in $E(G') \cap E(Y)$ the same direction as in $O$, and assigning each edge in $E(G')-E(Y)$ a direction such that every vertex in $J$ has in-degree and out-degree one.
Note that $O_1$ exists since $J$ is an independent set in $G$.
If $u \in V(G')$ is a support vertex in $G$, then either $\deg_Y(u)=2$ or $u \in J$, so $u$ has in-degree and out-degree one in $O_1$.
If $u \in V(G')$ is a nonsupport vertex in $G$, then $\deg_Y(u) \geq 2$ by Claim 7, so $u$ has in-degree and out-degree at least one in $O_1$.
Furthermore, if $v \in V(G')$ is a nonsupport vertex in $G$ with $\deg_{G'}(v)=3$ adjacent to a support vertex in $G$, then $v$ is an end of some $X_i$ with $\lvert V(X_i) \rvert=2$, and the other two neighbors of $v$ are of degree two in $G$.
In other words, $\deg_Y(v)=2$, and the two neighbors of $v$ in $Y$ are nonsupport vertices by Claim 3.
Therefore, $v$ is pointed by a nonsupport vertex in $G$, and $v$ points to a nonsupport vertex in $G$.
$\Box$

Let $O_2$ be the orientation of $G'$ obtained from $O_1$ by reversing the direction of each edge.
Define $D_1$ and $D_2$ to be the digraph whose underlining graph is $G'$ equipped with the orientation $O_1$ and $O_2$, respectively.

Let $F$ defined on $V(G)$ obey the empty set, and we may assume that $\lvert F(v) \rvert = 6-2\deg_G(v)$ for every $v \in V(G)$.
Since $X$ is a union of paths, we can write $V(G) = \{v_i: 1 \leq i \leq \lvert V(G) \rvert\}$ such that every vertex in $X$ has at most one neighbor in $G$ with smaller index, and the index of every vertex in $X$ is smaller than the index of any vertex not in $X$, and the index of any non-leaf of $G$ is smaller than the index of any leaf of $G$.
Note that $V(G') = \{v_i: 1 \leq i \leq \lvert V(G') \rvert\}$.
Now, for $i=1$ and $2$, we consecutively define functions $r_i: V(D_i) \rightarrow {\mathbb N}$ and $T_i: \{$support vertices of $G\} \rightarrow 2^{[14]}$ such that for each $v$ in the domain of $T_i$, either $\lvert T_i(v) \rvert \leq 2$, or $\lvert T_i(v) \rvert=5$ and $T_i(w) = \emptyset$ for every $w$ in the domain of $T_i$ and adjacent to $v$, according to the above ordering of $V(G')$ and the following rules:
	\begin{enumerate}
		\item Assume that $v$ is the vertex in some component $X_k$ of $X$ with $\lvert V(X_k) \rvert =1$.
If $\deg_{D_i}^+(v)=1$, then define $r_i(v)=6$; otherwise, define $r_i(v)=2$.
		\item Assume that $\lvert V(X) \rvert \geq 2$, and $v$ is an end of a component of $X$, and $v$ is in $I$.
So $v$ is a support vertex in $G$, and we let $v'$ be the leaf adjacent to $v$.
Define $r_i(v)=3$.
In $D_i$, if $v$ is pointed by a vertex of degree $3$ in $G$, then define $T_i(v)=\emptyset$; otherwise, $v$ is pointed by a vertex of degree $2$ in $G$, and we define $T_i(v)$ to be a subset of $F(v')-T_i(u)$ of size $2$, where $u$ is the neighbor of $v$ in $G$ such that $T_i(u)$ is defined (if any).
		\item Assume that $\lvert V(X) \rvert \geq 2$, and $v$ is an end of a component of $X$, but $v$ is not in $I$.
So $v$ is an end of a component $X_k$ of $X$ with $\lvert V(X_k) \rvert = 2$.
Assume that $v$ is not a support vertex in $G$.
In $D_i$, if $v$ is pointed by a vertex of degree $3$ in $G$, then define $r_i(v)=6$; otherwise, define $r_i(v)=4$.
		\item Assume that $\lvert V(X) \rvert \geq 2$, and $v$ is an end of a component of $X$, but $v$ is not in $I$.
Assume that $v$ is a support vertex in $G$, and let $v'$ be the leaf adjacent to $v$.
In $D_i$, if $v$ is pointed by a vertex of degree $3$ in $G$, then define $r_i(v)=3$ and $T_i(v)=\emptyset$; otherwise, define $r_i(v)=2$, and $T_i(v)$ to be a superset of $F(v')$ of size $5$ such that $T_i(v) \cap F(u) \neq \emptyset$, where $u$ is the neighbor of $v$ of degree $2$ in $G$.

(Notice that when $\lvert T_i(v) \rvert=5$, $v$ points to a vertex $w$ of degree $3$ in $G$.
And $w$ is the unique neighbor of $v$ having degree $3$ in $G$.
Also, $w$ satisfies rule $2$.
Since $w$ is pointed by $v$, $T_i(w) = \emptyset$.)

		\item Assume that $\lvert V(X) \rvert \geq 3$, and $v$ is an internal vertex of a component of $X$.
So $v$ is a support vertex in $G$, and we let $v'$ be the leaf adjacent to $v$.
Define $r_i(v)=3$, and $T_i(v)$ to be a subset of $F(v')-T_i(u)$ with size $2$, where $u$ is the neighbor of $v$ such that $T_i(u)$ is defined (if any).
		\item Define $r_i(v)=6$ for every vertex $v$ with $\deg_G(v)=2$.
	\end{enumerate}

It is clear that $r_1(v)+r_2(v)+ \lvert T_1(v) \rvert + \lvert T_2(v) \rvert = 16-2\deg_G(v)-2 \cdot 1_I(v)$ for every vertex $v$ in $G'$ (i.e. every non-leaf in $G$).
Define $L_i:V(D_i) \rightarrow 2^{[14]}$ by letting $L_i(v) = [14]-(F(v) \cup \bigcup_{u \in N_{G'}[v]} T_i(u))$ for every vertex $v$ in $G'$.
So $\lvert L_i(v) \rvert \geq r_i(v) + \sum_{u \in N^+_{D_i}} r_i(u)$ for every vertex $v$ in $G'$.
By Lemma \ref{choosable}, there are $F$-avoiding colorings $g_1$ and $g_2$ of $G'$ such that $g_i(v) \subseteq L_i(v)$ and $\lvert g_i(v) \rvert \geq r_i(v)$ for $i=1,2$, so $\lvert g_1(v) \rvert + \lvert g_2(v) \rvert \geq r_1(v) + r_2(v)$.
Define $f_i(v)=g_i(v) \cup T_i(v)$ for $i=1,2$ and $v \in V(G')$ (we let $T_i(v)$ be the empty set if $T_i(v)$ is not defined), so $\lvert f_1(v) \rvert + \lvert f_2(v) \rvert \geq r_1(v) + r_2(v) + \lvert T_1(v) \rvert + \lvert T_2(v) \rvert = 16-2\deg_G(v)-2 \cdot 1_I(v)$ for $v \in V(G')$.
In addition, for every pair of a support vertex $v$ and a leaf $v'$ such that $v$ is adjacent to $v'$, since either $T_i(v) \subseteq F(v')$, or $\lvert T_i(v) \rvert=5$ and $r_i(v)=2$ and $F(v') \subseteq T_i(v)$, we have that $\lvert f_i(v)-F(v') \rvert \leq \lvert g_i \rvert+1 \leq r_i(v)+1$, and the equalities hold only when $r_i(v)=2$ and $\lvert T_i(v) \rvert =5$.
So $\lvert f_i(v) - F_i(v') \rvert \leq 3$ for every pair of a support vertex $v$ and a leaf $v'$ adjacent to $v$.
Therefore, we can define $f_i(v')$ to be a subset of $[14]-(f_i(v) \cup F(v'))$ such that $\lvert f_i(v') \rvert=7$.
As a result, $(f_1,f_2)$ is an $F$-avoiding pair of $G$ with respect to $I$.
This completes the proof.
\end{pf}

\section{Penetrations and cooperations}
In this section, we will prove remaining lemmas from Section 3.

\bigskip

\noindent{\bf Proof of Lemma \ref{admit odd path}:}
Let $C = v_1v_2...v_k$, where $k$ is odd.
If $k=1$, then let $B=(N_G(C), \emptyset)$.
Clearly, for every $(B,\emptyset)$-compatible function $F$, $\lvert F(v_1) \rvert \leq 6$, so there exists an $F$-avoiding coloring $f:V(C) \rightarrow 2^{[14]}$ of $C=v_1$ such that $\lvert f(v_1) \rvert = 8$.
That is, $(B,\emptyset)$ penetrates $C$.

If $k \geq 5$, then let $B=(N_G(C), \{u_1u_2, u_{k-1}u_k: u_i \in N_G(v_i) \cap N_G(C)$ for $i=1,2,k-1,k\})$.
Note that $B$ is loopless since $G$ is triangle-free.
Therefore, for every $(B,\emptyset)$-compatible function $F$, $\lvert F(v_1) \rvert \leq 4$, $\lvert F(v_k) \rvert \leq 4$, $\lvert F(v_i) \rvert=2$ for every $2 \leq i \leq k-1$, and $F(v_1) \cap F(v_2) = F(v_{k-1}) \cap F(v_k) = \emptyset$.
By Lemma \ref{odd path}, there exists an $F$-avoiding coloring $f: V(C) \rightarrow 2^{[14]}$ of $C$ such that $\lvert f(v_i) \rvert = 8-\deg_C(v_i)$ for $1 \leq i \leq k$.
In other words, $(B, \emptyset)$ penetrates $C$.

When $k=3$, let $B_1$ and $B_2$ be graphs such that $V(B_1)=V(B_2)=N_G(C)$ and $E(B_1) = \{x_1x_2: x_i \in N_G(v_i) \cap N_G(C), i=1,2\}$, $E(B_2) = \{x_2x_3: x_i \in N_G(v_i) \cap N_G(C), i=2,3\}$.
Then it is easy to check that $(B_1, \emptyset)$ cooperates with $(B_2, \emptyset)$.
$\Box$.

\bigskip

\noindent{\bf Proof of Lemma \ref{admit even path and H}:}
First, we prove this lemma when $C$ is a path on even number of vertices.
Let $C=v_1v_2...v_{2k}$, and let $B_1, B_2$ be graphs with $V(B_1)=V(B_2)=N_G(C)$.
Now we define $E(B_1)$ and $E(B_2)$.
If $k=1$, then define $E(B_1)=E(B_2)=\{xy: x \in N_G(v_1) \cap N_G(C), y \in N_G(v_2) \cap N_G(C)\}$.
By Statement 1 of Lemma \ref{even path}, $(B_1 \emptyset)$ and $(B_2, \emptyset)$ satisfy Statements 1-6 of this lemma.
So we may assume that $k \geq 2$.
If $\deg_G(v_i) \leq 2$ for some $1 \leq i \leq 2k$ and $i \neq 2, 2k-1$, then define $E(B_1)=E(B_2)=\{u_ju_{j+1}: u_j \in N_G(v_j) \cap N_G(C), j=1,2k-1\}$.
Note that $\deg_{B_i}(v) \leq 2n_C(v)$ for every vertex $v$ in $B_i$.
By Statement 12 of Lemma \ref{even path}, both $(B_1, \emptyset)$ and $(B_2, \emptyset)$ penetrate $C$.
So we may assume that $\deg_G(v_i)=3$ for $i \neq 2,2k-1$.

Let $N_G(v_1)-C=\{u_1,u'_1\}$, $N_G(v_{2k})-C=\{u_{2k},u'_{2k}\}$ and $N_G(v_i)-C = \{u_i\}$ for $2 \leq i \leq 2k-1$.
Note that each $u_2$ and $u_{2k-1}$ may not exist, and we just ignore it when we mention it in the remaining of this proof if it does not exist.
If $k=2$, $N_G(v_1) \cap N_G(v_4) = \emptyset$ and $n_C(u_1)=n_C(u_1')=1$, then define $E(B_1)=\{u_1u_2, u'_1u_2, u_3u_4, u_3u'_4, u_1u_4, u_1u'_4\}$ and $E(B_2)=\{u_1u_2, u'_1u_2, u_3u_4, u_3u'_4, u_1'u_4, u'_1u_4'\}$.
If $k=2$, $N_G(v_1) \cap N_G(v_4) = \emptyset$, and some of $u_1$ and $u_1'$ has $n_C$-value at least two, say $u_1$, then $u_1=u_3$ by the triangle-freeness, and we define $E(B_1)=E(B_2)=\{u_1u_2, u_1u_4,u_1u_4',u_1'u_2\}$.
If $k=2$ and $N_G(v_1) \cap N_G(v_4) \neq \emptyset$, then $\lvert N_G(v_1) \cap N_G(v_4) \rvert \leq 1$ since there is no $f$-rainbow copy of $L_0$, and we assume that $u'_1=u'_4$ and define $E(B_1)=E(B_2)=\{u_1u_2, u'_1u_2, u_3u_4, u_3u'_4, u_1u_4\}$.
Then $(B_1,\emptyset)$ and $(B_2,\emptyset)$ satisfy Statements 1-6 of this lemma by Statements 1 and 2 of Lemma \ref{even path}.

Now we assume that $k \geq 3$.
Let $s=-\infty$ if $N_G(v_1) \cap N_G(v_4) = \emptyset$; otherwise, let $s=\max\{i \in {\mathbb Z}: N_G(v_{2j+1}) \cap N_G(v_{2j+4}) \neq \emptyset, 0 \leq j \leq i\}$.
Note that either $s=-\infty$ or $0 \leq s \leq k-2$.
If $s \geq 0$, then let $u_1 \in N_G(v_4)$.
If $s=k-2$, then let $u_{2k} \in N_G(v_{2k-3})$.

\noindent {\bf Case 1:} Assume that $s \geq 0$, $n_C(u_1') \geq 2$, $n_C(u_3) \geq 2$, and either $n_C(u_1')=3$ or $n_C(u_3)=3$.
Let $v \in \{u_1',u_3\}$ such that $n_C(v)=3$.
If $v \in N_G(v_{2k-1})$ and $\{u_1,u_1',u_3\}-\{v\} \subseteq N_G(v_{2k})$, then $n_C(u_1)=3$, and we define $J=\{u_1,u_1',u_3\}-\{v\}$; if $v \not \in N_G(v_{2k-1})$, then define $J=\{v,u\}$, where $u \in \{u_1,u_1',u_3\}-(N_G(v_{2k-1}) \cup \{v\})$; if $v \in N_G(v_{2k-1})$ but $\{u_1,u_1',u_3\}-\{v\} \not \subseteq N_G(v_{2k})$, then define $J=\{v,u\}$, where $u \in \{u_1,u_1',u_3\} - (N_G(v_{2k}) \cup \{v\})$.
Define $E(B_1)=E(B_2)=\{u_1u_2,u_1'u_2,u_{2k-1}u_{2k},u_{2k-1}u_{2k}'\}$.
Then it is straight forward to check that $(B_1,\{J\})$ and $(B_2,\{J\})$ satisfy Statements 1-6 of this lemma by Statements 7 and 8 of Lemma \ref{even path}.

\noindent {\bf Case 2:} Assume that $k-4 \geq s \geq 0$ but $u_1' \not \in N_G(v_{2s+6})$.
Define $E(B_1)=\{u_1u_2, u'_1u_2, u_{2j-1}u_{2j}, u_{2k-1}u_{2k}, u_{2k-1}u'_{2k}, \allowbreak u'_1u_{2s+6}, \allowbreak u_{2s+3}u_{2s+6}, \allowbreak u_{2s+5}u_{2s+6}: \allowbreak 2 \leq j \leq s+1\}$, and define $E(B_2)$ later.
Then $(B_1, \emptyset)$ penetrates $C$ by Statement 6 of Lemma \ref{even path}, and at most one vertex $v$ in $B_1$ has degree more than $2n_C(v)$.
Moreover, if such $v$ exists, then $v$ is adjacent to a vertex $u$ in $B_1$ with $\deg_{B_1}(u) \leq 2n_C(u)-1$, and $v$ is adjacent to a vertex $v_i$ in $C$ with an even index $i=2s+6$; if $\deg_{B_1}(v) \geq 2n_C(v)+1$, then either $n_C(v) \neq 2$, or $v \in N_G(v_2) \cup N_G(v_{2k-1})$.
If there is no vertex $v$ such that $\deg_{B_1}(v)=2n_C(v)+1$, or the vertex $v$ with $\deg_{B_1}(v)=2n_C(v)+1$ satisfies $n_C(v)=3$, then we further define $E(B_2)=E(B_1)$.

\noindent {\bf Case 3:} Assume that $s \geq 0$, $n_C(u_1)=3$, $n_C(u_1')=n_C(u_3)=2$, and either $u_1 \not \in N_G(v_{2k-1})$, or $u_1 \in N_G(v_{2k-1})$ but $\{u_1', u_3\} \not\subseteq N_G(v_{2k})$.
If $u_1 \not \in N_G(v_{2k-1})$, then define $J=\{u_1,v\}$, where $v \in \{u_1',u_3\}-N_G(v_{2k-1})$; otherwise, define $J=\{u_1,v\}$, where $v \in \{u_1',u_3\}-N_G(v_{2k})$.
Define $E(B_1)=E(B_2)=\{u_1u_2,u_1'u_2,u_{2k-1}u_{2k},u_{2k-1}u_{2k}'\}$.
It is straight forward to see that $(B_1,\{J\})$ and $(B_2,\{J\})$ satisfy Statements 1-6 of this lemma by Statements 7 and 8 of Lemma \ref{even path}.

\noindent {\bf Case 4:} Assume that $k-4 \geq s \geq 0$, $n_C(u_1)=3$, $n_C(u_1')=n_C(u_3)=2$, $u_1 \in N_G(v_{2k-1})$, and $\{u_1', u_3\} \subseteq N_G(v_{2k})$.
Observe that $N_G(u_1') \cap C = \{v_1,v_{2k}\}$.
Define $E(B_1)=\{u_1u_2, u'_1u_2, u_{2j-1}u_{2j}, \allowbreak u_{2k-1}u_{2k}, u_{2k-1}u'_{2k}, \allowbreak u'_1u_{2s+6}, \allowbreak u_{2s+3}u_{2s+6}, \allowbreak u_{2s+5}u_{2s+6}: \allowbreak 2 \leq j \leq s+1\}$, and define $E(B_2)$ later.
Then $(B_1, \emptyset)$ penetrates $C$ by Statement 6 of Lemma \ref{even path}, and at most one vertex $v$ in $B_1$ has degree more than $2n_C(v)$.
Moreover, if such $v$ exists, then $v$ is adjacent to a vertex $u$ in $B_1$ with $\deg_{B_1}(u) \leq 2n_C(u)-1$, and $v$ is adjacent to a vertex $v_i$ in $C$ with an even index $i=2s+6$; if $\deg_{B_1}(v) \geq 2n_C(v)+1$, then either $n_C(v) \neq 2$, or $v \in N_G(v_2) \cup N_G(v_{2k-1})$.
If there is no vertex $v$ such that $\deg_{B_1}(v)=2n_C(v)+1$, or the vertex $v$ with $\deg_{B_1}(v)=2n_C(v)+1$ satisfies $n_C(v)=3$, then we further define $E(B_2)=E(B_1)$.

\noindent {\bf Case 5:} Assume that $s \geq 0$, $n_C(u_1)=n_C(u_1')=n_C(u_3)=2$, and $\{u_1',u_3\} \cap N_G(v_{2k-1}) = \emptyset$.
Observe that $u_1 \not \in N_G(v_{2k-1})$ since $s \geq 0$ and $n_C(u_1)=2$.
Define $J=\{u_1,u_1',u_3\}$, and define $E(B_1)=E(B_2)=\{u_1u_2,u_1'u_2,u_{2k-1}u_{2k},u_{2k-1}u_{2k}'\}$.
Then it is straight forward to check that $(B_1,\{J\})$ and $(B_2,\{J\})$ satisfy Statements 1-6 of this lemma by Statements 7 and 8 of Lemma \ref{even path}.

\noindent {\bf Case 6:} Assume that $s=k-2$.
So $N_C(u_3) \geq 2$.
We may assume that Cases 1, 3, and 5 do not hold.
We claim that $u_1' \neq u_{2k}'$.
Suppose that $u_1'=u_{2k}'$, then $n_C(u_1') \geq 2$.
Since Case 1 does not hold, $n_C(u_1')=n_C(u_3)=2$.
So $\{u_1',u_3\} \cap N_G(v_{2k-1}) = \emptyset$.
Since Case 5 does not hold, $n_C(u_1)=3$.
Since Case 3 does not hold, $u_1 \in N_G(v_{2k-1})$ and $\{u_1',u_3\} \subseteq N_G(v_{2k})$, so $k=3$.
But it implies that $u_1v_4v_5$ forms a triangle in $G$, a contradiction.
This proves that $u_1' \neq u_{2k}'$.
Define $E(B_1)=E(B_2)=\{u_1u_2, u'_1u_2, \allowbreak u_{2k-1}u_{2k}, \allowbreak u_{2k-1}u'_{2k}, u'_1u'_{2k}, u_{2j-1}u_{2j}: 2 \leq j \leq k-1\}$.
Then $(B_1, \emptyset)$ and $(B_2, \emptyset)$ satisfy Statements 1-6 of this lemma by Statement 3 of Lemma \ref{even path}.

\noindent {\bf Case 7:} Assume that $s=k-3$.
Define $E(B_1)=E(B_2)=\{u_1u_2, u'_1u_2, \allowbreak u_{2j-1}u_{2j}, \allowbreak u_{2k-1}u_{2k}, u_{2k-1}u'_{2k}, u_{2k-3}u_{2k}, u_{2k-3}u'_{2k}: 2 \leq j \leq k-2\}$.
Then $(B_1, \emptyset)$ and $(B_2, \emptyset)$ satisfy Statements 1-6 of this lemma by Statement 5 of Lemma \ref{even path}.

\noindent {\bf Case 8:} Assume that $s = -\infty$.
Define $E(B_1)=\{u_1u_2,u'_1u_2,u_{2k-1}u_{2k}, \allowbreak u_{2k-1}u'_{2k},u_1u_4,u'_1u_4,u_3u_4\}$, and define $E(B_2)$ later.
Similarly, $(B_1, \emptyset)$ penetrates $C$ by Statement 4 of Lemma \ref{even path}, and at most one vertex $v$ in $B_1$ has degree more than $2n_C(v)$.
Once such $v$ exists, $v$ is adjacent to $v_4$; if $\deg_{B_1}(v)=2n_C(v)+1$ and $n_C(v) \geq 2$, then either $n_C(v)=3$, or $v \in N_G(v_2) \cup N_G(v_{2k-1})$.
If there is no vertex $v$ such that $\deg_{B_1}(v)=2n_C(v)+1$, or the vertex $v$ with $\deg_{B_1}(v)=2n_C(v)+1$ satisfies $n_C(v)=3$, then we further define $E(B_2)=E(B_1)$.

\noindent {\bf Case 9:} Assume that $s=0$, $u_1' \in N_G(v_6)$, and either $n_C(u_3)=1$, or $n_C(u_3)=2$ and $u_3 \in N_G(v_{2k-1})$.
Note that it is the remaining case that $E(B_1)$ has not been defined, since $s \geq 1$ implies that $n_C(u_3) \geq 2$.
We consider the following subcases.
	\begin{itemize}
		\item When $k \geq 5$ (so $u_3 \neq u_8$) and $u_5 \neq u_8$, define $E(B_1)=\{u_1u_2, u'_1u_2, \allowbreak u_3u_8, u_5u_8, u_7u_8, u_{2k-1}u_{2k}, \allowbreak u_{2k-1}u'_{2k}\}$, then $(B_1, \emptyset)$ satisfies Statements 1-6 of this lemma by Statement 10 of Lemma \ref{even path}.
Note that every vertex $v$ other than $u_8$ in $B_1$ has degree at most $2n_C(v)$.
And $\deg_{B_1}(u_8) \geq 2n_C(u_8)+1$ implies that $n_C(u_8) \neq 2$ or $n_C(u_8) = 2$ with $u_8 \in N_G(v_2) \cup N_G(v_{2k-1})$.
We further define $E(B_2)=E(B_1)$ if $\deg_{B_1}(u_8) \leq 2n_C(u_8)$ or $n_C(u_8)=3$.
		\item When $k \geq 4$ and $u_5=u_8$, let $t = \max\{i \in {\mathbb Z}: N(v_{2j+5}) \cap N(v_{2j+8}) \neq \emptyset, 0 \leq j \leq i\}$, so $0 \leq t \leq k-4$.
If $t=k-4$, then we assume that $u_{2k} \in N_G(v_{2k-3})$.
Define $E(B_1)=\{u_1u_2,u'_1u_2,u_5u_6, \allowbreak u_{2i+7}u_{2i+8},u_{2k-1}u_{2k},u_{2k-1}u'_{2k}: 0 \leq i \leq t-1\} \cup T$, where $T$ is defined as follows:
			\begin{itemize}
				\item $T=\{u_3u'_{2k}\}$ if $t=k-4$;
				\item $T = \{u_{2k-3}u_{2k},u_{2k-3}u'_{2k}\}$ if $t=k-5$;
				\item $T=\{u_3u_{2t+10}, u_{2t+7}u_{2t+10}, u_{2t+9}u_{2t+10}\}$ if $0 \leq t \leq k-6$.
			\end{itemize}
Notice that $B_1$ does not have loops since $n_C(u_3) \geq 2$ implies that $N_G(u_3) \cap V(C)=\{v_3,v_{2k-1}\}$.
Hence $(B_1, \emptyset)$ penetrates $C$ by Statement 11 of Lemma \ref{even path}.
Also, every vertex $v$ other than $u_{2t+10}$ in $B_1$ has degree at most $2n_C(v)$, and $\deg_{B_1}(u_{2t+10}) \geq 2n_C(u_{2t+10})+1$ only if $0 \leq t \leq k-6$.
And if $\deg_{B_1}(u_{2t+10}) \geq 2n_C(u_{2t+10})+1$, then $\deg_{B_1}(u_{2t+10})=2n_C(u_{2t+10})+1$, and either $n_C(u_{2t+10}) \neq 2$, or $n_C(u_{2t+10}) = 2$ with $u_{2t+10} \in N_G(v_2) \cup N_G(v_{2k-1})$.
We further define $E(B_2)=E(B_1)$ if $t \geq k-5$ or $\deg_{B_1}(u_{2t+10}) \leq 2n_C(u_{2t+10})$ or $n_C(u_{2t+10})=3$.
		\item When $k=4$ and $u_5 \not \in N_G(v_8)$, define $E(B_1)=E(B_2) = \{u_1u_2,u'_1u_2,\allowbreak u_5u_8,u_5u'_8,u_7u_8,u_7u'_8\}$, then $(B_1, \emptyset)$ and $(B_2, \emptyset)$ satisfy Statements 1-6 of this lemma by Statement 9 of Lemma \ref{even path}.
		\item When $k=3$, define $E(B_1)=E(B_2)=\{u_1u_2,u'_1u_2, u_3u_6, u_3u'_6, \allowbreak u_5u_6,u_5u'_6\}$, then $(B_1, \emptyset)$ and $(B_2, \emptyset)$ satisfy Statements 1-6 of this lemma by Statement 5 of Lemma \ref{even path}.
	\end{itemize}

Consequently, $E(B_1)$ are defined in all cases.
In fact, if $B_2$ is not defined, then $k \geq 3$ and there exist an even number $p \not \in \{2,2k\}$ and a vertex $w_1 \in N(v_p)$ of $B_1$ such that the following hold: 
	\begin{itemize}
		\item $\deg_{B_1}(w_1) = 2n_C(w_1)+1$; 
		\item $n_C(w_1) \leq 2$;
		\item if $n_C(w_1)=2$, then $w_1 \in N_G(v_p) \cap N_G(v_{2k-1})$ or $w_1 \in N_G(v_p) \cap N_G(v_2)$;
		\item if $v \in V(B_1)-\{w_1\}$, then $\deg_{B_1}(v) \leq 2n_C(v)$;
		\item if $v \in V(B_1)-\{w_1\}$ and $\deg_{B_1}(v)=2n_C(v)$, then $N_G(v) \cap V(C) \subseteq \{v_1,v_2,v_{2k-1}\}$.
	\end{itemize}
Let $x_i = v_{2k+1-i}$ for $1 \leq i \leq 2k$.
So $x_1x_2...x_{2k}$ is the same path as $v_1v_2...v_{2k}$, but the order of the indices of vertices is reversed.
Now, for those cases that $E(B_2)$ have not been defined, we define other boundary-graphs $B_1',B_2'$ for the path $x_1x_2...x_{2k}$ with $V(B_1')=V(B_1)=V(B_2')=V(B_2)$ by applying the same rules as we defined $B_1, B_2$.
The same argument shows that either $E(B_1')$ and $E(B_2')$ were defined, or $E(B_2')$ was not defined but $E(B_1')$ was defined in a way such that every but one vertex $w_2 \in N(v_{2k+1-q})$ of $B_1'$ with an even number $q \not \in \{2,2k\}$ has degree at most $2n_C(v)$ such that the following hold: 
	\begin{itemize}
		\item $\deg_{B_1'}(w_2) = 2n_C(w_2)+1$;
		\item $n_C(w_2) \leq 2$;
		\item if $n_C(w_2)=2$, then $w_2 \in N_G(v_{2k+1-q}) \cap N_G(v_2)$ or $w_2 \in \allowbreak N_G(v_{2k+1-q}) \cap N_G(v_{2k-1})$;
		\item if $v \in V(B_1') - \{w_2\}$ and $\deg_{B_1'}(v) = 2n_C(v)$, then $N_G(v) \cap V(C) \subseteq \{v_2,v_{2k-1},v_{2k}\}$.
	\end{itemize}
Finally, we complete the process of defining $B_1$ and $B_2$ in the following way.
If $E(B_1)$ and $E(B_2)$ were defined, then we are done; if $E(B_2)$ was not defined but $E(B_2')$ was defined, then replacing $E(B_1)$ by $E(B_1')$, and defining $E(B_2)=E(B_2')$; otherwise, defining $E(B_2)$ to be $E(B_1')$.
It remains to check that $\deg_{B_{3-i}}(w_i) \leq 2n_C(w_i)-1$ for $i=1,2$.
Suppose to the contrary, and without loss of generality, we may assume that $\deg_{B_2}(w_1)>2n_C(w_1)-1$.
If $w_1 \neq w_2$, then $\deg_{B_2}(w_1)=2n_C(w_1)$, so $v_p \in N_G(w_1) \cap V(C) \subseteq \{v_2,v_{2k-1},v_{2k}\}$, a contradiction.
If $w_1=w_2$, then since $p$ and $q$ are even, $n_C(w_1)=2$, so $w_2=w_1 \in N_G(v_p) \cap (N_G(v_2) \cup N_G(v_{2k-1}))$.
However, $2k+1-q$ is odd but not equal to $2k-1$, a contradiction.
This proves the lemma when $C$ is a path on an even number of vertices.

Now, we assume that $C$ is a graph in ${\mathcal H}$.
Denote $V(C)$ by $\{x_i,y_j: 1 \leq i \leq a, 1 \leq j \leq b\}$ and denote $E(C)$ by $\{x_ix_{i+1},y_jy_{j+1}, x_2y_2: 1 \leq i \leq a-1, 1 \leq j \leq b-1\}$ for some even positive integers $a$ and $b$.
Let $V(B'''_1)=V(B'''_2)=N_G(\{x_i: 1 \leq i \leq a, i \neq 2\}) \cap N_G(C)$.
Define $E(B'''_1)$, $E(B'''_2)$, ${\mathcal J'''_1}$, and ${\mathcal J'''_2}$ as in the case of path on an even number of vertices.
Then we define $B_1,B_2,\J_1,\J_2$.
Define $V(B_1)=V(B_2)=N_G(C)$, and let $S=\{uz:u \in N_G(y_{b-1}) \cap N_G(C), z \in N_G(y_b) \cap N_G(C)\}$.
First, assume that $\J_1'''=\J_2'''=\emptyset$.
In this case, we define $\J_1=\J_2=\emptyset$.
If there exist $1 \leq i \leq 2$ and $w_i \in V(B_i''')$ such that $\deg_{B_i'''}(w_i)=2n_C(w_i)+1$, and either $n_C(w_i)=3$, or there exists $1 \leq j \leq b$ such that $w_i \in N_G(y_j)$ and either $j \neq b-1$, or $j=b-1$ but $N_{B_i'''}(w_i) \cap N_G(y_b) \cap N_G(C) \neq \emptyset$, then defining $E(B_1)=E(B_2)=E(B_i''') \cup S$.
Otherwise, defining $E(B_1)=E(B_1''') \cup S$ and $E(B_2)=E(B_2''') \cup S$.
Now, we assume that $\J_1''' \neq \emptyset \neq \J_2'''$.
Recall that $B_1'''=B_2'''$ and $\J_1'''=\J_2'''$ in this case.
Then we define $E(B_1)=E(B_1''') \cup S$ and $E(B_2)=E(B_2''') \cup S$.
If $\lvert \J_1''' \rvert=2$ or $\J_1''' \cap (N_G(y_{b-1}) \cup N_G(y_b))=\emptyset$, then defining $\J_1=\J_2=\J_1'''$; otherwise, defining $\J_1=\J_2$ to be a subset $I$ of $\J_1'''$ of size two such that $I$ is an independent set in $B_1$ and $B_2$, and $I$ contains a vertex with $n_C$-value three.
It is straight forward to check that $(B_1,\J_1)$ and $(B_2,\J_2)$ satisfy Statements 1-6.
This proves the lemma.
$\Box$

\bigskip

\noindent{\bf Proof of Lemma \ref{admit non-path non-H}:}
Since $f$ is a good $3$-coloring, $C$ is a good graph.
By Lemma \ref{non-path non-H}, there is a matching $M$ of $C$ saturating no support vertex of $C$ and an independent set $I \subseteq \{v: \deg_C(v)=3\}$ such that given $F:V(C) \rightarrow 2^{[14]}$ that obeys $M$, there are two $F$-avoiding colorings $f_1,f_2$ of $C$, such that $\lvert f_1(v) \rvert + \lvert f_2(v) \rvert = 16 - 2\deg_C(v) -2 \cdot 1_I(v)$ for all $v \in V(C)$.
Let $S = \{v'u': v' \in N_G(v) \cap N_G(C), u' \in N_G(u) \cap N_G(C), \deg_C(v)=1, u \in N_C(v)\}$.
Let $T = \{x'y': xy \in M, x' \in N_G(x) \cap N_G(C), y' \in N_G(y) \cap N_G(C)\}$.
Define $V(B)=N_G(C)$ and $E(B)=S \cup T$.
Then by Lemma \ref{non-path non-H}, for every $(B,\emptyset)$-compatible function $F$, there exist an independent set $I$ and two $F$-avoiding colorings $f_1,f_2$ such that $\lvert f_1(v) \rvert + \lvert f_2(v) \rvert = 16 - 2\deg_C(v) -2 \cdot 1_I(v)$ for all $v \in V(C)$.
In other words, $(B, \emptyset)$ penetrates $C$.
$\Box$

\section{Bounded number of colors}
In this section, we prove that every triangle-free subcubic graph has a $(516:180)$-coloring.
To achieve this objective, we need a variation of Lemma \ref{2 cut set} to deal with vertex-cuts of size two that does not induced an edge.

\begin{lemma} \label{stronger form for 2-cut}
Let $k$ be a positive integer.
Let $G_0, G_1, ...G_k$ be induced subgraphs of a graph $G$ such that $G=G_0 \cup G_1 \cup ... \cup G_k$, $\lvert V(G_0) \cap V(G_i) \rvert =2$ and $V(G_i) \cap V(G_j) \subseteq V(G_0)$ for every different integers $i,j$ with $i,j \in \{1, 2, ..., k\}$.
Given $1 \leq i \leq k$, denote $V(G_0 \cap G_i)$ by $\{u_i,v_i\}$.
Assume that $u_iv_i$ is not an edge of $G_i$ for every $1 \leq i \leq k$.
Let $a,b,c,d$ be positive integers such that $a/b \geq c/d$.
If $G_0$ has an $(a:b)$-coloring, and $G_i+u_iv_i$ and $G_i/u_iv_i$ have $(c:d)$-colorings for every $1 \leq i \leq k$, then $G$ has an $(ad,bd)$-coloring.
\end{lemma}

\begin{pf}
Let $f_0$ be an $(a:b)$-coloring of $G_0$, and let $x_i = \lvert f_0(u_i) \cap f_0(v_i) \rvert$ for every $1 \leq i \leq k$.
So $G_0$ has an $(ad:bd)$-coloring $f_0'$ such that $\lvert f_0'(u_i) \cap f_0'(v_i) \rvert = x_id$.
For every $1 \leq i \leq k$, let $g_i$ and $h_i$ be a $(c:d)$-coloring of $G_i/u_iv_i$ and $G_i+u_iv_i$, respectively.
Since $u_iv_i$ is not an edge of $G_i$, there exist $(c:d)$-colorings $g_i',h_i'$ of $G_i$ such that $\lvert g'_i(u_i) \cap g'_i(v_i) \rvert = d$ and $\lvert h'_i(u_i) \cap h'_i(v_i) \rvert=0$.
Define $f_i:V(G_i) \rightarrow 2^{[bd]}$ to be a $(bc:bd)$-coloring of $G_i$ by setting $f_i(v) = \{y+pd,z+qd: y \in g_i(v), z \in h_i(v), 0 \leq p \leq x_i-1, x_i \leq q \leq b-1\}$ for every vertex $v$ of $G_i$.
Observe that $\lvert f_i(u_i) \cap f_i(v_i) \rvert = x_id$.
Therefore, by swapping colors, we may assume that $f_i(v)=f_0'(v)$ for every $v \in \{u_i,v_i\}$.
Define $f$ on $V(G)$ such that $f(v)=f'_0(v)$ for every vertex $v$ of $G_0$ and $f(v)=f_i(v)$ for every vertex $v$ of $G_i$.
Consequently, $f$ is an $(ad:bd)$-coloring of $G$.
\end{pf}

\begin{theorem}
Every triangle-free subcubic graph $G$ has a $(516:180)$-coloring.
\end{theorem}

\begin{pf}
By Theorem \ref{strong main} and Lemma \ref{stronger form for good 3-coloring}, every $2$-connected $\{K_3,R_i:1 \leq i \leq 7\} \cup \L'$-free subcubic graph has an $(172:60)$-coloring.
Observe that it implies that every $\{K_3,R_i:1 \leq i \leq 7\} \cup \L'$-free subcubic graph has an $(172:60)$-coloring.
Recall that $R_1, R_2,...,R_7$ and $\L'$ are defined in Section 4.

Given $i=0,1,2$, we say that a graph $H$ is an {\it $i$-extension} if it is a triangle-free subcubic graph such that a graph in $\{R_i:1 \leq i \leq 7\} \cup \L'$ can be obtained from $H$ by deleting $i$ vertices.
Observe that every graph in $\{R_i:1 \leq i \leq 7\} \cup \L'$ has minimum degree at least two and contains at most two vertices of degree less than three.
Recall that for $i=0,1,2$, every $i$-extension has an $(8:3)$-coloring by Lemmas \ref{extend 8:3} and \ref{L'-free} and the coloring in Figure \ref{R's}.

Let $\F$ be a maximal set of pairwise disjoint induced subgraphs of $G$, where each member of $\F$ belongs to $\{K_3,R_i: 1 \leq i \leq 7\} \cup \L'$.
Note that if $G$ is a member of $\F$, then $G$ is a $0$-extension and we are done.
Let $G_0 = G-\bigcup_{H \in \F} V(H)$.
Then $G_0$ has an $(172:60)$-coloring since it is $\{K_3,R_i: 1 \leq i \leq 7\} \cup \L'$-free.
Let $G'$ be a maximal induced subgraph of $G$ such that it can be repeatedly obtained from $G_0$ and an $i$-extension $H'$ of a member $H$ of $\F$ by identifying $i$ vertices of $G_0$ that induce a clique in $G$ with $V(H')-V(H)$, for some $i=1,2$, and then renaming the resulting graph as $G_0$.
Clearly, $G'$ has an $(172:60)$-coloring.
Then $G$ has a $(516:180)$-coloring by Lemma \ref{stronger form for 2-cut}.
\end{pf}

\section{Concluding remarks}
In this paper, we proved that $\chi_f(G) \leq 43/15$ for every triangle-free subcubic graph $G$.
In fact, we proved the stronger statement that for every fractionally $t$-critical triangle-free subcubic graph $G$ with $t \geq 8/3$, there is an independent set $I$ and a function $f:V(G) \rightarrow 2^{[168]}$ such that $\lvert f(v) \rvert = 72 - 4 \deg_G(v) - 4 \cdot 1_I(v)$ for every vertex $v$.
Hence, $\chi_f(G) \leq 14/5$ if we can get rid of the error term $1_I(v)$ from the above equation.
Recall that every vertex in $I$ is adjacent to three vertices which receive the same color from a good $3$-coloring.
So one possible way to remove the error term is to prove that every fractionally $14/5$-critical graph has a proper $3$-coloring without any rainbow copy of $L_0$ such that every pair of color classes induces a graph of maximum degree at most two.
Note that it is not hard to show that such a $4$-coloring exists by using a result of the author and Yu about linear colorings \cite{ly}.

Instead of looking for a proper coloring such that every pair of color classes induce a graph of maximum degree at most two, it is also interesting to study a proper coloring such that every pair of color classes induces a graph whose $2$-connected components have simple structures.
More precisely, we conjecture that every subcubic graph other than $K_4$ has a proper $3$-coloring such that every pair of color classes induce a disjoint union of cacti.
A cactus is a graph whose blocks are edges and cycles.
At the moment when we tried to simplify the proof in this paper, we thought that the positive answer of the conjecture might be helpful.

Finally, we note that Heckman and Thomas \cite{ht2} also conjectured that every planar triangle-free subcubic graph has fractional chromatic number at most $8/3$.
In Section 4, we investigated the structure of fractionally $t$-critical graphs with $t \geq 8/3$.
So it might be helpful in dealing with this conjecture.

\bigskip

\noindent{\bf Acknowledgement.}
The author thanks Professor Robin Thomas and Albert Bush for many suggestions for preparing this paper.
The author is also very grateful for the referees for careful reading and valuable comments.

\end{document}